\documentclass[pdflatex,sn-mathphys,Numbered]{sn-jnl}

\usepackage{multirow}%
\usepackage{mathrsfs}%
\usepackage[title]{appendix}%
\usepackage{textcomp}%
\usepackage{manyfoot}%
\usepackage{algpseudocode}%
\usepackage{listings}%

\usepackage[utf8]{inputenc}
\usepackage{amsmath}
\usepackage{dutchcal}
\usepackage{amsfonts}
\usepackage{amssymb}
\usepackage{amsthm}
\usepackage[dvipsnames]{xcolor}
\usepackage{graphicx}

\usepackage{tikz}
\usepackage{tikz-3dplot}

\usetikzlibrary{spy}
\usetikzlibrary{patterns}

\usetikzlibrary{positioning,
                shapes,
                arrows,
                calc,
                3d}

\usepackage{booktabs}

\usepackage{pgfplots}
\usetikzlibrary{pgfplots.groupplots,pgfplots.fillbetween}
\usepackage{subcaption}

\usepackage{bm}
\usepackage{mathtools}
\usepackage{enumerate}
\usepackage{hyperref}
\hypersetup{
    colorlinks = true,
}
\usepackage[capitalize]{cleveref}

\usepackage{algorithm}
\usepackage{algcompatible}

\newcommand{\tikzmark}[1]{\tikz[overlay,remember picture] \node (#1) {\phantom{X}};}

\newcommand*{\SpaceReservedForComments}{2.5cm}%
\newcommand*{\HorizontalOffset}{-1.0em}%
\newcommand*{\VerticalOffset}{0.6ex}%
\newcommand*{\AddNote}[4][]{%
    \begin{tikzpicture}[overlay, remember picture]
        \draw [#1]
            ($(#3.south)+(\HorizontalOffset,\VerticalOffset)$) --  ($(#2.north)+(\HorizontalOffset,\VerticalOffset)$)
            node [align=left, text width=\SpaceReservedForComments-1.0em, pos=0.5, anchor=east] {#4};
    \end{tikzpicture}
}%

\makeatletter
    \algrenewcommand\alglinenumber[1]{\tiny#1:\tikzmark{\arabic{ALG@line}}}
\makeatother

\usetikzlibrary{decorations.pathreplacing}

\usepackage{physics}

\definecolor{TUcol1}{cmyk}{1,0,0,0}
\definecolor{TUcol2}{cmyk}{0.54,0,0.32,0}
\definecolor{TUcol3}{cmyk}{1,0.15,0.4,0}
\definecolor{TUcol4}{cmyk}{1,0.66,0,0.4}
\definecolor{TUcol5}{cmyk}{0.98,1,0,0.35}
\definecolor{TUcol6}{cmyk}{0.82,0,0.21,0.08}
\definecolor{TUcol7}{cmyk}{0.45,0,0.06,0.06}
\definecolor{TUcol8}{cmyk}{0.45,0.2,0,0.07}
\definecolor{TUcol9}{cmyk}{0.02,0.56,0.84,0}
\definecolor{TUcol10}{cmyk}{0.58,1,0,0.02}
\definecolor{TUcol11}{cmyk}{0.19,1,0,0.19}
\definecolor{TUcol12}{cmyk}{0.36,0,1,0}
\definecolor{TUcol13}{cmyk}{0.02,0,0.54,0}

\definecolor{col1}{HTML}{0480B0}
\definecolor{col2}{HTML}{FFB600}
\definecolor{col3}{HTML}{FF2C00}
\definecolor{col4}{HTML}{AA6F39}

\definecolor{chop1}{RGB}{92,76,161}
\definecolor{chop2}{RGB}{87,172,188}
\definecolor{chop3}{RGB}{161,217,163}
\definecolor{chop4}{RGB}{235,245,156}
\definecolor{chop5}{RGB}{250,230,154}
\definecolor{chop6}{RGB}{253,163,93}
\definecolor{chop7}{RGB}{227,81,74}
\definecolor{chop8}{RGB}{156,0,63}

\pgfplotscreateplotcyclelist{MyCyclelist}{%
thick,solid,col1,mark=*,mark size=1,mark size=2,every mark/.append style={solid,line width = 0.5,fill opacity=0.5,col1}\\
thick,dashed,col2,mark=square*,mark size=1,mark size=2,every mark/.append style={solid,line width = 0.5,fill opacity=0.5,col2}\\
thick,densely dashed,col3,mark=triangle*,mark size=1,mark size=2,every mark/.append style={solid,line width = 0.5,fill opacity=0.5,col3}\\
thick,dashdotted,col4,mark=diamond*,mark size=1,mark size=2,every mark/.append style={solid,line width = 0.5,fill opacity=0.5,col4}\\
}

\tikzset{%
p2*/.style = {thick,solid,            col1,    mark=*,         mark size=2,every mark/.append style={solid,line width = 0.5,fill opacity=0.75,col1}},
p2/.style = {thick,solid,            col1,    mark=o,         mark size=2,every mark/.append style={solid,line width = 0.5,fill opacity=0.75,col1}},
p3*/.style = {thick,densely dashed, col2,      mark=square*,   mark size=2,every mark/.append style={solid,line width = 0.5,fill opacity=0.75,col2}},
p3/.style = {thick,densely dashed, col2,      mark=square,   mark size=2,every mark/.append style={solid,line width = 0.5,fill opacity=0.75,col2}},
p4*/.style = {thick,densely dotted,   col3,  mark=triangle*, mark size=2,every mark/.append style={solid,line width = 0.5,fill opacity=0.75,col3}},
p4/.style = {thick,densely dotted,   col3,  mark=triangle, mark size=2,every mark/.append style={solid,line width = 0.5,fill opacity=0.75,col3}},
p5*/.style = {thick,dashdotted,   col4,    mark=diamond*,     mark size=2,every mark/.append style={solid,line width = 0.5,fill opacity=0.75,col4}},
p5/.style = {thick,dashdotted,   col4,    mark=diamond,     mark size=2,every mark/.append style={solid,line width = 0.5,fill opacity=0.75,col4}},
}

\newcommand{\logLogSlopeTriangle}[5]
{
    \small

    \pgfplotsextra
    {
        \pgfkeysgetvalue{/pgfplots/xmin}{\xmin}
        \pgfkeysgetvalue{/pgfplots/xmax}{\xmax}
        \pgfkeysgetvalue{/pgfplots/ymin}{\ymin}
        \pgfkeysgetvalue{/pgfplots/ymax}{\ymax}

        \pgfmathsetmacro{\xArel}{#1}
        \pgfmathsetmacro{\yArel}{#3}
        \pgfmathsetmacro{\xBrel}{#1-#2}
        \pgfmathsetmacro{\yBrel}{\yArel}
        \pgfmathsetmacro{\xCrel}{\xArel}

        \pgfmathsetmacro{\lnxB}{\xmin*(1-(#1-#2))+\xmax*(#1-#2)} 
        \pgfmathsetmacro{\lnxA}{\xmin*(1-#1)+\xmax*#1} 
        \pgfmathsetmacro{\lnyA}{\ymin*(1-#3)+\ymax*#3} 
        \pgfmathsetmacro{\lnyC}{\lnyA+#4*(\lnxA-\lnxB)}
        \pgfmathsetmacro{\yCrel}{\lnyC-\ymin)/(\ymax-\ymin)} 

        \coordinate (A) at (rel axis cs:\xArel,\yArel);
        \coordinate (B) at (rel axis cs:\xBrel,\yBrel);
        \coordinate (C) at (rel axis cs:\xCrel,\yCrel);

        \draw[#5]   (A)-- node[pos=0.5,anchor=north] {1}
                    (B)--
                    (C)-- node[pos=0.5,anchor=west] {#4}
                    cycle;
    }
}

\newcommand{\logLogSlopeTriangleRev}[5]
{
\small

    \pgfplotsextra
    {
        \pgfkeysgetvalue{/pgfplots/xmin}{\xmin}
        \pgfkeysgetvalue{/pgfplots/xmax}{\xmax}
        \pgfkeysgetvalue{/pgfplots/ymin}{\ymin}
        \pgfkeysgetvalue{/pgfplots/ymax}{\ymax}

        \pgfmathsetmacro{\xArel}{#1-#2}
        \pgfmathsetmacro{\yArel}{#3}
        \pgfmathsetmacro{\xBrel}{#1}
        \pgfmathsetmacro{\yBrel}{\yArel}
        \pgfmathsetmacro{\xCrel}{\xArel}

        \pgfmathsetmacro{\lnxB}{\xmin*(1-(#1-#2))+\xmax*(#1-#2)} 
        \pgfmathsetmacro{\lnxA}{\xmin*(1-#1)+\xmax*#1} 
        \pgfmathsetmacro{\lnyA}{\ymin*(1-#3)+\ymax*#3} 
        \pgfmathsetmacro{\lnyC}{\lnyA-#4*(\lnxA-\lnxB)}
        \pgfmathsetmacro{\yCrel}{\lnyC-\ymin)/(\ymax-\ymin)} 

        \coordinate (A) at (rel axis cs:\xArel,\yArel);
        \coordinate (B) at (rel axis cs:\xBrel,\yBrel);
        \coordinate (C) at (rel axis cs:\xCrel,\yCrel);

        \draw[#5]   (A)-- node[pos=0.5,anchor=south] {1}
                    (B)--
                    (C)-- node[pos=0.5,anchor=east] {#4}
                    cycle;
    }
}

\pgfplotsset{
    cycle list name = MyCyclelist,
    width = 0.9\linewidth,
    every axis plot/.append style={thick,mark=square},
    legend cell align={left},
    ylabel near ticks,
    xlabel near ticks,
    legend style={fill=white, fill opacity=0.6, draw opacity=1,text opacity=1},
}

\theoremstyle{thmstyleone}%
\theoremstyle{thmstyletwo}%

\crefname{rem}{Remark}{Remarks}

\crefname{ex}{Example}{Examples}

\theoremstyle{thmstylethree}%

\newcommand{\RevONE}[1]{\textcolor{Green}{#1}}
\newcommand{\RevTWO}[1]{\textcolor{Blue}{#1}}

\begin{document}


\title{Goal-Adaptive Meshing of Isogeometric Kirchhoff-Love Shells}


\author*[1,2]{\fnm{H.M.} \sur{Verhelst}}\email{h.m.verhelst@tudelft.nl}

\author[3]{\fnm{A.} \sur{Mantzaflaris}}\email{angelos.mantzaflaris@inria.fr}

\author[2]{\fnm{M.} \sur{Möller}}\email{m.moller@tudelft.nl}

\author[1]{\fnm{J.H.} \spfx{Den} \sur{Besten}}\email{henk.denbesten@tudelft.nl}

\affil*[1]{\orgdiv{Department of Maritime and Transport Technology}, \orgname{Delft University of Technology}, \orgaddress{\street{Mekelweg 2}, \city{Delft}, \postcode{2628 CD}, \country{The Netherlands}}}

\affil[2]{\orgdiv{Department of Applied Mathematics}, \orgname{Delft University of Technology}, \orgaddress{\street{Mekelweg 4}, \city{Delft}, \postcode{Delft 2628 CD}, \country{The Netherlands}}}

\affil[3]{\orgdiv{AlgebRa, geOmetry, Modeling and AlgoriTHms}, \orgname{Inria Sophia Antipolis - Méditerranée, Université Côte d’Azur}, \orgaddress{\street{2004 route des Lucioles}, \city{Sophia Antipolis cedex}, \postcode{06902}, \country{France}}}


\abstract{Mesh adaptivity is a technique to provide detail in numerical solutions without the need to refine the mesh over the whole domain. Mesh adaptivity in isogeometric analysis can be driven by Truncated Hierarchical B-splines (THB-splines) which add degrees of freedom locally based on finer B-spline bases. Labeling of elements for refinement is typically done using residual-based error estimators. In this paper, an adaptive meshing workflow for isogeometric Kirchhoff-Love shell analysis is developed. This framework includes THB-splines, mesh admissibility for combined refinement and coarsening and the Dual-Weighted Residual (DWR) method for computing element-wise error contributions. The DWR can be used in several structural analysis problems, allowing the user to specify a goal quantity of interest which is used to mark elements and refine the mesh. This goal functional can involve, for example, displacements, stresses, eigenfrequencies etc. The proposed framework is evaluated through a set of different benchmark problems, including modal analysis, buckling analysis and non-linear snap-through and bifurcation problems, showing high accuracy of the DWR estimator and efficient allocation of degrees of freedom for advanced shell computations.}

\keywords{Isogeometric Analysis,Buckling,Dual-Weighted Residual Method,Adaptive Meshing,Kirchhoff-Love Shell}



\maketitle


\section{Introduction}
The idea behind isogeometric analysis (IGA) \cite{Hughes2005} is to bridge the gap between computer aided design (CAD) and finite element analysis (FEA). By employing B-splines or Non-uniform rational B-splines (NURBS) as the basis for FEA, IGA not only provides geometrically exact analysis, but the high smoothness of the spline bases also provides high accuracy per degree of freedom \cite{Sande2020}. The close link with conventional engineering fields such as automotive, offshore, aircraft or civil engineering makes structural analysis with isogeometric analysis a particular field of interest. Besides the performance of the different isogeometric element formulations for thin Kirchhoff-Love shells \cite{Kiendl2009,Kiendl2015a,Verhelst2021,Alaydin2021}, moderately thick Reissner-Mindlin shells \cite{Benson2011,Hu2020,Kiendl2015b,Sobota2017,Benson2013} or thicker solid-like shells \cite{Hosseini2013,Leonetti2018} in static and dynamic simulations, conventional engineering disciplines also rely on accurate modal and (post-)buckling simulations. In addition, the ability to handle complex (multipatch) CAD geometries via trimming \cite{Coradello2020,Guo2018,Leidinger2019} or patch coupling methods \cite{Herrema2019,Leonetti2020b,Bouclier2017,Guo2015,Coradello2021} improves the applicability of IGA in structural engineering. For problems with a large number of degrees of freedom or problems with a large number of load/time steps, mesh adaptivity can play a key role in providing efficient simulations for industrial applications.\\

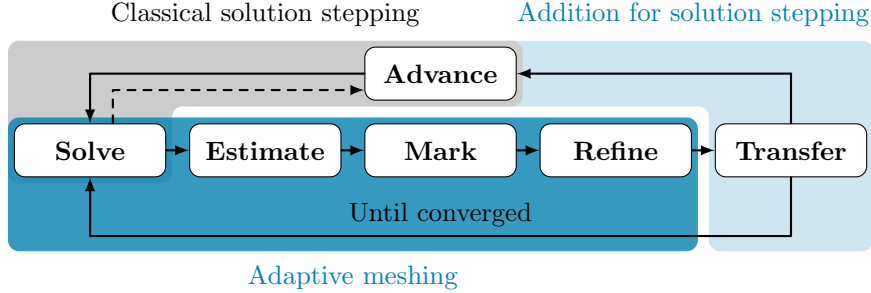
\begin{figure}
\centering
\pgfdeclarelayer{bg}    
\pgfsetlayers{bg,main}  

\tikzstyle{arrow} = [draw, -latex,thick,line cap=round,line join=round]
\tikzstyle{line} = [draw,thick,line cap=round,line join=round]
\tikzstyle{block} = [rectangle, draw, fill=white,
    text width=5em, text centered, rounded corners, minimum height=2em]
\tikzstyle{dummyblock} = [rectangle, fill=none,
    text width=5em, text centered, rounded corners, minimum height=2em]
\tikzstyle{noblock} = [rectangle, draw=none, fill=none,
    text width=5em, rounded corners, minimum height=0em,node distance = 0.05\linewidth, text = blue]
\tikzstyle{textblock} = [text width=7em,text centered]
\tikzstyle{textblock2} = [text width=14em,text centered]
\tikzstyle{textblock4} = [text width=28em,text centered]

\begin{tikzpicture}[node distance = 0.3cm]
    \node[block                      ] (solve)        {\textbf{Solve}};
	\coordinate[below= of solve,inner sep=0] (belowsolve);

    \node[block,  right = of solve   ] (estimate)     {\textbf{Estimate}};
    \path [arrow] (solve) -- (estimate);

    \node[block,  right= of estimate] (mark)         {\textbf{Mark}};
    \path [arrow] (estimate) -- (mark);

    \node[block,  right= of mark    ] (refine)       {\textbf{Refine}};
    \path [arrow] (mark) -- (refine);

    \node[dummyblock,  above= of solve    ] (empty2)       {};

    \node[block,  right= of refine    ] (transfer)       {\textbf{Transfer}};
    \node[dummyblock,  above= of transfer    ] (empty)       {};
    \path [arrow] (refine) -- (transfer);


    \node[block,  above= of mark  ] (advance)       {\textbf{Advance}};
    \path [line] (transfer) -- (empty.center);
	\path [arrow] (empty.center) -- (advance);
    \path [arrow] (advance) -| (solve);

	\node[noblock,  below= of solve] (blockbelowsolve)       {};
	\node[noblock,  below= of refine] (blockbelowrefine)       {};
	\node[noblock,  below= of estimate] (blockbelowestimate)       {};
	\node[noblock,  below= of transfer] (blockbelowtransfer)       {};
	\node[noblock,  above= of empty] (blockaboveempty)       {};
	\node[noblock,  above= of empty2] (blockaboveempty2)       {};
    \path [arrow] (transfer)--(blockbelowtransfer.center)--node[above,midway]{Until converged}(blockbelowsolve.center)--(solve);

    \path [arrow,dashed] ([xshift=-20]solve.north east) |- ([yshift=4]advance.south west);

\begin{pgfonlayer}{bg}    
	\draw[draw=col1!20,fill=col1!20][rounded corners,name=B]
([xshift=-2,yshift=2]transfer.north west) --
([xshift=-2,yshift=-2]blockbelowtransfer.south west) --
([xshift=2,yshift=-2]blockbelowtransfer.south east) --
([xshift=2,yshift=2]transfer.north east) --
([xshift=2,yshift=2]empty.north east) --
([xshift=2,yshift=2]advance.north west) --
([xshift=2,yshift=-2]advance.south west) --
([xshift=-2,yshift=-2]empty.south west)
--cycle;
	\draw[draw=black!20,fill=black!20][rounded corners,name=B]
([xshift=-2,yshift=-2]solve.south west) --
([xshift=2,yshift=-2]solve.south east) --
([xshift=2,yshift=-2]empty2.south east) --
([xshift=2,yshift=-2]advance.south east) --
([xshift=2,yshift=2]advance.north east) -- node[yshift=2,midway](steppinglabel){}
([xshift=-2,yshift=2]empty2.north west)
--cycle;
	\draw[draw=col1,fill=col1,opacity=0.8][rounded corners,name=B] ([xshift=-2,yshift=2]solve.north west) --
([xshift=2,yshift=2]refine.north east) --
([xshift=2,yshift=-2]blockbelowrefine.south east)--node[yshift=-2,midway](meshinglabel){}
([xshift=-2,yshift=-2]blockbelowsolve.south west)
--cycle;
\end{pgfonlayer}

\node[textblock2,anchor=south,color=col1,yshift=4,xshift=-8] at (empty.north west) {Addition for solution stepping};
\node[textblock4,anchor=north,color=col1] at (meshinglabel) {Adaptive meshing};
\node[textblock2,anchor=south,color=black] at (steppinglabel) {Classical solution stepping};

\end{tikzpicture}
\caption{A typical flowchart for an adaptive meshing routine. The \emph{classical solution stepping} depicts a process without adaptive meshing. Here, a solution is obtained by the \emph{solve} and the solution is advanced (e.g. in time or load step) and recomputed. The \emph{adaptive meshing step} denotes the additional operations for mesh adaptivity and the \emph{Addition for solution stepping} includes an additional \emph{transfer} step in case the adaptive meshing method is applied to solution. The \emph{Estimate} block provides an error estimation with local contributions per element or per degree of freedom (DoF). The \emph{Mark} block contains a marking rule that marks regions for refinement based on a specific rule. The \emph{Refine} block transforms the current mesh to a new mesh, where regions are refined and coarsened based on the marking rule. The block \emph{Transfer} transfers the previous solution to the new mesh, so that it can be used to recompute the present interval on a modified mesh. This recomputation is performed again in the \emph{Solve} block and follows through the subsequent blocks, until an adaptivity criterion is reached. For example, a criterion that requires the total error in the mesh to be within certain bounds.}
\label{fig:flowchart}
\end{figure}

A loop in an adaptive isogeometric method (AIGM) consists of the steps \emph{solve} the Partial Differential Equation (PDE) at hand, \emph{estimate} element-wise error contributions, \emph{mark} regions for refinement, \emph{refine (coarsen)} marked regions \cite{Buffa2016}, see \cref{fig:flowchart}. Here, localised regions can be defined element-wise or function-wise. The AIGM process can be repeated in an iterative manner (e.g. for static, buckling or modal analysis), until satisfactory accuracy is achieved or it can be applied (iteratively) within a time/load stepping procedure. A broad overview of the mathematical foundations of AIGMs is given in \cite{Buffa2022}. In previous works, AIGMs are developed for different applications (\emph{solve}), using different \emph{estimation} strategies, \emph{marking} strategies and often for mesh \emph{refinement}, with a few also providing \emph{coarsening} strategies \cite{DAngella2020a,Carraturo2019,Lorenzo2017,Hennig2018,Garau2018}. \\

\begin{description}
 \item[Solve] The \emph{solve} block contains the partial differential equation (PDE) at hand. It can be a physics-based problem, e.g., to solve shell \cite{Antolin2020,Coradello2020,Coradello2020a}, linear elasticity \cite{Verhoosel2015} or free-surface flow \cite{Kuru2014} problems. Alternatively, the solve step can involve a non-physics PDE, e.g. for mesh generation \cite{Hinz2020}.
 \item[Estimate] Determination of localised errors is done in the \emph{estimate} block. In the works \cite{Antolin2020,Coradello2020,Coradello2020a}, an error estimator based on a residual-like variational problem in the so-called \emph{bubble-space} was presented for Kirchhoff plates, Kirchhoff-Love shells and trimmed domains. This method has proven a large decrease in CPU time compared to a residual-based error estimator in the strong form, due to its easy parallelisation and the small block-structure of the linear system to solve. As an alternative to this method, error estimation can also be performed in a goal-oriented fashion, e.g., by the Dual-Weighted Residual (DWR) method. This method has been applied in the FEA context in various works \cite{Rannacher2004,Gedicke2013,Hartmann2010,Hartmann2002,Hartmann2003,Hartmann2002a,Moller2006,Cliffe2010} and was used in the works \cite{VanDerZee2011,Kuru2014} for Poisson and free-boundary problems, in \cite{Dede2012} for a geometrically non-linear rod, in \cite{Hinz2020} \RevTWO{for} PDE-based domain parametrisations and in \cite{Verhoosel2015} for micromechanical modelling of trabecular bone. Goal-oriented refinement in general provides localised error estimates by solving a linear adjoint problem on the current space and an enriched space.
 \item[Mark] As soon as localised error contributions are known, regions can be marked for refinement. This marking is mostly done using the Dörfler marking strategy  \cite{Dorfler1996}, as in \cite{Carraturo2019,Buffa2016,Gedicke2013}, which involves marking the regions with the largest error contributions until their sum exceeds a certain percentage of the total error. An alternative is to mark the regions with an error higher than a threshold (an absolute threshold based on the maximum error) \cite{Antolin2020,Giannelli2016} or based on a relative threshold taking a fixed percentage of the total number of cells for refinement. In \cite{Giannelli2016} the latter two strategies are discussed.
 \item[Refine] Local refinement for adaptive meshing in isogeometric analysis is enabled by the use of Hierarchical B-splines (HB-splines) \cite{Vuong2011}, Truncated Hierarchical B-splines (THB-splines) \cite{Giannelli2012,Giannelli2016} or T-splines \cite{Bazilevs2010} amongst other spline constructions, which are reviewed in \cite{Hennig2016}. HB-splines provides a nested, linear dependent space that violates the partition of unity property. To preserve the latter, THB-splines have been introduced in \cite{Giannelli2012}. For (T)HB-splines \emph{suitable grading} to generate admissible meshes should be taken into account in order to guarantee a bounded error \cite{Buffa2016}, for which algorithms have been presented in \cite{Bracco2018a}. In the work of \cite{Hennig2017}, a distinction is made between \textit{greedy} and \textit{safe} refinement, the former being a refinement of cells with a 1-level difference with adjacent cells and the latter being a refinement complying with the refinement neighborhoods defined in \cite{Buffa2016}. Besides for adaptive meshing for solving PDEs \cite{Hinz2020,DAngella2020a,Carraturo2019,Lorenzo2017,Hennig2018,Garau2018}, THB-splines have also been succesfully applied in the context of fitting \cite{Bracco2022,Kiss2014}.
 \item[Transfer] The \emph{transfer} from previous time/load-steps onto a new mesh can be done using different methodologies. In the work of \cite{Hennig2018} different least-squares approaches are provided. Furthermore, quasi-interpolation \cite{Speleers2016,Giust2020} is a technique that can be used to transfer solutions between hierarchical meshes.
\end{description}

In this paper, we employ goal-oriented adaptive refinement for isogeometric thin shell analysis to facilitate THB-adaptive meshing for a variation of structural analysis problems. \RevTWO{The developed framework is versatile in terms of the goal functional being used and provides an adaptive meshing strategy for linear an non-linear static, modal, buckling and post-buckling problems. In brief, the contribution of the paper is threefold. Firstly, we use the Dual-Weighted Residual (DWR) method to derive novel error estimators for structural shell analysis, given goal functionals based on displacements, (principal) stresses and strains, forces and moments. Secondly, we employ the eigenvalue DWR from \cite{Lathouwers2011,Rannacher2004} for error estimations for modal and buckling analyses. Thirdly, the goal functionals are used to drive an adaptive meshing strategy with suitable grading and efficient transfer of solutions by quasi-interpolation method on hierarchical spline spaces \cite{Speleers2016,Giust2020}. This adaptive meshing strategy is applied to non-linear shell analysis with focus on buckling problems with snap-through and bifurcation instabilities - being new applications in the realm of adaptive meshing research for nonlinear shell problems. It should be noted that the present framework is developed for  isogeometric Kirchhoff-Love shells - since it provides a natural separation of bending and membrane terms - but it is easily adapted for other shell formulations. By defining a frame work for 2-dimensional parametric domains and by presenting a wide range of mechanics-inspired goal functionals, the present work extends an earlier work by \cite{Dede2012} for geometrically non-linear rods.}

The paper is structured as follows. In \cref{sec:shell}, the isogeometric Kirchhoff-Love shell analysis proposed by \cite{Kiendl2009} is briefly revised and some basic concepts for structural analysis computations are given. In \cref{sec:DWR}, the Dual-Weighted Residual (DWR) method is provided for the isogeometric Kirchhoff-Love shell using the membrane and flexural strain split. However, it can be used for general elasticity problems. Moreover, the section provides the DWR method for eigenvalue problems to compute error estimators for modal and buckling analyses. Thereafter, \cref{sec:refinement} provides the details for adaptivity for isogeometric analysis. This includes the concept of Truncated-Hierarchical B-splines (THB-splines) and admissible refinement. Furthermore, the \emph{mark} and \emph{transfer} operations are described. In \cref{sec:algorithm}, a summary of the preceeding sections is provided by means of a global algorithm for the AIGM for structural analysis computations with load-stepping. In \cref{sec:results} the present work is evaluated on numerical benchmark problems, ranging from linear problems with analytical solutions to non-linear shell problems. Finally, \cref{sec:conclusions} provides conclusions and and outlook based on this work.

\section{Isogeometric Kirchhoff-Love Shell Analysis}\label{sec:shell}
In this section, we provide a brief background on the Kirchhoff-Love shell formulation. For more details on this formulation, we refer to \cite{Kiendl2009,Kiendl2015a,Roohbakhshan2017,Sauer2017,Verhelst2021}.

\RevTWO{Check punctuation!!}\\

\subsection{Shell Kinematics}
Since Kirchhoff-Love shells satisfy the Kirchhoff Hypothesis \cite{Reddy2014}, the coordinates $\vb{x}$ of any parametric point $\bm{\theta} = (\theta^1, \theta^2, \theta^3)$ in the shell surface can be represented by the surface position $\vb{r}(\theta^1,\theta^2)$ and contribution in normal direction $\theta^3\vb{a}_3$ as
\begin{equation}\label{eq:coordinate}
 \vb{x}(\bm{\theta}) = \vb{r}(\theta^1,\theta^2) + \theta^3\vb{a}_3,
\end{equation}
Given the covariant basis of the surface $\vb{r}$, defined by $\vb{a}_{\alpha}, \alpha=1,2$ and the orthogonal unit normal $\vb{a}_3$, the covariant basis of $\vb{x}$ is defined as follows:
\begin{equation}
 \vb{g}_\alpha = \vb{x}_{,\alpha} =  \vb{a}_\alpha + \theta^3 \vb{a}_{3,\alpha},\quad \vb{g}_3 = \vb{x}_{,3} = \vb{a}_3.
\end{equation}
Given the second fundamental form $b_{\alpha\beta} = \vb{a}_3\cdot\vb{a}_{\alpha,\beta}=-\vb{a }_{3,\beta}\cdot\vb{a}_\alpha$ and the metric coefficients defined as
\begin{equation}\label{eq:metric}
 g_{\alpha\beta} = \vb{g}_\alpha \cdot \vb{g}_\beta= a_{\alpha\beta}-2\theta^3b_{\alpha\beta},
\end{equation}
the contravariant basis vectors $\vb{g}^\alpha$ can simply be obtained by $\vb{g}^{\alpha} = g^{\alpha\beta}\vb{g}_\beta$. The undeformed configuration $\vb{r}$ and the deformed configuration $\mathring{\vb{r}}$ of the surface are related by $\vb{r}=\mathring{\vb{r}}+\vb{u}$. From the defintion of the deformation gradient $\vb{F} = \vb{g}_i \otimes \mathring{\vb{g}}^i$, the deformation tensor $\vb{C}$ can be obtained:
\begin{equation}
        \vb{C} = \vb{F}^\top\vb{F} = \vb{g}_i\cdot\vb{g}_j \: \mathring{\vb{g}}^i\otimes\mathring{\vb{g}}^j = g_{ij} \: \mathring{\vb{g}}^i\otimes\mathring{\vb{g}}^j.
\end{equation}
Note that the deformation tensor is defined in the contravariant undeformed basis $\mathring{\vb{g}}^i\otimes\mathring{\vb{g}}^j$. For Kirchhoff-Love shells, it is known that $g_{\alpha3}=g_{3\alpha}=0$, hence this implies $C_{\alpha3}=C_{3\alpha}=0$, since $g_{33}=1$, which implies $C_{33}$ to be one and meaning that the thickness remains constant under deformation. As a result, the Green-Lagrange strain tensor $\vb{E}=E_{\alpha\beta}\:\vb{\mathring{g}}^\alpha\otimes \vb{\mathring{g}}^\beta$ and its decomposition to membrane and bending contributions ($\varepsilon$ and $\kappa$, respectively) is \cite{Kiendl2009,Kiendl2015a}:
\begin{equation}\label{eq:strain}
\begin{aligned}
 E_{\alpha\beta} &= \frac{1}{2}\qty(g_{\alpha\beta}-\mathring{g}_{\alpha\beta}) = \frac{1}{2}\qty( (a_{\alpha\beta} - \mathring{a}_{\alpha\beta}) - 2\theta^3\qty( b_{\alpha\beta} - \mathring{b}_{\alpha\beta} ) ) \\ &= \varepsilon_{\alpha\beta} + \theta_3\kappa_{\alpha\beta}.
\end{aligned}
\end{equation}

\subsection{Constitutive relation}\label{subsec:kinematics}
The constitutive relations for the Kirchhoff-Love shell relate the Green-Lagrange strain tensor $\vb{E}$ to the second Piola-Kirchhoff stress tensor $\vb{S}$. For linear elastic materials, this is achieved by:
\begin{equation}
    S^{\alpha\beta} = \mathbb{C}^{\alpha\beta\gamma\delta}E_{\gamma\delta}
\end{equation}
where $\mathbb{C}=\mathbb{C}^{\alpha\beta\gamma\delta}\:\mathring{\vb{g}}_i\otimes\mathring{\vb{g}}_j\otimes\mathring{\vb{g}}_k\otimes\mathring{\vb{g}}_l$ is the material tensor, which takes for linear materials the form $\mathbb{C}^{\alpha\beta\gamma\delta} = 4\frac{\lambda\mu}{\lambda+2\mu}\mathring{g}^{\alpha\beta}\mathring{g}^{\gamma\delta} + 2\mu\qty(\mathring{g}^{\alpha\delta}\mathring{g}^{\beta\gamma}+\mathring{g}^{\alpha\gamma}\mathring{g}^{\beta\delta})$ \cite{Goyal2015}. For non-linear hyperelastic constitutive relations, the stress and material tensors are derived from the 3D constitutive relations for (in)compressible materials and due to through-thickness integration, the shell normal force and bending moment tensors $\vb*{n}=n^{\alpha\beta}\:\vb{\mathring{g}}_\alpha\otimes \vb{\mathring{g}}_\beta$ and $\vb*{m}=m^{\alpha\beta}\:\vb{\mathring{g}}_\alpha\otimes \vb{\mathring{g}}_\beta$, respectively, are defined as
\begin{equation}
    n^{\alpha\beta}(\vb{u}) = \int_{T} S^{\alpha\beta}(\vb{u})\dd{\theta^3}, \quad
    m^{\alpha\beta}(\vb{u}) = \int_{T} \theta^3 S^{\alpha\beta}(\vb{u})\dd{\theta^3},
\end{equation}
where $T=[-t/2,t/2]$ is the through-thickness domain. For more details on hyperelastic material models, the reader is referred to \cite{Kiendl2015a,Roohbakhshan2017} and specifically for stretch-based ones to \cite{Verhelst2021}.\\

\subsection{Variational Formulation}\label{subsec:variational}
The shell internal and external equilibrium equations in variational form are derived by the principle of virtual work \cite{Kiendl2009,Kiendl2015a}. The weak formulation follows from the principle of virtual work with virtual displacements $\bm{\phi}$:
\begin{equation}\label{eq:variational_form}
\begin{aligned}
\text{Find }\vb{u}\in\mathbcal{S}\text{ s.t.}
    \mathbcal{W}(\vb{u},\bm{\phi}) &:= \delta W^{\text{int}} - \delta W^{\text{ext}} \\ & = \int_\Omega \vb*{n}(\vb{u}):\bm{\varepsilon}'(\vb{u},\bm{\phi}) + \vb*{m}(\vb{u}):\bm{\kappa}'(\vb{u},\bm{\phi}) \dd{\Omega} - \int_\Omega \vb{f}(\vb{u})\cdot\bm{\phi} \dd{\Omega},\\
    \forall \bm{\phi}\in\mathbcal{S}
\end{aligned}
\end{equation}
With $\vb{f}(\vb{u})$ the surface load acting on the mid-surface, for the sake of generality defined as a function of the displacements $\vb{u}$ (e.g. a follower pressure $p$ gives $\vb{f}(\vb{u}) = p\vb{a}_3(\vb{u})$). Furthermore, $\bm{\varepsilon}'(\vb{u},\bm{\phi})$ and $\bm{\kappa}'(\vb{u},\bm{\phi})$ are the virtual strain components given displacements $\vb{u}$ and being linear with respect to variation $\bm{\phi}$, hence $\mathbcal{W}(\vb{u},\bm{\phi})$ is also linear in its second argument. The coefficients of the variations of the membrane force and bending moment tensors are
\begin{equation}
\begin{aligned}
(n')^{\alpha\beta}(\vb{u},\bm{\phi}) &= \int_{T} \mathbb{C}^{\alpha\beta\gamma\delta}(\vb{u})\dd{\theta^3}\varepsilon'_{\gamma\delta}(\vb{u},\bm{\phi}) + \int_{T} \theta^3\mathbb{C}^{\alpha\beta\gamma\delta}(\vb{u})\dd{\theta^3}\kappa'_{\gamma\delta}(\vb{u},\bm{\phi}),\\
(m')^{\alpha\beta}(\vb{u},\bm{\phi}) &= \int_{T} \theta^3\mathbb{C}^{\alpha\beta\gamma\delta}(\vb{u})\dd{\theta^3}\varepsilon'_{\gamma\delta}(\vb{u},\bm{\phi}) + \int_{T} \qty(\theta^3)^2\mathbb{C}^{\alpha\beta\gamma\delta}(\vb{u})\dd{\theta^3}\kappa'_{\gamma\delta}(\vb{u},\bm{\phi}).
\end{aligned}
\end{equation}
Linearizing the virtual work from \cref{eq:variational_form} provides the continuous equivalent of the Jacobian or \emph{tangential stiffness matrix} for Newton iterations which will be performed to solve the non-linear weak formulation \cref{eq:variational_form} in a discrete setting \cite{Sauer2017}:
\begin{equation}\label{eq:linearisationN}
    \begin{aligned}
    \mathbcal{W}'(\vb{u},\bm{\phi},\bm{\psi}) &:= \int_{\Omega} \vb*{n}'(\vb{u},\bm{\psi}):\bm{\varepsilon}'(\vb{u},\bm{\phi})+\vb*{n}(\vb{u}):\bm{\varepsilon}''(\vb{u},\bm{\phi},\bm{\psi}) \\
    &\quad\quad+ \vb*{m}'(\vb{u},\bm{\psi}):\bm{\kappa}'(\vb{u},\bm{\phi}) + \vb*{m}:\bm{\kappa}''(\vb{u},\bm{\phi},\bm{\psi})\dd{\Omega}\\
    &\quad- \int_{\Omega} \vb{f}'(\vb{u},\bm{\psi})\cdot\bm{\phi}\dd{\Omega},
    \end{aligned}
\end{equation}
where $\bm{\varepsilon}''(\vb{u},\bm{\phi},\bm{\psi})$ and $\bm{\kappa}(\vb{u},\bm{\phi},\bm{\psi})$ are the second variations of the membrane and bending strains and $\vb{f}'$ is the first variation of the applied force, being nonzero when the force is depending on the displacements $\vb{u}$. For the details on these formulations, we refer to previous publications \cite{Sauer2017,Goyal2015,Kiendl2009}. It should be noted that in the undeformed case, $\vb{u}=\vb{0}$, the internal membrane forces and bending forces, $\vb*{n}(\vb{u})$ and $\vb*{m}(\vb{u})$, respectively, vanish. As a result, the continuous equivalent for the linear stiffness matrix is:
\begin{equation}\label{eq:linearisationNlinear}
    \begin{aligned}
    \mathring{\mathbcal{W}}'(\bm{\phi},\bm{\psi})
     &= \int_{\Omega} \mathring{\vb*{n}}'(\bm{\psi}):\mathring{\bm{\varepsilon}}'(\bm{\phi}) + \mathring{\vb*{m}}'(\bm{\psi}):\mathring{\bm{\kappa}}'(\bm{\phi})\dd{\Omega} - \int_{\Omega} \mathring{\vb{f}}'(\bm{\psi})\cdot\bm{\phi}\dd{\Omega},
    \end{aligned}
\end{equation}
where the $\mathring{\cdot}$ denotes tensors and functions on the undeformed geometry, i.e. with $\vb{u}=\vb{0}$.\\

In our implementation, the tangent stiffness matrix is computed using appropriate Gauss-Lengendre quadrature for each element in the hierarchical mesh. We note that more efficient numerical integration approaches exist~\cite{Pan2021,Pan2022,Giannelli2022} for hierarchical splines that might further reduce the computational cost.\\

\subsection{Structural Analysis}
\label{subsec:structuralAnalysis}
In the present paper, we will provide different goal functionals for different structural analysis applications. Therefore, we briefly recall the different structural analysis types. Firstly, in case of \emph{static analysis}, the problem as in \cref{eq:variational_form} is solved. In case of \emph{quasi-static analysis}, load and/or displacement steps are performed successively and in each step, a static solve is performed. Typically, one writes \cref{eq:variational_form} for load-control as
\begin{equation}\label{eq:variational_formALM}
\begin{aligned}
\text{Find }\vb{u}\in\mathbcal{S}\text{ s.t.}
    \mathbcal{W}(\vb{u},\bm{\phi},\lambda) &:= \delta W^{\text{int}} - \lambda\delta W^{\text{ext}} \\ & = \int_\Omega \vb*{n}(\vb{u}):\bm{\varepsilon}'(\vb{u},\bm{\phi}) + \vb*{m}(\vb{u}):\bm{\kappa}'(\vb{u},\bm{\phi}) \dd{\Omega} - \int_\Omega \lambda\vb{f}_0\cdot\bm{\phi} \dd{\Omega},\\
    \forall \bm{\phi}\in\mathbcal{S}, \lambda\in\mathbb{R},
\end{aligned}
\end{equation}
where $\lambda$ is the load factor scaling the reference load $\vb{f}_0$. Quasi-static simulations can be solved using simple load or displacement controlled schemes, using arc-length continuation such as Riks' method or Crisfield's method \cite{Riks1972,Crisfield1981}. When quasi-static analysis is performed for post-buckling analysis, one or multiple bifurcation points are passed by definition. On a bifurcation point, the determinant of the tangential stiffness matrix $K$ is equal to zero, hence this matrix is singular. To cope with instabilities, \emph{a priori} perturbations can be applied to the geometry, or a procedure for approximating singular points \cite{Wriggers1988} can be used. In our previous work, we provide more details on arc-length continuation for post-buckling analysis without providing \emph{a priori} perturbations  \cite{Verhelst2019}.\\

In the case of \emph{modal analysis} and \emph{buckling analysis}, a generalised eigenvalue problem needs to be solved. These eigenvalue problems have the general form
\begin{equation}
    \text{find}\:(\mu,\vb{v})\in\mathbb{R}\times\mathbcal{S}\:\text{s.t.}\:\mathbcal{A}(\vb{v},\bm{\phi}) = \mu \mathbcal{B}(\vb{v},\bm{\phi})\:\forall\bm{\phi}\in\mathbcal{S}
\end{equation}
Where $\mu$ provides the eigenfrequency in modal analysis and the critical load factor in buckling analysis and where $\vb{v}$ denotes the vibration or buckling mode shape. The operators $\mathbcal{A}$ and $\mathbcal{B}$ are bi-linear. For buckling analysis, $\mathbcal{A}(\vb{v},\bm{\phi})=\mathbcal{W}'(\vb{u}_L,\vb{v},\bm{\phi})$ and $\mathbcal{B}(\vb{v},\bm{\phi})=\mathbcal{W}'(\vb{0},\vb{v},\bm{\phi})$ with $\vb{u}_L$ the pre-buckling solution given load $\lambda_L$. For modal analysis, $\mathbcal{A}(\vb{v},\bm{\phi})=\mathbcal{W}'(\vb{0},\vb{v},\bm{\phi})$ and $\mathbcal{B}(\vb{v},\bm{\phi})=\mathbcal{M}(\vb{v},\bm{\phi})$ with $\mathbcal{M}$ the mass operator:
\begin{equation}\label{eq:mass}
 \mathbcal{M}(\vb{v},\bm{\phi}) = \int_\Omega \rho \vb{v}\bm{\phi}\dd{\Omega}
\end{equation}
Where $\rho$ is the density function over the surface.

\section{Dual-Weighted Residual method}\label{sec:DWR}
This section elaborates on the Dual-Weighted Residual (DWR) method \cite{Becker2001,Bangerth2003}, which is used in the \emph{Estimate} step of \cref{fig:flowchart}. The DWR is a method to compute the \emph{a posteriori} error of a solution in terms of a given goal functional of interest, by solving a linear dual problem. The DWR provides a global estimate of the error, but given a partition of unity of the spline space, it can be used to provide an error contribution per basis function.

\subsection{General Framework}
\label{subsec:DWR_general}
The general framework of the dual weighted residual (DWR) is presented by \cite{Becker2001,Bangerth2003,Hartmann2002}. For the sake of completeness, we provide a brief overview of the DWR here. Consider the following non-linear problem to solve
\begin{equation}\label{eq:nlproblem}
    \text{find}\:\vb{u}\in\mathbcal{S}\: \text{s.t.}\:\mathbcal{W}(\vb{u},\bm{\phi})= 0\:\forall\bm{\phi}\in\mathbcal{S},
\end{equation}
where $\mathbcal{W}(\vb{u})$ is a semi-linear operator, $\vb{u}$ is the solution, $\bm{\phi}$ is a test function and $\mathbcal{S}$ is a suitably chosen vector space including $\vb{u}\in\mathbcal{S}$. The approximation of $\vb{u}$, denoted by $\vb{u}_h$ can be found by solving the discrete counterpart of \cref{eq:nlproblem}
\begin{equation}\label{eq:nlproblem_discretised}
    \text{find}\:\vb{u}_h\in\mathbcal{S}_h^p\: \text{s.t.}\:\mathbcal{W}(\vb{u}_h,\bm{\phi}_h)= 0\:\forall\bm{\phi}\in\mathbcal{S}_h^p,
\end{equation}
where $\vb{u}_h$ and $\bm{\phi}_h$ are the discrete counterparts of $\vb{u}$ and $\bm{\phi}$, respectively, and the space $\vb{S}_h^p\subset\vb{S}$ is a function space on the (isogeometric) mesh $\mathbcal{T}_h^p(\Omega)$ with mesh size $h$ and order $p$ covering the computational domain $\Omega$. The solution to this problem is typically obtained by iteratively solving
\begin{equation}
    \text{find}\:\bm{\phi}_h\in\mathbcal{S}_h^p\: \text{s.t.}\:\mathbcal{W}'(\vb{u}_h,\bm{\phi}_h,\bm{\psi}_h) = \mathbcal{R}(\vb{u}_h,\bm{\psi}_h) \:\forall\bm{\psi}_h\in\mathbcal{S}_h^p,
\end{equation}
\RevONE{
while updating the discrete solution. Here, the residual is defined as
\begin{equation}
    \mathbcal{R}(\vb{u}_h,\bm{\phi}_h) = -\mathbcal{W}(\vb{u}_h,\bm{\phi}_h).
\end{equation}
Let us now define a non-linear and differentiable goal functional $\mathbcal{L}(\vb{u})$ or quantity of interest, such that
\begin{equation}\label{eq:goalfunction}
    \Delta \mathbcal{L}(\vb{u}_h) = \mathbcal{L}(\vb{u}) - \mathbcal{L}(\vb{u}_h).
\end{equation}
Then, from Proposition 4.1 of \cite{Hartmann2006}, it follows that:
\begin{equation}\label{eq:errorest}
    \Delta\mathbcal{L}(\vb{u}_h) = \mathbcal{R}(\vb{u}_h,\bm{\xi} - \bm{\xi}_h) \approx \mathbcal{R}(\vb{u}_h,\tilde{\bm{\xi}} - \bm{\xi}_h).
\end{equation}
Here, the solutions $\bm{\xi}\in\mathbcal{S}$ and $\bm{\xi}_h\in\mathbcal{S}^p_h$ are the exact and discrete solutions to the adjoint problem defined using the mean value linearizations of $\mathbcal{W}$ and $\mathbcal{L}$, see Equations 10 to 12 of \cite{Hartmann2006}. Sinze the exact dual solution $\bm{\xi}$ is not available, it is approximated by $\tilde{\bm{\xi}}\in\tilde{\mathbcal{S}}$. The discrete dual solution $\bm{\xi}_h\in\mathbcal{S}^p_h$ is obtained by solving the following discrete adjoint problem:
\begin{equation}\label{eq:adjoint_discrete}
    \text{find}\:\bm{\xi}_h\in\mathbcal{S}_h^p\: \text{s.t.}\:\mathbcal{W}'(\vb{u}_h,\bm{\zeta}_h,\bm{\xi}_h) = \mathbcal{L}'(\vb{u}_h,\bm{\zeta}_h) \:\forall\bm{\zeta}_h\in\mathbcal{S}_h^p.
\end{equation}
The approximation $\tilde{\bm{\psi}}\in\tilde{\mathbcal{S}}$ is obtained by solving the adjoint problem in an enriched space, i.e.
\begin{equation}\label{eq:adjointenriched}
    \text{find}\:\tilde{\bm{\xi}}_h\in\tilde{\mathbcal{S}}_h^{p}\: \text{s.t.}\:\mathbcal{W}'(\vb{u}_h,\tilde{\bm{\zeta}}_h,\tilde{\bm{\xi}}_h) = \mathbcal{L}'(\vb{u}_h,\tilde{\bm{\zeta}}_h) \:\forall\tilde{\bm{\zeta}}_h\in\tilde{\mathbcal{S}}_h^p,
\end{equation}
with $\tilde{\bm{\xi}}_h$ and $\tilde{\bm{\zeta}}_h$ the dual solution and test functions on the enriched space $\tilde{\mathbcal{S}}_h^p$, respectively. A choice for $\tilde{\mathbcal{S}}_h^p$ is to use the same mesh as for $\mathbcal{S}_h^p$, with the same regularity but with a higher degree, i.e. $\tilde{\mathbcal{S}}_h^p = \mathbcal{S}_h^{p+1}$. This is easily achieved using spline bases. When using B-Splines, one can repeat all knots of the knot vector an extra time compared to the original basis, such that $\mathbcal{S}_h^p\subset\mathbcal{S}_h^{p+1}\subset\mathbcal{S}$ is a nested space.\\
}

Finally, using \cref{eq:errorest} together with the dual solution $\bm{\xi}_h\in\mathbcal{S}_h^p$ and the enriched dual solution $\tilde{\bm{\xi}}_h\in\tilde{\mathbcal{S}}_h^p$, an estimate for the global error with respect to the goal functional $\mathbcal{L}$ can be obtained. To obtain the local element-wise error estimations $r_i$ for element $\omega_i\in\mathcal{T}_h^p(\Omega)$, such that
\begin{equation}\label{eq:funccontributions}
    \Delta\mathbcal{L}(\vb{u}_h) = \mathbcal{R}(\vb{u}_h,\tilde{\bm{\xi}}_h-\bm{\xi}_h) = \sum_{\omega_i\in\mathcal{T}_h^p(\Omega)} r_i,
\end{equation}
element-wise integration of \cref{eq:errorest} is simply performed to obtain $r_i$. However, as discussed in \cref{subsec:labeling} it can be beneficial to have strictly positive element error contributions for element labeling. One can either take the absolute values of $r_i$ or one can integrate the squared norm of the integrant in \cref{eq:errorest} to ensure positivity of element error contributions. Obviously, the sum of the element errors would not be equal to $\Delta\mathbcal{L}$.\\

For Kirchhoff-Love shells specifically, the operator $\mathbcal{W}(\vb{u},\bm{\phi})$ and its linearisation $\mathbcal{W}'(\vb{u},\bm{\phi},\bm{\psi})$ are used to perform the DWR analysis.


\subsection{Eigenvalue problems}\label{subsec:eigenvalueproblem}
When the problem of interest is an eigenvalue problem, the DWR routine is slightly different. Here, we follow the works \cite{Gedicke2013,Lathouwers2011,Becker2001,Giani2012,Cliffe2010,Bangerth2003} to give a brief overview of the DWR for eigenvalue problems.\\

Let us consider the following eigenvalue problem
\begin{equation}\label{eq:eigenvalueproblem}
    \text{find}\:(\mu,\vb{v})\in\mathbb{R}\times\mathbcal{S}\:\text{s.t.}\:\mathbcal{A}(\vb{v},\bm{\phi}) = \mu \mathbcal{B}(\vb{v},\bm{\phi})\:\forall\bm{\phi}\in\mathbcal{S}.
\end{equation}
Here, $\mathbcal{A}$ and $\mathbcal{B}$ are bi-linear operators. For uniqueness of the problem, the discrete eigenvectors $\vb{v}_h$ are normalised by the condition \cite{Lathouwers2011}
\begin{equation}\label{eq:normalisation}
    \mathbcal{B}(\vb{v},\vb{v}) = 1.
\end{equation}
Typically, discretizing the system gives the following:
\begin{equation}\label{eq:eigenvalueproblem_discrete}
    \text{find}\:(\mu_h,\vb{v}_h)\in\mathbb{R}\times\mathbcal{S}_h^p\:\text{s.t.}\:\mathbcal{A}(\vb{v}_h,\bm{\phi}_h) = \mu \mathbcal{B}(\vb{v}_h,\bm{\phi}_h)\:\forall\bm{\phi}_h\in\mathbcal{S}_h^p,
\end{equation}
where the eigenpairs $\hat{\vb{v}}_h = (\mu_h,\vb{v}_h)$ are the solutions of the eigenvalue problem. \RevTWO{
In addition, the adjoint eigenvalue problem is defined by the eigenvalue problem \cite{Bangerth2003}:
\begin{equation}\label{eq:eigdual}
    \text{find}\:(\eta,\bm{\psi})\in\mathbb{R}\times\mathbcal{S}\:\text{s.t.}\:\mathbcal{A}(\bm{\psi},\bm{\phi}) = \eta \mathbcal{B}(\bm{\psi},\bm{\phi})\:\forall(\bm{\phi})\in\mathbcal{S},
\end{equation}
Of for which the normalization similar to \cref{eq:normalisation} is used for the dual eigenvectors $\bm{\psi}$
\begin{equation}\label{eq:normalisation_adj}
    \mathbcal{B}(\vb{v},\bm{\psi}) = 1.
\end{equation}
}

To derive the DWR method for the eigenvalue problem in \cref{eq:eigenvalueproblem}, the functional $\mathbcal{V}(\cdot,\cdot)$ is defined, such that the following problem should be solved
\begin{equation}\label{eq:Gfunctional}
    \begin{aligned}
        &\text{Find}\:\hat{\vb{v}}=(\mu,\vb{v})\in\mathbb{R}\times\mathbcal{S}\:\text{s.t.}\\
        &\quad  \mathbcal{V}(\hat{\vb{v}},\hat{\bm{\phi}}) = \mu \mathbcal{B}(\vb{v},\bm{\phi}) - \mathbcal{A}(\vb{v},\bm{\phi}) + \tau\qty(\mathbcal{B}(\vb{v},\vb{v})-1) = 0,\\
        &\forall\hat{\bm{\phi}}=(\tau,\vb{\phi})\in\mathbb{R}\times\mathbcal{S},
    \end{aligned}
\end{equation}
where the normalisation condition from \cref{eq:normalisation} is enforced weakly. The discrete counterpart of this equation reads:
\begin{equation}
    \begin{aligned}
        &\text{Find}\:\hat{\vb{v}}_h=(\mu_h,\vb{v}_h)\in\mathbb{R}\times\mathbcal{S}_h^p\:\text{s.t.}\\
        &\quad  \mathbcal{V}(\hat{\vb{v}}_h,\hat{\bm{\phi}}_h) = \mu \mathbcal{B}(\vb{v}_h,\bm{\phi}_h) - \mathbcal{A}(\vb{v}_h,\bm{\phi}_h) + \tau_h\qty(\mathbcal{B}(\vb{v}_h,\vb{v}_h)-1) = 0,\\
        &\forall\hat{\bm{\phi}}_h=(\tau_h,\bm{\phi}_h)\in\mathbb{R}\times\mathbcal{S}_h^p.
    \end{aligned}
\end{equation}

Furthermore, a goal-function for the eigenvalues is defined as follows:
\begin{equation}\label{eq:eiggoal}
    \mathbcal{L}(\hat{\vb{v}}) = \mu = \mu \mathbcal{B}(\vb{v},\vb{v}),
\end{equation}
giving
\begin{equation}\label{eq:deiggoal}
    \Delta\mathbcal{L}(\hat{\bm{\phi}}_h) = \mu - \mu_h.
\end{equation}
Using the non-linear functional $\mathbcal{V}$ and the goal functional $\mathbcal{L}$, the same derivations as in \cref{subsec:DWR_general} can be followed to find a system of equations to solve the DWR eigenvalue problem. The Gateaux derivative of $\mathbcal{V}$, denoted by $\mathbcal{V}'$ is given by:
\begin{equation}
    \mathbcal{V}'(\hat{\vb{v}},\hat{\bm{\phi}},\hat{\bm{\psi}}) = \eta \mathbcal{B}(\vb{v},\bm{\psi}) + \mu \mathbcal{B}(\bm{\psi},\bm{\phi}) - \mathbcal{A}(\bm{\psi},\bm{\phi}) + \tau ( \mathbcal{B}(\vb{v},\bm{\psi}) + \mathbcal{B}(\bm{\psi},\vb{v}) ),
\end{equation}
where the derivatives $\mathbcal{A}'(\bm{\psi},\bm{\phi})$ and $\mathbcal{B}'(\bm{\psi},\bm{\phi})$ are equal to the bi-linear operators $\mathbcal{A}(\vb{u},\bm{\phi})$ and $\mathbcal{B}(\vb{u},\bm{\phi})$ themselves. Furthermore, the solution around which the linerisation is performed is denoted by $\hat{\vb{v}}=(\mu,\vb{v})$, the test functions are denoted by $\hat{\bm{\phi}}=(\tau,\bm{\phi})$ and the trial functions are denoted by $\hat{\bm{\psi}}=(\eta,\bm{\psi})$. Furthermore, the linearisation of the goal functional \cref{eq:eiggoal} is
\begin{equation}
    \mathbcal{L}'(\hat{\vb{v}},\hat{\bm{\psi}}) = \eta \mathbcal{B}(\vb{v},\vb{v}) + \mu \qty[ \mathbcal{B}(\vb{v},\bm{\psi}) + \mathbcal{B}(\bm{\psi},\vb{v}) ],
\end{equation}
\RevTWO{
such that the adjoint eigenvalue problem, analoguously to \cref{eq:adjoint_discrete}, given by
\begin{equation}
\text{Find}\:\hat{\bm{\phi}}=(\tau,
    \bm{\phi})\in\mathbb{R}\times\mathbcal{S}\:\text{s.t.}\:
\mathbcal{V}'(\hat{\vb{v}},\hat{\bm{\phi}},\hat{\bm{\psi}}) = \mathbcal{L}'(\hat{\vb{v}},\hat{\bm{\psi}})
\forall\hat{\bm{\psi}}=(\eta,\vb{\psi})\in\mathbb{R}\times\mathbcal{S},
\end{equation}
becomes \cite{Bangerth2003,Lathouwers2011}:
\begin{equation}
    \begin{aligned}
        &\text{Find}\:\hat{\bm{\phi}}=(\tau,
    \bm{\phi})\in\mathbb{R}\times\mathbcal{S}\:\text{s.t.}\\
    &\quad \eta \mathbcal{B}(\vb{v},\bm{\phi}) + \mu \mathbcal{B}(\bm{\psi},\bm{\phi}) - \mathbcal{A}(\bm{\psi},\bm{\phi}) + \tau ( \mathbcal{B}(\vb{v},\bm{\psi}) + \mathbcal{B}(\bm{\psi},\vb{v}) ) = &\eta \mathbcal{B}(\vb{v},\vb{v}) + \\ &\mu \qty[ \mathbcal{B}(\vb{v},\bm{\psi}) + \mathbcal{B}(\bm{\psi},\vb{v}) ],\\
    &\forall\hat{\bm{\psi}}=(\eta,
    \vb{\psi})\in\mathbb{R}\times\mathbcal{S}.
\end{aligned}
\end{equation}
This equation can be simplified to obtain the following \cite{Bangerth2003,Lathouwers2011}:
\begin{equation}\label{eq:eigadjoint2}
    \begin{aligned}
        &\text{Find}\:\hat{\bm{\phi}}=(\tau,
    \bm{\phi})\in\mathbb{R}\times\mathbcal{S}\:\text{s.t.}\\
    &\quad \mu \mathbcal{B}(\bm{\psi},\bm{\phi}) - \mathbcal{A}(\bm{\psi},\bm{\phi}) + \eta \qty[\mathbcal{B}(\vb{v},\bm{\phi}) - \mathbcal{B}(\vb{v},\vb{v}) ] + (\tau - \mu) \qty[ \mathbcal{B}(\vb{v},\bm{\psi}) + \mathbcal{B}(\bm{\psi},\vb{v}) ] = 0, \\
    &\forall\hat{\bm{\psi}}=(\eta,
    \vb{\psi})\in\mathbb{R}\times\mathbcal{S}.
\end{aligned}
\end{equation}
Using the normalizations from \cref{eq:normalisation,eq:normalisation_adj} and the fact that \cref{eq:eigdual} solves the same equation as \cref{eq:eigenvalueproblem}, it follows that \cref{eq:eigadjoint2} is solved by \cref{eq:eigdual} \cite{Bangerth2003}.\\
}
Using \cref{eq:errorest,eq:deiggoal,eq:nlproblem} with $\mathbcal{W}=\mathbcal{V}$ according to \cref{eq:Gfunctional} and with $\bm{\psi}$ denoting the dual eigenvector and $\eta$ the dual eigenvalue, the error estimation according to the DWR method for an eigenvalue problem is
\begin{equation}
    \Delta\mathbcal{L}(\hat{\vb{v}}_h)= \mathbcal{A}(\vb{v}_h,\bm{\psi}-\bm{\psi}_h)- \mu_h \mathbcal{B}(\vb{v}_h,\bm{\psi}-\bm{\psi}_h) + (\eta-\eta_h)(\mathbcal{B}(\vb{v}_h,\vb{v}_h)-1),
\end{equation}
for $\hat{\vb{v}}_h=(\mu_h,\vb{v}_h)\in\mathbb{R}\times\mathbcal{S}_h^p$, $\hat{\bm{\psi}}_h=(\eta_h,\bm{\psi}_h)\in\mathbb{R}\times\mathbcal{S}_h^p$ and $\hat{\bm{\psi}}=(\eta,\bm{\psi})\in\mathbb{R}\times\mathbcal{S}$ . The exact adjoint solution $\hat{\bm{\psi}}_h$ is again approximated by solving \cref{eq:eigdual} on an \emph{enriched} space $\tilde{\mathbcal{S}}_h^p\subset \mathbcal{S}$, $\tilde{\mathbcal{S}}_h^p\supset \mathbcal{S}_h^p$, providing $(\tilde{\eta}_h,\tilde{\bm{\psi}_h})\in\mathbb{R}\times\tilde{\mathbcal{S}}_h^p$. In \cite{Hartmann2003} different choices for constructing $\tilde{\mathbcal{S}}_h^p$ are given, including an $h$-refinement and a $p$-refinement. As in the work of \cite{Hinz2020}, the second approach is used in the present paper, with the same mesh as for $\mathbcal{S}_h^p$, but with a higher order and with the same regularity, i.e. $\tilde{\mathbcal{S}}_h^p = \mathbcal{S}_h^{p+1}$, \RevTWO{as it introduces less degrees of freedom compared to an $h$-refinement.}\\

As specificed in the end of \cref{subsec:structuralAnalysis}, the DWR method for modal analysis requires $\mathbcal{A}(\vb{v},\bm{\phi})=\mathbcal{W}'(\vb{0},\vb{v},\bm{\phi})$ and $\mathbcal{B}(\vb{v},\bm{\phi})=\mathbcal{M}(\vb{v},\bm{\phi})$. For buckling analysis, $\mathbcal{A}(\vb{v},\bm{\phi})=\mathbcal{W}'(\vb{u}_L,\vb{v},\bm{\phi})$ and $\mathbcal{B}(\vb{v},\bm{\phi})=\mathbcal{W}'(\vb{0},\vb{v},\bm{\phi})$ with the first operator defined about a pre-buckling solution $\vb{u}_L$.


%
%
%
%

\subsection{Goal functionals for Isogeometric Kirchhoff-Love Shells}\label{subsec:goalfunctionals}
The remainder of this section focusses on defining the goal functional $\mathbcal{L}(\vb*{u})$, see \cref{eq:goalfunction}, together with its variation $\mathbcal{L}'(\vb*{u})$ such that the dual problem (\cref{eq:adjoint_discrete}) can be solved and the error estimate (\cref{eq:errorest}) can be computed.\\

In general, the goal functional can be defined in a point, on a boundary or over the domain:
\begin{align}
    \mathbcal{L}(\cdot) &= \int_\Omega \mathbcal{l}(\cdot,\vb{x})\dd{\Omega}, && \text{Domain-wise},\\
    \mathbcal{L}(\cdot) &= \int_{\partial\Omega} \mathbcal{l}(\cdot,\vb{x})\dd{\Gamma}, && \text{Boundary-wise},\\
    \mathbcal{L}(\cdot) &= \sum_{i\in\mathbcal{I}} \mathbcal{l}(\cdot,\vb{x}_i), && \text{Point-wise},
\end{align}
where $\Omega$ denotes the integration domain, $\partial\Omega$ a side of $\Omega$ and $\mathbcal{I}$ a set of indices corresponding to points $\vb{x}_i\in\Omega$. Furthermore, $\mathbcal{l}(\cdot,\vb{x}_i)$ denotes a goal functional summant or integrant, which has a variation denoted by $\mathbcal{l}'(\cdot,\bm{\phi}\vb{x}_i)$. The variation of $\mathbcal{L}$, denoted by $\mathbcal{L}'(\cdot,\bm{\phi},\vb{x}_i)$, directly follows from $\mathbcal{l}'(\cdot,\bm{\phi}\vb{x}_i)$ due to linearity of integrals and summation. In addition, we classify two different types of goal functional integrants, resulting in norm-based and component-based goal functionals. In the former case, $\mathbcal{l}$ is of the form $\mathbcal{l}=\Vert\vb{A}\Vert^2$ with variation $\mathbcal{l}'=2\vb{A}\cdot\vb{A}'$. For component-based goal functionals, we define $\mathbcal{l}=\vb{A}\cdot\vb{e}_i$ with variation $\mathbcal{l}'=A'\cdot\vb{e}_i$. Here, $\vb{e}_i$ is a unit vector in direction $i$. \RevTWO{It should be noted that the goal functional $\mathbcal{L}(\cdot)$ needs to be bounded in all cases, therefore making point-wise goal functionals not always suitable, e.g. in case of a stress singularity in a point.}\\

In \cref{tab:overviewGFs} we provide some goal functional integrants or summants $\mathbcal{l}(\cdot,\vb{x}_i)$. Together with their variations $\mathbcal{l}'(\cdot,\bm{\phi}\vb{x}_i)$, these provide $\mathbcal{L}'(\cdot,\bm{\phi},\vb{x}_i)$ due to linearity of integrals and summation. The \emph{tensor-based} goal functionals refer to goal functionals that could be used for any second-order tensor, e.g. the membrane strain tensor $\bm{\varepsilon}(\vb{u})$ or the flexural moment tensor $\vb{m}(\vb{u})$.\\

\begin{table}
    \centering
    \caption{Overview of the goal functionals. Here $\vb{u}_h=\vb{u}(\vb{x})$ is the discrete deformation tensor depending on position coordinate $\vb{x}$, $\vb{C}_h=\vb{C}(\vb{u}_h)$ is the deformation tensor based on $\vb{u}_h$ and $\vb{C}(\bm{\phi}'_h)=\vb{C}'(\vb{u}_h,\bm{\phi})$ is its variation, $\mathbcal{T}(\vb{A})$ is the transformation of a second-order tensor $\vb{A}$ from the undeformed contravariant basis to the basis spanned by the principal directions. Note that the variation $[\mathbcal{T}(\vb{A})]'$ of $\mathbcal{T}(\vb{A})$ is $\mathbcal{T}(\vb{A}')$, since the spectral decomposition of the deformation tensor itself is just a linear change of tensor basis.}
    \label{tab:overviewGFs}
    \begin{tabular}{@{}lll@{}}
        \toprule
        Displacement-norm & $\mathbcal{l}_{\Vert\vb{u}\Vert}(\vb{u}_h) = \Vert \vb{u}_h(\vb{x}) \Vert^2$ & $\mathbcal{l}'_{\Vert\vb{u}\Vert}(\vb{u}_h,\bm{\phi}_h) = 2 \vb{u}_h\cdot\bm{\phi}_h$\\
        Displacement-component & $\mathbcal{l}_{\vb{u}_i}(\vb{u}_h) = \vb{u}_h(\vb{x})\cdot \vb{e}_i$ & $\mathbcal{l}'_{\Vert\vb{u}\Vert}(\vb{u}_h,\bm{\phi}_h) = \vb{e}_i\cdot\bm{\phi}_h$\\
        \midrule
        Stretch-norm & $\mathbcal{l}_{\Vert\bm{\lambda}\Vert}(\vb{u}_h) = \Vert \mathbcal{T}(\vb{C}_h) \Vert^2$ & $\mathbcal{l}'_{\Vert\bm{\lambda}\Vert}(\vb{u}_h,\bm{\phi}_h) = 2 \mathbcal{T}(\vb{C}_h)\cdot\bm{\phi}_h$\\
        Stretch-component & $\mathbcal{l}_{\bm{\lambda}_i}(\vb{u}_h) = \mathbcal{T}(\vb{C}_h)\cdot \vb{e}_i$ & $\mathbcal{l}'_{\Vert\bm{\lambda}\Vert}(\vb{u}_h,\bm{\phi}_h) = \vb{e}_i\cdot\mathbcal{T}(\vb{C}'_h)$\\
        \midrule
        Tensor-norm & $\mathbcal{l}_{\Vert\vb{A}\Vert}(\vb{u}_h) = \Vert \vb{A} \Vert^2$ & $\mathbcal{l}'_{\Vert\vb{A}\Vert}(\vb{u}_h,\bm{\phi}_h) = 2 \vb{A}\cdot\bm{\phi}_h$\\
        Tensor-component & $\mathbcal{l}_{\vb{A}_i}(\vb{u}_h) = \vb{A}\cdot \vb{e}_i$ & $\mathbcal{l}'_{\Vert\vb{A}\Vert}(\vb{u}_h,\bm{\phi}_h) = \vb{e}_i\cdot\vb{A}$\\
        \botrule
    \end{tabular}
\end{table}

\section{Coarsening and Refinement using THB Splines}\label{sec:refinement}
This section elaborates on coarsening and refinement of isogeometric meshes using THB-splines. In particular, this section elaborates on the \emph{Mark}, \emph{Refine} and \emph{Transfer} blocks of \cref{fig:flowchart}. Firstly, \cref{subsec:THB} will provide a brief background on THB-splines, which enable the \emph{Refine} step of the adaptive meshing flowchart. Then, \cref{subsec:grading} elaborates on methods for suitable grading for refinement meshes; which is required to provide admissible refinement with (Truncated) Hierarchical B-spline ((T)HB) bases, given labelled elements. \Cref{subsec:labeling} elaborates on the labeling method for the \emph{Mark} step, given an element-wise error distribution, taking admissibility into account. Lastly, \cref{subsec:interpolate} elaborates on the quasi-interpolation method that is used to \emph{Transfer} the solution of one solution step to the next. The notations in this section will be closely related to those used in \cite{Bracco2018a,Carraturo2019}.\\

\subsection{(T)HB-Splines}\label{subsec:THB}
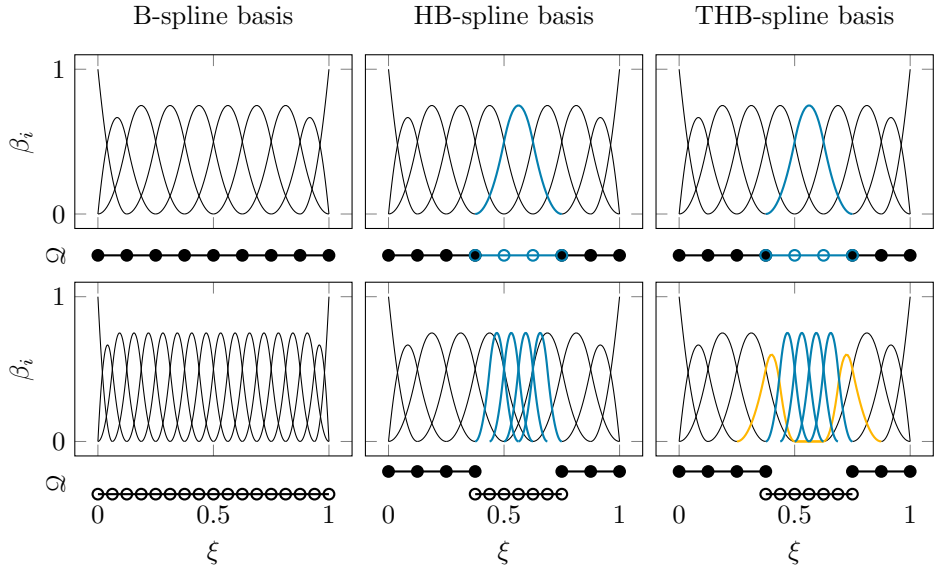
\begin{figure}
\centering
 \begin{tikzpicture}
\begin{groupplot}
[
height=0.2\textheight,
width=0.4\linewidth,
ytick = {0,1}, 
legend pos = outer north east,
group style={
                group name=my plots,
                group size=3 by 4,
                xlabels at=edge bottom,
                xticklabels at=edge bottom,
                ylabels at=edge left,
                yticklabels at=edge left,
                vertical sep=5pt,
                horizontal sep=5pt},
ylabel={$\beta_i$},
xlabel={$\xi$},
]
\nextgroupplot[no markers,title=B-spline basis,]
\foreach \k in {0,...,9}
{
\pgfmathtruncatemacro{\y}{2*\k+1};
\pgfmathtruncatemacro{\x}{2*\k};
\addplot+[black,thin,solid] table[col sep = comma, header=false, y index = {\y}, x index = {\x}] {Data/Refinement_illustration/init.csv};
}

\nextgroupplot[no markers,title=HB-spline basis,]
\foreach \k in {0,...,4}
{
\pgfmathtruncatemacro{\y}{2*\k+1};
\pgfmathtruncatemacro{\x}{2*\k};
\addplot+[black,thin,solid] table[col sep = comma, header=false, y index = {\y}, x index = {\x}] {Data/Refinement_illustration/init.csv};
}

\foreach \k in {6,...,9}
{
\pgfmathtruncatemacro{\y}{2*\k+1};
\pgfmathtruncatemacro{\x}{2*\k};
\addplot+[black,thin,solid] table[col sep = comma, header=false, y index = {\y}, x index = {\x}] {Data/Refinement_illustration/init.csv};
}
\def\k{5}
\pgfmathtruncatemacro{\y}{2*\k+1};
\pgfmathtruncatemacro{\x}{2*\k};
\addplot+[col1,thick,solid] table[col sep = comma, header=false, y index = {\y}, x index = {\x}] {Data/Refinement_illustration/init.csv};

\nextgroupplot[no markers,title=THB-spline basis,]
\foreach \k in {0,...,4}
{
\pgfmathtruncatemacro{\y}{2*\k+1};
\pgfmathtruncatemacro{\x}{2*\k};
\addplot+[black,thin,solid] table[col sep = comma, header=false, y index = {\y}, x index = {\x}] {Data/Refinement_illustration/init.csv};
}

\foreach \k in {6,...,9}{
\pgfmathtruncatemacro{\y}{2*\k+1};
\pgfmathtruncatemacro{\x}{2*\k};
\addplot+[black,thin,solid] table[col sep = comma, header=false, y index = {\y}, x index = {\x}] {Data/Refinement_illustration/init.csv};
}
\def\k{5}
\pgfmathtruncatemacro{\y}{2*\k+1};
\pgfmathtruncatemacro{\x}{2*\k};
\addplot+[col1,thick,solid] table[col sep = comma, header=false, y index = {\y}, x index = {\x}] {Data/Refinement_illustration/init.csv};

\nextgroupplot
[
height=0.1\textheight,      
axis line style={draw=none},
tick style={draw=none},
ylabel={$\mathcal{Q}$},
ytick=\empty,
ymin=0.9,ymax=1.1
]
\addplot[mark=*,black] coordinates {(0,1) (0.125,1) (0.250,1) (0.375,1) (0.500,1) (0.625,1) (0.750,1) (0.875,1) (1.000,1)};

\nextgroupplot
[
height=0.1\textheight,      
axis line style={draw=none},
tick style={draw=none},
ymin=0.9,ymax=1.1
]
\addplot[mark=*,black] coordinates {(0,1) (0.125,1) (0.250,1) (0.375,1)};
\addplot[mark=*,black] coordinates {(0.750,1) (0.875,1) (1.000,1)};
\addplot[mark=o,col1] coordinates { (0.375,1) (0.500,1) (0.625,1) (0.750,1)};

\nextgroupplot
[
height=0.1\textheight,      
axis line style={draw=none},
tick style={draw=none},
ymin=0.9,ymax=1.1
]
\addplot[mark=*,black] coordinates {(0,1) (0.125,1) (0.250,1) (0.375,1)};
\addplot[mark=*,black] coordinates {(0.750,1) (0.875,1) (1.000,1)};
\addplot[mark=o,col1] coordinates { (0.375,1) (0.500,1) (0.625,1) (0.750,1)};

\nextgroupplot[no markers]
\foreach \k in {0,...,17}
{
\pgfmathtruncatemacro{\y}{2*\k+1};
\pgfmathtruncatemacro{\x}{2*\k};
\addplot+[black,thin,solid] table[col sep = comma, header=false, y index = {\y}, x index = {\x}] {Data/Refinement_illustration/BSpline.csv};
}

\nextgroupplot[no markers]
\foreach \k in {0,...,8}
{
\pgfmathtruncatemacro{\y}{2*\k+1};
\pgfmathtruncatemacro{\x}{2*\k};
\addplot+[black,thin,solid] table[col sep = comma, header=false, y index = {\y}, x index = {\x}] {Data/Refinement_illustration/HBspline.csv};
}

\foreach \k in {9,...,12}
{
\pgfmathtruncatemacro{\y}{2*\k+1};
\pgfmathtruncatemacro{\x}{2*\k};
\addplot+[col1,thick,solid] table[col sep = comma, header=false, y index = {\y}, x index = {\x}] {Data/Refinement_illustration/HBspline.csv};
}

\nextgroupplot[no markers]
\foreach \k in {0,1,2,3,6,7,8}
{
\pgfmathtruncatemacro{\y}{2*\k+1};
\pgfmathtruncatemacro{\x}{2*\k};
\addplot+[black,thin,solid] table[col sep = comma, header=false, y index = {\y}, x index = {\x}] {Data/Refinement_illustration/THBspline.csv};
}

\foreach \k in {4,5}
{
\pgfmathtruncatemacro{\y}{2*\k+1};
\pgfmathtruncatemacro{\x}{2*\k};
\addplot+[col2,thick,solid] table[col sep = comma, header=false, y index = {\y}, x index = {\x}] {Data/Refinement_illustration/THBspline.csv};
}

\foreach \k in {9,...,12}
{
\pgfmathtruncatemacro{\y}{2*\k+1};
\pgfmathtruncatemacro{\x}{2*\k};
\addplot+[col1,thick,solid] table[col sep = comma, header=false, y index = {\y}, x index = {\x}] {Data/Refinement_illustration/THBspline.csv};
}

\nextgroupplot
[
height=0.1\textheight,      
axis line style={draw=none},
tick style={draw=none},
ylabel={$\mathcal{Q}$},
ytick=\empty,
ymin=-0.1,ymax=1.1
]
\addplot[mark=o,black] coordinates {(0,0) (0.0625,0) (0.125,0) (0.1875,0) (0.250,0) (0.3125,0) (0.375,0) (0.4375,0) (0.500,0) (0.5625,0) (0.625,0) (0.6875,0) (0.750,0) (0.8125,0) (0.875,0) (0.9375,0) (1.000,0)};

\nextgroupplot
[
height=0.1\textheight,      
axis line style={draw=none},
tick style={draw=none},
ymin=-0.1,ymax=1.1
]
\addplot[mark=*,black] coordinates {(0,1) (0.125,1) (0.250,1) (0.375,1)};
\addplot[mark=*,black] coordinates {(0.750,1) (0.875,1) (1.000,1)};
\addplot[mark=o,black] coordinates {(0.375,0) (0.4375,0) (0.500,0) (0.5625,0) (0.625,0) (0.6875,0) (0.750,0)};

\nextgroupplot
[
height=0.1\textheight,      
axis line style={draw=none},
tick style={draw=none},
ymin=-0.1,ymax=1.1
]
\addplot[mark=*,black] coordinates {(0,1) (0.125,1) (0.250,1) (0.375,1)};
\addplot[mark=*,black] coordinates {(0.750,1) (0.875,1) (1.000,1)};
\addplot[mark=o,black] coordinates {(0.375,0) (0.4375,0) (0.500,0) (0.5625,0) (0.625,0) (0.6875,0) (0.750,0)};
\end{groupplot}
\end{tikzpicture}
\caption{Principles of refinement for different spline bases. The top figures represent the basis on level $\mathbcal{B}^0$, optionally with refined basis functions coloured \textcolor{col1}{blue}. Bottom pictures represent refined bases. left) uniform refinement (hence $\mathbcal{B}^1$; middle) HB-refinement; right) THB-refinement, with truncated basis functions coloured \textcolor{col2}{yellow}. Note that the refinement basis functions are from $V^1$. The unrefined unique knot vector in all cases is $\Xi = \{0,1/8,2/8,\dots,7/8,1\}$ and the degree of the basis is 2. The bases are generated in \texttt{G+Smo} \cite{Juttler2014}. }
\label{fig:IGArefinement}
\end{figure}

Refinement of B-spline meshes can be done using (Truncated) Hierarchical B-splines ((T)HB-splines), of which the details can be found in \cite{Vuong2011,Giannelli2012}. The conceptual idea behind (T)HB-splines is that they are constructed from a sequence of $N$ nested tensor B-spline spaces in different levels $l=0,...,N-1$, denoted by $V^0\subset V^1\subset,...,V^{N-1}$ with associated bases $\mathbcal{B}^\ell$ of degree $p$ on a grid $G^\ell$ with elements $Q$. The parametric domains are defined as $\Omega=\Omega^0\supseteq\Omega^1\supseteq...\supseteq\Omega^{N-1}=\emptyset$. By defining the set of active cells by $\mathbcal{G}^\ell:=\{Q\in G^\ell: Q\subset\Omega^\ell \wedge Q\not\subset \Omega^{\ell+1}\}$, the hierarchical mesh is defined as $\mathbcal{Q}=\{Q\in \mathbcal{G}^\ell:\ell=0,...,N-1\}$.  In \cref{fig:IGArefinement}, an illustration is given for a refined B-spline basis (left), a refined HB-spline basis (middle) and a refined THB-spline basis (right). For the (T)HB-spline basis, this picture depicts the refinement of a single basis function, corresponding to the elements in its support. The (T)HB-spline bases show that for THB-splines a truncation is performed to ensure partition of unity, which is discussed in more detail in \cite{Giannelli2012}.

\subsection{Admissible Meshing}\label{subsec:grading}
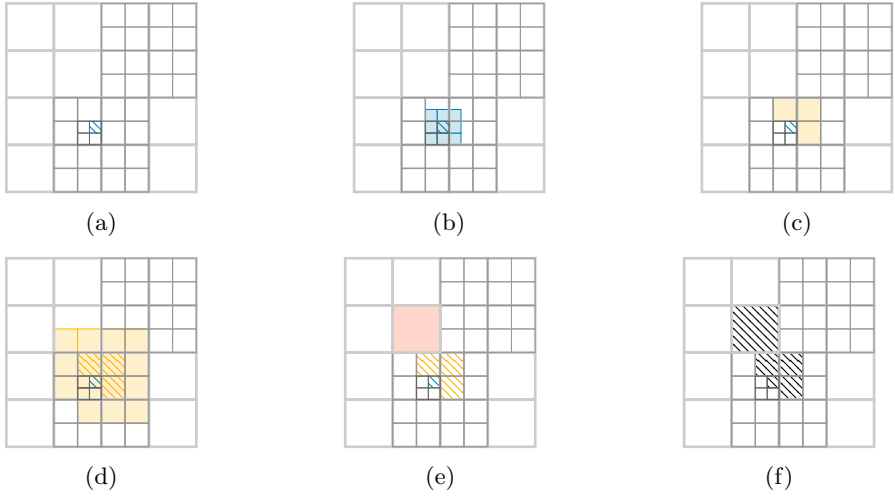
\begin{figure}
 \centering
 {
\newcommand\Grid[7]{
	\pgfmathsetmacro\h{2^#1}
	\pgfmathsetmacro\imin{#2}
	\pgfmathsetmacro\jmin{#3}
	\pgfmathsetmacro\imax{#4-1}
	\pgfmathsetmacro\jmax{#5-1}
	\pgfmathsetmacro\level{#6}
	\pgfmathsetmacro\linv{1/#1}
	\foreach \i in {\imin,...,\imax}
	{
		\pgfmathsetmacro\ipp{\i+1}
		\foreach \j in {\jmin,...,\jmax}
		{
			\pgfmathsetmacro\jpp{\j+1}
			\draw[#7] (\i*\level/\h,\j*\L/\h) rectangle (\ipp*\level/\h,\jpp*\level/\h);
		}
	}
}

\def\L{5}

\begin{subfigure}{0.3\linewidth}
\centering
\begin{tikzpicture}[scale=0.5]
\Grid{2}{0}{0}{4}{4}{\L}{black!20,line width=1}
\Grid{3}{2}{0}{6}{4}{\L}{black!40,line width=1/2}
\Grid{3}{4}{4}{8}{8}{\L}{black!40,line width=1/3}
\Grid{4}{6}{4}{8}{6}{\L}{black!60,line width=1/4}
\Grid{4}{7}{5}{8}{6}{\L}{draw=none,line width=1/4,pattern={north west lines},pattern color=col1}
\end{tikzpicture}
\caption{}
\end{subfigure}
\hfill
\begin{subfigure}{0.3\linewidth}
\centering
\begin{tikzpicture}[scale=0.5]
\Grid{4}{6}{4}{9}{7}{\L}{fill=col1!20,draw=col1,line width=1/5}

\Grid{2}{0}{0}{4}{4}{\L}{black!20,line width=1}
\Grid{3}{2}{0}{6}{4}{\L}{black!40,line width=1/2}
\Grid{3}{4}{4}{8}{8}{\L}{black!40,line width=1/3}
\Grid{4}{6}{4}{8}{6}{\L}{black!60,line width=1/4}
\Grid{4}{7}{5}{8}{6}{\L}{draw=none,line width=1/4,pattern={north west lines},pattern color=col1}
\end{tikzpicture}
\caption{}
\end{subfigure}
\hfill
\begin{subfigure}{0.3\linewidth}
\centering
 \begin{tikzpicture}[scale=0.5]
\Grid{3}{3}{3}{5}{4}{\L}{fill=col2!20,draw=col2,line width=1/5}
\Grid{3}{4}{2}{5}{3}{\L}{fill=col2!20,draw=col2,line width=1/5}

\Grid{2}{0}{0}{4}{4}{\L}{black!20,line width=1}
\Grid{3}{2}{0}{6}{4}{\L}{black!40,line width=1/2}
\Grid{3}{4}{4}{8}{8}{\L}{black!40,line width=1/3}
\Grid{4}{6}{4}{8}{6}{\L}{black!60,line width=1/4}
\Grid{4}{7}{5}{8}{6}{\L}{draw=none,line width=1/4,pattern={north west lines},pattern color=col1}
\end{tikzpicture}
\caption{}
\end{subfigure}

\begin{subfigure}{0.3\linewidth}
\centering
 \begin{tikzpicture}[scale=0.5]
\Grid{3}{3}{3}{5}{4}{\L}{draw=black!60,line width=1/5,pattern={north west lines},preaction={fill,  col2!20},pattern color=col2}
\Grid{3}{4}{2}{5}{3}{\L}{draw=black!60,line width=1/5,pattern={north west lines},preaction={fill, col2!20},pattern color=col2}
\Grid{4}{7}{5}{8}{6}{\L}{draw=black!80,line width=1/4,pattern={north west lines},preaction={fill, col2!20},pattern color=col1}
\Grid{4}{6}{5}{7}{6}{\L}{draw=black!80,line width=1/4,fill=col2!20}
\Grid{4}{6}{4}{8}{5}{\L}{draw=black!80,line width=1/4,fill=col2!20}
\Grid{3}{3}{1}{6}{2}{\L}{fill=col2!20,draw=col2,line width=1/5}
\Grid{3}{5}{2}{6}{5}{\L}{fill=col2!20,draw=col2,line width=1/5}
\Grid{3}{2}{4}{6}{5}{\L}{fill=col2!20,draw=col2,line width=1/5}
\Grid{3}{2}{2}{3}{4}{\L}{fill=col2!20,draw=col2,line width=1/5}

\Grid{2}{0}{0}{4}{4}{\L}{black!20,line width=1}
\Grid{3}{2}{0}{6}{4}{\L}{black!40,line width=1/2}
\Grid{3}{4}{4}{8}{8}{\L}{black!40,line width=1/3}
\Grid{4}{6}{4}{8}{6}{\L}{black!60,line width=1/4}
\end{tikzpicture}
\caption{}
\end{subfigure}
\hfill
\begin{subfigure}{0.3\linewidth}
\centering
\begin{tikzpicture}[scale=0.5]
\Grid{3}{3}{3}{5}{4}{\L}{draw=none,line width=1/5,pattern={north west lines},pattern color=col2}
\Grid{3}{4}{2}{5}{3}{\L}{draw=none,line width=1/5,pattern={north west lines},pattern color=col2}
\Grid{4}{7}{5}{8}{6}{\L}{draw=none,line width=1/4,pattern={north west lines},pattern color=col1}

\Grid{2}{0}{0}{4}{4}{\L}{black!20,line width=1}
\Grid{3}{2}{0}{6}{4}{\L}{black!40,line width=1/2}
\Grid{3}{4}{4}{8}{8}{\L}{black!40,line width=1/3}
\Grid{4}{6}{4}{8}{6}{\L}{black!60,line width=1/4}
\Grid{2}{1}{2}{2}{3}{\L}{draw=black!20,fill=col3!20,line width=1}
\end{tikzpicture}
\caption{}
\end{subfigure}
\hfill
\begin{subfigure}{0.3\linewidth}
\centering
 \begin{tikzpicture}[scale=0.5]
\Grid{3}{3}{3}{5}{4}{\L}{draw=black!60,line width=1/5,pattern={north west lines},pattern color=black}
\Grid{3}{4}{2}{5}{3}{\L}{draw=black!60,line width=1/5,pattern={north west lines},pattern color=black}
\Grid{4}{7}{5}{8}{6}{\L}{draw=black!80,line width=1/4,pattern={north west lines},pattern color=black}

\Grid{2}{0}{0}{4}{4}{\L}{black!20,line width=1}
\Grid{3}{2}{0}{6}{4}{\L}{black!40,line width=1/2}
\Grid{3}{4}{4}{8}{8}{\L}{black!40,line width=1/3}
\Grid{4}{6}{4}{8}{6}{\L}{black!60,line width=1/4}
\Grid{2}{1}{2}{2}{3}{\L}{draw=black!20,line width=1,pattern={north west lines},pattern color=black}
\end{tikzpicture}
\caption{}
\end{subfigure}
}
 \caption{Recursive marking strategy on the marked element of level $\ell$ on the initial mesh represented in (a). As a first step, the support extension of the marked element is obtained (b), from which the parents that are active on level $\ell-1$ define the $\mathcal{T}$-neighborhood of the marked cell (c). Starting the same procedure on the marked cells of level $\ell-1$, the support extension can again be obtained (d) with their corresponding parents on level $\ell-2$, marking the $\mathcal{T}$-neighborhood of the marked elements of level $\ell-1$ (e). The complete recursive marking from the marked element in (a) is depicted in (f).}
 \label{fig:Neighborhoods_ref}
\end{figure}

The concept of admissible meshing was discussed in \cite{Buffa2016,Bracco2018a,Buffa2022}. An admissible mesh of class $m$ is a mesh of which the truncated basis functions belong to at most $m$ successive levels and mesh admissibility ensures that the number of basis functions acting on a mesh elements does not depend on the number of levels in the hierarchy, but on the parameter $m$. In order to guarantee mesh admissibility for refinement and coarsening operations, refinement and coarsening neighborhoods are defined such that admissible meshes can be constructed recursively, which is discussed in more detail in \cite{Buffa2016,Bracco2018a,Buffa2022}. In \cref{fig:Neighborhoods_ref}, we illustrate a simple mesh together with the refinement neighborhood of some selected elements. The $\mathcal{T}$-refinement neighborhood $\mathbcal{N}_r(\mathbcal{Q},Q,m)$ of element $Q$ is defined as
\begin{equation}
 \mathbcal{N}_r^\text{THB}(\mathbcal{Q},Q,m) = \left\{Q'\in\mathbcal{G}^{\ell-m+1}: \exists Q''\in S(Q,\ell-m+2),Q''\subseteq Q'\right\},
\end{equation}
where $S(Q,k)$ is the \emph{multi-level} support extension with respect to level $k$.\\

The coarsening neighborhood $\mathbcal{N}_c(Q)$ of element $Q\in G^\ell$ is defined by \cite{Carraturo2019}. When coarsening element $Q^\ell$, the coarsening neighborhood ensures that the newly activated basis functions are not active on the surrounding basis functions of level $\ell+m$. In other words, if element $Q$ of level $\ell$ is the element to be coarsened, then the coarsening neighborhood defined by
\begin{equation}\label{eq:crsNeighborhood}
 \mathbcal{N}_c(\mathbcal{Q},Q,m):=\left\{Q'\in\mathbcal{G}^{\ell+m-1}:\exists Q''\in\mathbcal{G}^{\ell}\text{ and }Q''\subset P(Q),\text{ with }Q'\in S(Q'',\ell)\right\},
\end{equation}
must be empty. Here, $P(Q)$ denotes the parent of $Q$, i.e. the unique cell $Q'\in G^{\ell-1}$ such that $Q\subset Q'$. Note the small difference with respect to the definintion given in \cite{Carraturo2019}, since the coarsening neighborhood in their work is defined for the element $\hat{Q}$ of level $\ell$ which will be activated, i.e. $\hat{Q}$ is the parent of $Q$ for which the coarsening neighborhood is defined here. Given the definition in \cref{eq:crsNeighborhood} and given a set of elements marked for refinement $\mathbcal{M}_r$, a coarsening neighborhood checking elements marked for refinement, can be defined:
\begin{multline}\label{eq:crsNeighborhood2}
 \mathbcal{N}^r_c(\mathbcal{Q},Q,m,\mathbcal{M}_r):=\big\{Q'\in\mathbcal{G}^{\ell+m-2}:Q'\in\mathbcal{M}_r,\exists Q''\in\mathbcal{G}^{\ell-1}\text{ and }Q''\subset P(Q),\\ \text{ with }Q'\in S(Q'',\ell-1)\big\}.
\end{multline}
In other words, this is the coarsening neighborhood that checks whether for element $Q$ of level $\ell$ to be coarsened there are elements in the marked set $\mathbcal{M}_r$ that will be part of the coarsening neighborhoord as soon as they are refined; thus it uses \cref{eq:crsNeighborhood} with $\ell-1$. This neighborhood ensures that coarse labeling can be performed conforming with the Dörfler marking strategy and without refining first. This avoids to compute element-error contributions on an in-between mesh which has been refined first. Obviously, if another element with the same parent as $Q$ is marked for refinement, no coarsening should take place. An element can be coarsened if $ \mathbcal{N}^r_c(\mathbcal{Q},Q,m,\mathbcal{M}_r)=\emptyset$. Combining both neighborhoods, an element $Q$ of level $\ell$ can be coarsened if and only if $\hat{\mathbcal{N}}_c(\mathbcal{Q},Q,m,\mathbcal{M}_r)=\mathbcal{N}_c(\mathbcal{Q},Q,m)\cup\mathbcal{N}_c^r(\mathbcal{Q},Q,m,\mathbcal{M}_r)=\emptyset$. In \cref{fig:neighborhood_crs}, the coarsening neighborhood is illustrated for a simple mesh.
\begin{figure}
 \centering
 {
\newcommand\Grid[7]{
	\pgfmathsetmacro\h{2^#1}
	\pgfmathsetmacro\imin{#2}
	\pgfmathsetmacro\jmin{#3}
	\pgfmathsetmacro\imax{#4-1}
	\pgfmathsetmacro\jmax{#5-1}
	\pgfmathsetmacro\level{#6}
	\pgfmathsetmacro\linv{1/#1}
	\foreach \i in {\imin,...,\imax}
	{
		\pgfmathsetmacro\ipp{\i+1}
		\foreach \j in {\jmin,...,\jmax}
		{
			\pgfmathsetmacro\jpp{\j+1}
			\draw[#7] (\i*\level/\h,\j*\L/\h) rectangle (\ipp*\level/\h,\jpp*\level/\h);
		}
	}
}

\def\L{5}

\begin{subfigure}{0.3\linewidth}
\centering
 \begin{tikzpicture}[scale=0.5]
\Grid{3}{3}{6}{7}{7}{\L}{fill=col2!20,draw=black,line width=1/2}
\Grid{3}{3}{3}{7}{4}{\L}{fill=col2!20,draw=black,line width=1/2}
\Grid{3}{3}{4}{4}{6}{\L}{fill=col2!20,draw=black,line width=1/2}
\Grid{3}{6}{4}{7}{6}{\L}{fill=col2!20,draw=black,line width=1/2}
\Grid{3}{5}{5}{6}{6}{\L}{line width=1/4,pattern={north east lines},pattern color=col2}

\Grid{3}{3}{3}{5}{4}{\L}{draw=black,line width=1/5,pattern={north west lines},pattern color=black}
\Grid{3}{4}{2}{5}{3}{\L}{draw=black,line width=1/5,pattern={north west lines},pattern color=black}
\Grid{4}{7}{5}{8}{6}{\L}{draw=none,line width=1/4,pattern={north west lines},pattern color=black}

\Grid{2}{0}{0}{4}{4}{\L}{black!20,line width=1}
\Grid{3}{2}{0}{6}{4}{\L}{black!40,line width=1/2}
\Grid{3}{4}{4}{8}{8}{\L}{black!40,line width=1/3}
\Grid{4}{6}{4}{8}{6}{\L}{black!60,line width=1/4}
\Grid{2}{1}{2}{2}{3}{\L}{draw=none,line width=1/4,pattern={north west lines},pattern color=black}
\end{tikzpicture}
\caption{}
\end{subfigure}
\hfill
\begin{subfigure}{0.3\linewidth}
\centering
 \begin{tikzpicture}[scale=0.5]
\Grid{3}{5}{5}{8}{6}{\L}{fill=col2!20,draw=black,line width=1/2}
\Grid{3}{4}{5}{5}{8}{\L}{fill=col2!20,draw=black,line width=1/2}
\Grid{3}{7}{6}{8}{8}{\L}{fill=col2!20,draw=black,line width=1/2}
\Grid{3}{6}{6}{7}{7}{\L}{line width=1/4,pattern={north east lines},pattern color=col2}

\Grid{3}{3}{3}{5}{4}{\L}{draw=black,line width=1/5,pattern={north west lines},pattern color=black}
\Grid{3}{4}{2}{5}{3}{\L}{draw=black,line width=1/5,pattern={north west lines},pattern color=black}
\Grid{4}{7}{5}{8}{6}{\L}{draw=none,line width=1/4,pattern={north west lines},pattern color=black}

\Grid{2}{0}{0}{4}{4}{\L}{black!20,line width=1}
\Grid{3}{2}{0}{6}{4}{\L}{black!40,line width=1/2}
\Grid{3}{4}{4}{8}{8}{\L}{black!40,line width=1/3}
\Grid{4}{6}{4}{8}{6}{\L}{black!60,line width=1/4}
\Grid{2}{1}{2}{2}{3}{\L}{draw=none,line width=1/4,pattern={north west lines},pattern color=black}
\end{tikzpicture}
\caption{}
\end{subfigure}
\hfill
\begin{subfigure}{0.3\linewidth}
\centering
 \begin{tikzpicture}[scale=0.5]
\Grid{3}{5}{6}{8}{7}{\L}{fill=col2!20,draw=black,line width=1/2}
\Grid{3}{5}{3}{8}{4}{\L}{fill=col2!20,draw=black,line width=1/2}
\Grid{3}{5}{4}{6}{6}{\L}{fill=col2!20,draw=black,line width=1/2}
\Grid{3}{6}{4}{7}{5}{\L}{line width=1/4,pattern={north east lines},pattern color=col2}

\Grid{3}{3}{3}{5}{4}{\L}{draw=black,line width=1/5,pattern={north west lines},pattern color=black}
\Grid{3}{4}{2}{5}{3}{\L}{draw=black,line width=1/5,pattern={north west lines},pattern color=black}
\Grid{4}{7}{5}{8}{6}{\L}{draw=none,line width=1/4,pattern={north west lines},pattern color=black}

\Grid{2}{0}{0}{4}{4}{\L}{black!20,line width=1}
\Grid{3}{2}{0}{6}{4}{\L}{black!40,line width=1/2}
\Grid{3}{4}{4}{8}{8}{\L}{black!40,line width=1/3}
\Grid{4}{6}{4}{8}{6}{\L}{black!60,line width=1/4}
\Grid{2}{1}{2}{2}{3}{\L}{draw=none,line width=1/4,pattern={north west lines},pattern color=black}
\end{tikzpicture}
\caption{}
\end{subfigure}

\begin{subfigure}{0.3\linewidth}
\centering
\begin{tikzpicture}[scale=0.5]
\Grid{3}{4}{4}{8}{8}{\L}{fill=col2!20,draw=none,line width=1/2}
\Grid{4}{6}{4}{8}{6}{\L}{fill=col3!20,draw=none,line width=1/3}

\Grid{3}{3}{3}{5}{4}{\L}{draw=black,line width=1/5,pattern={north west lines},pattern color=black}
\Grid{3}{4}{2}{5}{3}{\L}{draw=black,line width=1/5,pattern={north west lines},pattern color=black}
\Grid{4}{7}{5}{8}{6}{\L}{draw=none,line width=1/4,pattern={north west lines},pattern color=black}
\Grid{2}{1}{2}{2}{3}{\L}{draw=none,line width=1/4,pattern={north west lines},pattern color=black}

\Grid{2}{0}{0}{4}{4}{\L}{black!20,line width=1}
\Grid{3}{2}{0}{6}{4}{\L}{black!40,line width=1/2}
\Grid{3}{4}{4}{8}{8}{\L}{black!40,line width=1/3}
\Grid{4}{6}{4}{8}{6}{\L}{black!60,line width=1/4}
\end{tikzpicture}
\caption{} 
\end{subfigure}
\hfill
\begin{subfigure}{0.3\linewidth}
\centering
\begin{tikzpicture}[scale=0.5]
\Grid{3}{4}{6}{8}{8}{\L}{fill=col2!20,draw=none,line width=1/2}
\Grid{3}{6}{4}{8}{6}{\L}{fill=col2!20,draw=none,line width=1/2}

\Grid{3}{6}{4}{8}{6}{\L}{fill=col2!20,draw=none,line width=1/2}

\Grid{3}{3}{3}{5}{4}{\L}{draw=none,line width=1/5,pattern={north west lines},pattern color=black}
\Grid{3}{4}{2}{5}{3}{\L}{draw=none,line width=1/5,pattern={north west lines},pattern color=black}
\Grid{4}{7}{5}{8}{6}{\L}{draw=none,line width=1/4,pattern={north west lines},pattern color=black}
\Grid{2}{1}{2}{2}{3}{\L}{draw=none,line width=1/4,pattern={north west lines},pattern color=black}

\Grid{2}{0}{0}{4}{4}{\L}{black!20,line width=1}
\Grid{3}{2}{0}{6}{4}{\L}{black!40,line width=1/2}
\Grid{3}{4}{4}{8}{8}{\L}{black!40,line width=1/3}
\Grid{4}{6}{4}{8}{6}{\L}{black!60,line width=1/4}
\end{tikzpicture}
\caption{}
\end{subfigure}
\hfill
\begin{subfigure}{0.3\linewidth}
\centering
 \begin{tikzpicture}[scale=0.5]
\Grid{2}{2}{3}{4}{4}{\L}{line width=1/4,pattern={north east lines},pattern color=black!10}
\Grid{2}{3}{2}{4}{3}{\L}{line width=1/4,pattern={north east lines},pattern color=black!10}
 
\Grid{3}{3}{3}{5}{4}{\L}{draw=none,line width=1/5,pattern={north west lines},pattern color=black!10}
\Grid{3}{4}{2}{5}{3}{\L}{draw=none,line width=1/5,pattern={north west lines},pattern color=black!10}
\Grid{4}{7}{5}{8}{6}{\L}{draw=none,line width=1/4,pattern={north west lines},pattern color=black!10}
\Grid{2}{1}{2}{2}{3}{\L}{draw=none,line width=1/4,pattern={north west lines},pattern color=black!10}

\Grid{2}{0}{0}{4}{4}{\L}{black!20,line width=1}
\Grid{3}{2}{0}{6}{6}{\L}{black!40,line width=1/2}
\Grid{4}{6}{4}{10}{8}{\L}{black!60,line width=1/4}
\Grid{5}{14}{10}{16}{12}{\L}{black!80,line width=1/5}
\end{tikzpicture}
\caption{}
\end{subfigure}
}
 \caption{Given the mesh from \cref{fig:Neighborhoods_ref}-(f), the coarsening neighborhoods are evaluated in (a)-(c) of this figure. The cell \protect\tikz{\protect\draw[pattern={north east lines},pattern color=black] (0,0) rectangle (1em,1ex);} marks the cell of level $\ell$ that is marked for coarsening to its parent and the cells \protect\tikz{\protect\draw[pattern={north west lines},pattern color=black] (0,0) rectangle (1em,1ex);} mark cells that are marked for refinement. The ring around the cell marked for coarsening depicts the region that should be checked for the coarsening neighborhood. That is, it defines the region that should not contain cells of level $\ell+1$ (for $\mathcal{N}_c$) or cells of level $\ell$ that are marked for refinement (for $\mathcal{N}'_c$). The cells for which $\mathcal{N}_c=\emptyset$ are marked in (d) and the cells with $\mathcal{N}'_c=\emptyset$ are marked in (e). The final mesh after refinement and coarsening is depicted in (f). The coarsened elements that satisfy $\mathcal{N}_c\cup\mathcal{N}'_c=\emptyset$ are marked as coarsened elements.}
 \label{fig:neighborhood_crs}
\end{figure}

%

\subsection{Labeling Methods}\label{subsec:labeling}
Let $\tilde{\mathbcal{Q}}$ be the ordered set of $\mathbcal{Q}$ such that $e_k\geq e_{k+1}\:\forall k\in\tilde{\mathbcal{Q}}$, where $e_k$ denotes the error of element $k$. Then, the Dörfler marking strategy \cite{Dorfler1996} is defined as the elements $Q_i\in\tilde{\mathbcal{Q}},\:i=0,...,k$ such that the sum of their respective errors is smaller than a fraction $\rho_r$ of the total element error $e=\sum_i e_i$:
\begin{equation}
    \mathbcal{M}^{\text{Dörfler}}_r = \qty{ Q_i\in\mathbcal{Q} : \sum_{k=0}^i e_k < \rho_r e}.
 \end{equation}
This marking strategy, however, does not take into account the contributions of the elements that are marked because they are part of a refinement neighborhood of a marked element $Q_i\in\mathbcal{M}_r$. Therefore, we define the index set $\mathbcal{I}^K_r$ as the set of element indices whose span contains elements $Q_k\in\tilde{\mathbcal{Q}}$ and their refinement neighborhoods:
\begin{equation}\label{eq:indexRef}
\mathbcal{I}^K_r:=\qty{k\in1,...,K:\:\mathbcal{N}_r(Q_k,\mathbcal{Q},m)\cup Q_k,\:Q_k\in\tilde{\mathbcal{Q}}}, 
\end{equation}
and we define $\kappa_r$ as the maximum index for which the sum of all elements with indices $i$ in $\mathbcal{I}^{\kappa_r}_r$ is smaller than the error tolerance $\rho_r e$:
\begin{equation}
\kappa_r:=\arg\max \sum_{i\in\mathbcal{I}^K_r}e_i<\rho_r e,
\end{equation}
such that the Dörfler marking including refinement neighborhoods is
\begin{equation}\label{eq:markRef}
  \mathbcal{M}_r = \qty{ Q_k\in\mathbcal{Q} : k\in\mathbcal{I}^{\kappa_r}_r}.
\end{equation}
For marking a set of coarsening elements, $\mathbcal{M}_c$ the Dörfler marking procedure can be followed again. The original Dörfler marking strategy would be coarsening the elements $Q_i\in\mathbcal{Q}$ such that their total element error is smaller than a fraction of the total element error $\rho_c e$, with coarsening parameter $\rho_c$:
\begin{equation}
    \mathbcal{M}^{\text{Dörfler}}_c = \qty{ Q_i\in\mathbcal{Q} : \sum_{k=i}^N e_k < \rho_c e}.
 \end{equation}
Similar to marking for refinement, the marking rule for coarsening can be specified more precisely by including admissible coarsening. In this case, the elements for which $\hat{\mathbcal{N}}_c(\mathbcal{Q},Q,m,\mathbcal{M}_r)=\emptyset$ holds are added in the sum of marked elements. Therefore, let us define the index set $\mathbcal{J}_K$ that contains all elements $Q_k\in\mathbcal{Q}$ for which the admissible coarsening condition holds, starting from the element with the smallest error, i.e. $Q_N$.
\begin{equation}\label{eq:indexCrs}
\mathbcal{I}^K_c:=\qty{k\in1,...,K:\:\hat{\mathbcal{N}}_c(\mathbcal{Q},Q,m,\mathbcal{M}_r)=\emptyset,\:Q_{N-k-1}\in\tilde{\mathbcal{Q}}}.
\end{equation}
Similar to $\kappa_r$, we define $\kappa_c$ as the maximum index for which the sum of all elements with indices $i$ in $\mathbcal{I}^{\kappa_c}_c$ is smaller than the error tolerance $\rho_c e$:
\begin{equation}
\kappa_c:=\arg\max \sum_{i\in\mathbcal{I}^K_r}e_i<\rho_c e,
\end{equation}
such that the Dörfler marking strategy taking into account coarsening admissibility is defined as
\begin{equation}\label{eq:markCrs}
  \mathbcal{M}_c = \qty{ Q_k\in\mathbcal{Q} : k\in\mathbcal{I}^{\kappa_c}_c}.
\end{equation}

An alternative to the Dörfler marking strategy is a strategy where a given fraction of the total number of elements is marked. In that case, the formulations from \cref{eq:markRef,eq:markCrs} would still hold, but in \cref{eq:indexRef,eq:indexCrs} the indices $\kappa_r$ and $\kappa_c$ are defined by the sum of the marked elements in respectively $\mathbcal{I}^K_r$ and $\mathbcal{I}^K_c$.\\

Whether to mark a set for refinement or coarsening, i.e. to construct $\mathbcal{M}_r$ and $\mathbcal{M}_c$ depends on the global error $\Delta\mathbcal{L}$ following from the DWR and user-defined tolerances for refinement and coarsening. Let $\text{tol}_r$ be the tolerance for refinement and $\text{tol}_c$ the tolerance for coarsening, such that $\mathbcal{M}_r\neq\emptyset$ if and only if $\Delta\mathbcal{L}>\text{tol}_r$ and $\mathbcal{M}_c\neq\emptyset$ if and only if $\Delta\mathbcal{L}<\text{tol}_r$. As a consequence, if $\text{tol}_r\geq\text{tol}_c$, refinement and coarsening are never performed simultaneously. If $\text{tol}_r<\text{tol}_c$ a band with bandwidth $\text{tol}_c-\text{tol}_r$ is defined, in which refinement and coarsening are performed simultaneously. In the present work, tolerances are defined such that the latter condition is satisfied, and the adaptivity iterations are terminated when $\Delta\mathbcal{L}\in[\text{tol}_r,\text{tol}_c]$, i.e.:
\begin{equation}\label{eq:markBand}
 \begin{dcases}
  \mathbcal{M}_r = \emptyset,\:\mathbcal{M}_c \neq \emptyset & \text{if}\:\Delta\mathbcal{L} < \text{tol}_r,\\
  \mathbcal{M}_r \neq \emptyset,\:\mathbcal{M}_c = \emptyset & \text{if}\:\Delta\mathbcal{L} > \text{tol}_c,\\
  \mathbcal{M}_r \neq \emptyset,\:\mathbcal{M}_c \neq \emptyset & \text{if}\:\text{tol}_r \geq \Delta\mathbcal{L} \geq \text{tol}_c,\\
 \end{dcases}
\end{equation}
given $\text{tol}_r \leq \text{tol}_c$.\\

Note that the total element error $e$ and the total estimated error of the system of equations $\Delta\mathbcal{L}$ are not necessarily the same, since the element error measure $e_k$ can be defined in different ways. In case of the DWR method, a natural choice is to choose $e_k$ as the element-wise integrals of $\Delta\mathbcal{L}$ from \cref{eq:errorest}. However, integrating the squared norm of the integrant from \cref{eq:errorest} would yield strictly positive element errors, making the ordering of the set of element errors simple.

\subsection{Quasi-Interpolation}\label{subsec:interpolate}
In the discrete setting, solution of the problem $\vb{u}_h$ is represented by the THB-spline basis $\bm{\phi}_i\in\mathbcal{S}_h^p$ together with the solution coefficients $\alpha_i\in\mathbb{R}$. In case of analyses with multiple solution steps (e.g. dynamic or quasi-static analysis), mesh refinements can be performed after each solution step. As a consequence, the solution at load step $k+1$ is defined on another set of basis functions $\{\bar{\bm{\phi}_i}\}\in\mathbcal{S}_h^p$ with corresponding coefficients $\bar{\alpha}^k_i$ compared to the previous solution at step $k$. In order to transfer the coefficients $\alpha_i^k$ to the new basis, an interpolation scheme needs to be used.\\

Interpolation on a spline basis can be a costly part of the simulation. Global interpolation implies that the contributions of all basis functions are taken into account in the interpolation. This requires solving a large dense system. An efficient way of interpolating spline coefficients for hierarchical basis is a so-called quasi-interpolation scheme \cite{Speleers2016,Giust2020}. Here, on each level of the hierarchical basis a quasi-interpolant is constructed. This quasi-interpolant interpolates a given function $f$ over the support of each basis function individually, to find the coefficient related to that basis function. More precisely, given a function $f\in C(\Omega^0)$, the quasi-interpolant for level $\ell$ is defined as
\begin{equation}
 \Lambda^\ell(f) = \sum_{i=1}^{N} \lambda_{i,\ell}(f)B_{i,\ell}, \quad \ell=0,\dots,n-1,
\end{equation}
where the coefficients $\lambda_{i,\ell}$ are suitable linear functionals on $C(\Omega^0)$. Across all levels $\ell=0,\dots,N-1$, the interpolant for the function becomes:
\begin{equation}\label{eq:QI}
 \Lambda(f) = \sum_{\ell=0}^{N-1}\sum_{i\in \mathbcal{I}_{\ell,\Omega_n}} \lambda_{i,\ell}(f)B_{i,\ell},
\end{equation}
where $B_{i,\ell,\Omega_n}^\mathbcal{T}$ is a THB spline of level $\ell$ constructed on domain $\Omega_n$. For any basis function $B_{i,\ell}$, the coefficient $\lambda_{i,\ell}$ are found by locally interpolating the function $f$ onto all active basis functions $B_{j,\ell}, j\in\mathbcal{J}$ in the support of $B_{i,\ell}$. This gives coefficients $\lambda_{j,\ell}, j\in\mathbcal{J}$ of which coefficient $i$ gives $\lambda_{i,\ell}$. \RevTWO{This quasi-interpolation scheme is used in the present framework to express the solution obtained from the previous load-step in terms of the newly, adaptively refined and coarsened basis. In the case of non-nested spaces -- which can occur when coarsening -- this implies that the quasi-interpolation scheme is not exact.}

\section{Algorithmic Overview}\label{sec:algorithm}
\begin{figure}
\centering
 \pgfdeclarelayer{bg}    
\pgfsetlayers{bg,main}  

\tikzstyle{arrow} = [draw, -latex,thick,line cap=round,line join=round]
\tikzstyle{line} = [draw,thick,line cap=round,line join=round]
\tikzstyle{block} = [rectangle, draw, fill=white,
    text width=0.15\linewidth, text centered, rounded corners, minimum height=11ex]
\tikzstyle{dummyblock} = [rectangle, fill=none,
    text width=0.15\linewidth, text centered, rounded corners, minimum height=11ex]
\tikzstyle{noblock} = [rectangle, draw=none, fill=none,
    text width=0.15\linewidth, rounded corners, minimum height=0em,node distance = 0.05\linewidth, text = blue]
\tikzstyle{textblock} = [text width=7em,text centered]
\tikzstyle{textblock2} = [text width=14em,text centered]
\tikzstyle{textblock4} = [text width=28em,text centered]

\begin{tikzpicture}[node distance = 0.5cm]
    \node[block                      ] (solve)        {\textbf{Solve} \\ Eq.~\eqref{eq:nlproblem_discretised}};

    \node[block,  right = of solve   ] (estimate)     {\textbf{Estimate} \\ Eqs.~\ref{eq:adjoint_discrete}, \ref{eq:adjointenriched}, \ref{eq:errorest}};
    \path [arrow] (solve) -- (estimate) node[above,midway]{$\vb{u}_h$};

    \node[block,  right= of estimate] (mark)         {\textbf{Mark}\\ Eqs.~\ref{eq:markCrs}, \ref{eq:markRef}, \ref{eq:markBand}};
    \path [arrow] (estimate) -- (mark) node[above,midway]{$e_k$};

    \node[block,  right= of mark    ] (refine)       {\textbf{Adapt} \\ Fig.~\ref{fig:IGArefinement}};
    \path [arrow] (mark) -- (refine) node[above,midway]{$\mathbcal{M}_r$};
    \path [arrow] (mark) -- (refine) node[below,midway]{$\mathbcal{M}_c$};

    \node[dummyblock,  above= of solve    ] (empty2)       {};

    \node[block,  right= of refine    ] (transfer)       {\textbf{Transfer} \\ Eq.~\eqref{eq:QI}};
    \node[dummyblock,  above= of transfer    ] (empty)       {};
    \path [arrow] (refine) -- (transfer);

    \node[block,  above= of mark  ] (advance)       {\textbf{Advance}};
    \path [line] (transfer) -- (empty.center);
	\path [arrow] (empty.center) -- (advance);
    \path [arrow] (advance) -| (solve);

   \node[noblock,  below= of solve] (blockbelowsolve)       {};
	\node[noblock,  below= of refine] (blockbelowrefine)       {};
	\node[noblock,  below= of estimate] (blockbelowestimate)       {};
	\node[noblock,  below= of transfer] (blockbelowtransfer)       {};
	\node[noblock,  above= of empty] (blockaboveempty)       {};
	\node[noblock,  above= of empty2] (blockaboveempty2)       {};
    \path [arrow] (transfer)--(blockbelowtransfer.center) -- node[midway,above]{\textbf{while} $\Delta\mathbcal{L}\not\in[\text{tol}_r,\text{tol}_c]$ \textbf{and } $i < I_{\text{max}}$ \textbf{do}} (blockbelowsolve.center)--(solve);

    \path [arrow,dashed] ([xshift=-20]solve.north east) |- ([yshift=4]advance.south west);

\begin{pgfonlayer}{bg}    
	\draw[draw=col1!20,fill=col1!20][rounded corners,name=B]
([xshift=-2,yshift=2]transfer.north west) --
([xshift=-2,yshift=-2]blockbelowtransfer.south west) --
([xshift=2,yshift=-2]blockbelowtransfer.south east) --
([xshift=2,yshift=2]transfer.north east) --
([xshift=2,yshift=2]empty.north east) --
([xshift=2,yshift=2]advance.north west) --
([xshift=2,yshift=-2]advance.south west) --
([xshift=-2,yshift=-2]empty.south west)
--cycle;
\draw[draw=black!20,fill=black!20][rounded corners,name=B]
([xshift=-2,yshift=-2]solve.south west) --
([xshift=2,yshift=-2]solve.south east) --
([xshift=2,yshift=-2]empty2.south east) --
([xshift=2,yshift=-2]advance.south east) --
([xshift=2,yshift=2]advance.north east) -- node[yshift=2,midway](steppinglabel){}
([xshift=-2,yshift=2]empty2.north west)
--cycle;
	\draw[draw=col1,fill=col1,opacity=0.8][rounded corners,name=B] ([xshift=-2,yshift=2]solve.north west) --
([xshift=2,yshift=2]refine.north east) --
([xshift=2,yshift=-2]blockbelowrefine.south east)--node[yshift=-2,midway](meshinglabel){}
([xshift=-2,yshift=-2]blockbelowsolve.south west)
--cycle;
\end{pgfonlayer}

\node[textblock2,anchor=south,color=col1,yshift=4] at (empty.north west) {Addition for solution stepping};
\node[textblock4,anchor=north,color=col1] at (meshinglabel) {Adaptive meshing};
\node[textblock2,anchor=south,color=black] at (steppinglabel) {Classical solution stepping};

\end{tikzpicture}
 \caption{A graphical summary of the adaptive meshing flowchart from \cref{fig:flowchart} used in the present work. The equations which are used in each step are indicated in the blocks. The adaptive meshing iterations are performed within each solution step until the total error $\Delta\mathbcal{L}$ is contained in the interval $[\text{tol}_r,\text{tol}_c]$, $\text{tol}_r<\text{tol}_c$, following the tolerances in \cref{subsec:labeling}. In case of convergence, the solution is advanced, e.g. with an arc-length iteration. \Cref{alg:flowchart2} provides an algorithm corresponding to this flow-chart.}
 \label{fig:flowchart2}
 \end{figure}
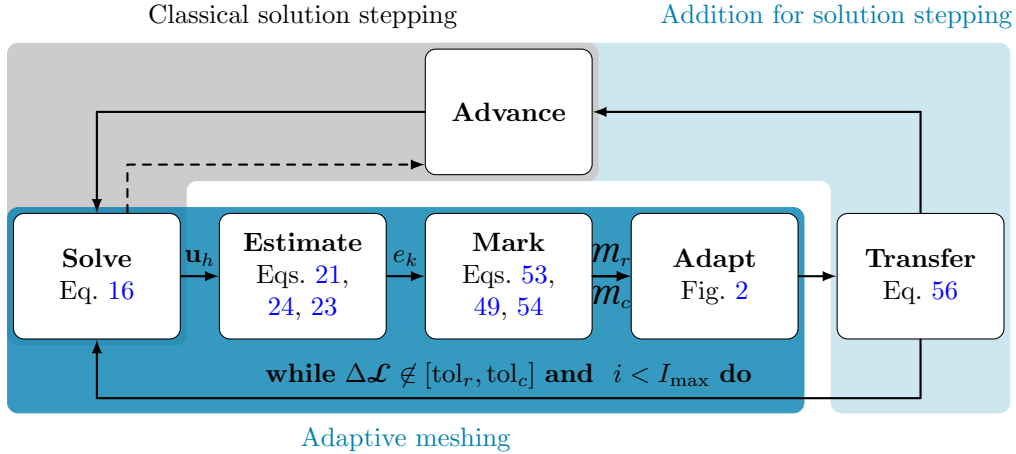

In \cref{fig:flowchart} the adaptive isogeometric method for solution stepping problems has been presented. Based on  \cref{sec:shell,sec:DWR,sec:refinement}, a summarised workflow for adaptive isogeometric shell analysis is depicted in \cref{fig:flowchart2,alg:flowchart2}.\\

The \emph{Solve} block involves solving the non-linear isogeometric Kirchhoff-Love shell equation from \cref{eq:nlproblem_discretised}. This variational formulation involves geometric and material non-linearities and can potentially also involve load non-linearities. After solving the Kirchhoff-Love shell problem, the discrete solution vector $\vb{u}_h$ is passed to the \emph{Estimate} block. Here, the DWR method is solved by computing the adjoint problem in the primal space (\cref{eq:adjoint_discrete}) and in the enriched space (\cref{eq:adjointenriched}). Then, the element-wise error estimate can be obtained by integrating \cref{eq:errorest} element-wise. The element-wise errors $e_k$ can be passed to the \emph{Mark} block, where elements are marked for refinement (\cref{eq:markCrs}) if the total error $\Delta\mathbcal{L}$ is larger than a lower (refinement) tolerance $\text{tol}_r$ and a coarsening marking (\cref{eq:markCrs}) is performed if the total error is above an upper (coarsening) tolerance $\text{tol}_c$. This implies that if $\text{tol}_c<\Delta\mathbcal{L}<\text{tol}_r$, a combined coarsening and refinement step is performed, as described in \cref{eq:markBand}. In this case, the coarsening marking from \cref{eq:markCrs} is performed given $\mathbcal{M}_r$. Given the elements marked for refinement and coarsening, collected in $\mathbcal{M}_r$ and $\mathbcal{M}_c$ respectively, the mesh can be \emph{Adapt}ed. In order to start the solution interval again, the start point should be \emph{Transfer}red to the new mesh and the governing equation can be solved again if the error is not in the interval $[\text{tol}_r,\text{tol}_c]$ or if the number of refinement iterations $i$ exceeds the maximum number of refinement iterations, $I_{\text{max}}$. If the total error is in the interval $[\text{tol}_r,\text{tol}_c]$ or if the number of refinement iterations $i$ exceeds the maximum number of refinement iterations, $I_{\text{max}}$, the solution can be advanced, e.g. using a load-stepping or an arc-length method. Thereafter, the governing equations can be \emph{Solve}d again. Note that if $I_{\text{max}}=1$, no inner iterations for adaptive meshing are performed.

\begin{algorithm}
    \hspace*{\SpaceReservedForComments}{}%
    \begin{minipage}{\dimexpr\linewidth-\SpaceReservedForComments\relax}
    \captionsetup{font=footnotesize} 
    \caption{An algorithmic summary of the goal-adaptive meshing routine employed in the present work. \Cref{fig:flowchart2} provides a graphical summary of this algorithm.}
    \label{alg:flowchart2}
    \begin{algorithmic}[1]
    \FOR {loadsteps}
    \WHILE{$\Delta\mathbcal{L}\not\in[\text{tol}_l,\text{tol}_u]\text{\textbf{ and }} i < I_{\text{max}}$}
        \STATE Compute the primal solution $\vb{u}_h$ (Eq.~\eqref{eq:nlproblem_discretised})
        \STATE Compute the dual solution, $\bm{\xi}_h$ given $\vb{u}_h$ (Eq.~\eqref{eq:adjoint_discrete})
        \STATE Compute the enriched dual solution, $\tilde{\bm{\xi}}_h$ given $\vb{u}_h$ (Eq.~\eqref{eq:adjointenriched})
        \STATE Compute the total error estimation $\Delta\mathbcal{L}$ according to Eq.~\eqref{eq:errorest} and find the element-wise errors $e_k$ and the total element error $e$.
        \IF{$\Delta\mathbcal{L}>\text{tol}_r$}
        \STATE Mark elements for refinement into $\mathbcal{M}_r$ using $e$, see Eq.~\eqref{eq:markRef}.
        \ENDIF
        \IF{$\Delta\mathbcal{L}<\text{tol}_c$}
        \STATE Mark elements for coarsening into $\mathbcal{M}_c$ using $e$, see Eq.~\eqref{eq:markCrs}.
        \ENDIF
        \STATE Refine all $Q\in\mathbcal{M}_r$ using THB-splines, see Fig.~\ref{fig:IGArefinement}
        \STATE Coarsen all $Q\in\mathbcal{M}_c$  using THB-splines, see Fig.~\ref{fig:IGArefinement}
        \STATE Transfer the solutions required to start the new solution step to the new mesh using Quasi-Interpolation, see Eq.~\eqref{eq:QI}.
    \ENDWHILE
    \STATE Advance the solution to the next solution step
    \ENDFOR
    \end{algorithmic}
    \AddNote[black]{3}{3}{\textbf{Solve}}
    \AddNote[black]{4}{6}{\textbf{Estimate}}
    \AddNote[black]{7}{12}{\textbf{Mark}}
    \AddNote[black]{13}{14}{\textbf{Adapt}}
    \AddNote[black]{15}{15}{\textbf{Transfer}}
    \AddNote[black]{17}{17}{\textbf{Advance}}
    \end{minipage}
\end{algorithm}

\section{Numerical examples}\label{sec:results}
In this section, several numerical examples are presented. The examples represent different applications of the theory presented in this paper and - without loss of generality - all employ Isogeometric Kirchhoff-Love shells. The first three examples illustrate the performance of the DWR error estimator and the last three example demonstrate the use of this error estimator for adaptive meshing. More precisely, the numerical examples performed in this section, as well as their purpose, are:
\begin{description}
\item[Linear static analysis of a square plate (\cref{subsec:linearShell})] A simple example of linear Kirchhoff-Love shell theory is presented. In this case, error estimators using the DWR are computed for different goal functionals and verified using the actual error computed from manufactured solutions. The goal of this benchmark problem is to evaluate the accuracy of the error estimators in linear static analysis.
\item[Modal analysis of a circular plate (\cref{subsec:modal})] Since the analytical eigenvalues and eigenmodes are known for this case, the goal of this benchmark problem is to verify the error estimator for a vibration eigenvalue problem, given in \cref{eq:eigenvalueproblem_discrete} with the stiffness and mass operators  $\mathbcal{A}(\vb{v},\bm{\phi})=\mathbcal{W}'(\vb{0},\vb{v},\bm{\phi})$ and $\mathbcal{B}(\vb{v},\bm{\phi})=\mathbcal{M}(\vb{v},\bm{\phi})$.
\item[Linear buckling analysis of a square plate (\cref{subsec:buckling})] Analytical critical buckling loads and mode shapes are also known for this case. Therefore, the goal of this benchmark problem is to provide verification for the buckling error estimator from \cref{eq:eigenvalueproblem_discrete} with the buckling operators $\mathbcal{A}(\vb{v},\bm{\phi})=\mathbcal{W}'(\vb{u}_L,\vb{v},\bm{\phi})$ and $\mathbcal{B}(\vb{v},\bm{\phi})=\mathbcal{M}(\vb{v},\bm{\phi})$.
\item[Non-linear analysis of a pinched thin plate (\cref{subsec:non-linearShell})] In this example a thin plate with very low bending stiffness subject to a out-of-plane load is analysed. The error estimator is used to provide mesh adaptivity to compare to uniform refinement. The goal of this benchmark problem is to evaluate the performance of the DWR as driver for adaptive meshing in a static load case with geometric non-linearities.
\item[Snap-through instability of a cylindrical roof (\cref{subsec:roof})] The snap-through behaviour of a cylindrical roof are considered in this example. The benchmark problem is a well-known application of arc-length methods and shells. The goal of solving this problem is to test the full adaptive solution stepping procedure from \cref{fig:flowchart2} on a benchmark problem.
\item[Wrinkling analysis (\cref{subsec:wrinkling})] In the last example, the procedure from \cref{fig:flowchart2} is applied to the modeling of membrane wrinkling. This problem contains geometric non-linearities and material non-linearities. The results are compared to uniformly refinements to evaluate the efficiency of adaptive meshing for such applications. The goal of this example is to demonstrate the use of the adaptive meshing routine from \cref{fig:flowchart2} on a complex load-stepping problem with geometric and material non-linearities.
\end{description}


In the following subsections, the short-hand notations $\mathbcal{L}_{\text{an}}=\mathbcal{L}(\vb*{u}_{\text{an}})$, $\mathbcal{L}_{\text{num}}=\mathbcal{L}(\vb*{u}_{\text{num}})$, $\Delta\mathbcal{L}_{\text{an}}=\mathbcal{L}_{\text{an}} - \mathbcal{L}_{\text{num}}$ and $\Delta\mathbcal{L}_{\text{num}}=\mathbcal{R}(\vb*{u}_{\text{num}},\tilde{\bm{\xi}}_h-\bm{\xi}_h)$ (see \cref{eq:errorest}) are used, given the analytical and numerical solutions $\vb*{u}_{\text{an}}$ and $\vb*{u}_{\text{num}}$, respectively. Furthermore, where relevant, the parameters $\rho_r$, $\rho_c$, $\text{tol}_r$ and $\text{tol}_c$ (see \cref{eq:markBand,eq:markCrs,eq:markRef}) are fixed per example. A study on finding optimal values for these parameters is out of the scope of this paper. Lastly, all simulations are performed using the open-source Geometry+Simulation modules \cite{Juttler2014}.

\subsection{Linear static analysis of a square plate}\label{subsec:linearShell}
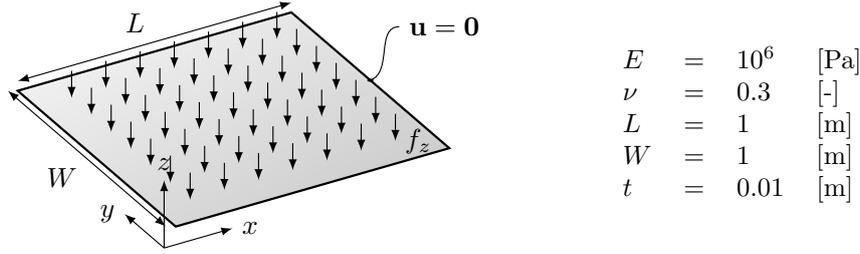
\begin{figure}
 \centering
 \begin{minipage}{0.55\linewidth}
 \resizebox{\linewidth}{!}
 {
 \def\N{20}
\tdplotsetmaincoords{60}{-30}
\begin{tikzpicture}[tdplot_main_coords,scale=2]
\node (A) at (-1,-1,0){};
\node (B) at (1,-1,0){};
\node (C) at (1,1,0){};
\node (D) at (-1,1,0){};
\node (Z) at (0,0,0){};
\filldraw[thick,black,top color=black!05,bottom color=black!20,shading angle=30](A.center)--node [midway](AB){} (B.center)-- node [midway](BC){} (C.center)-- node [midway](CD){}(D.center)--cycle node [midway](DA){};

\coordinate (O) at (-1.2,-1.2,0);
\draw[-latex](O)--++(0.5,0,0) node[right]{$x$};
\draw[-latex](O)--++(0,0.5,0) node[left]{$y$};
\draw[-latex](O)--++(0,0,0.5) node[above]{$z$};

\draw[latex-latex](D.west)--node[midway,below left]{$W$}(A.west);
\draw[latex-latex](C.north)--node[midway,above left]{$L$}(D.north);

\foreach \x in {-0.75,-0.5,...,0.75}
{
	\foreach \y in {-0.75,-0.5,...,0.75}
	{
		\draw[latex-] (\x,\y,0) --++(0,0,0.2);
	}
}
\node[right] at (0.75,-0.75,0) {$f_z$};

\node (label) at (1,0,0) {};
\fill (label)  circle [radius=0.02];
\draw (label.center) to[in=180,out=20] ++(0.5,0.5,0) node[right]{$\mathbf{u}=\mathbf{0}$};
\end{tikzpicture}
 }
 \end{minipage}
 \hfill
 \begin{minipage}{0.35\linewidth}
 \begin{tabular}{llll}
 $E$    &$=$& $10^6$ & $[\text{Pa}]$ \\
 $\nu$ &$=$& $0.3$ & $[\text{-}]$ \\
 $L$    &$=$& $1$ & $[\text{m}]$\\
 $W$    &$=$& $1$ & $[\text{m}]$\\
 $t$     &$=$& $0.01$ & $[\text{m}]$\\
\end{tabular}
 \end{minipage}
 \caption{Geometry and parameters for the example of a unit-square plate with a distributed vertical load $f_z$ given by \cref{eq:manufactured_load}. The displacements are fixed on all edges.}
 \label{fig:setupLinear}
\end{figure}

For the linear shell example, let us consider a flat plate with unit dimensions $L=W=1\:[\text{m}]$, a thickness of $t=10^{-2}\:[\text{m}]$ and with material parameters $E=10^6\:[\text{Pa}]$, $\nu=0.3$, which is clamped on all sides, see \cref{fig:setupLinear}. A load vector of
\begin{multline}\label{eq:manufactured_load}
    \vb{f} = \frac{2EAt^3}{1-\nu^2}\bigg(x^4 - 2x^3 + 3x^2 - 2x + y^4 - 2y^3 + 3y^2 - 2y + 12x^2y^2 \\- 12x^2y -12xy^2 + 12xy + \frac{1}{3}\bigg) \vb{e}_z
\end{multline}
is applied, based on the manufactured solution given by
\begin{equation}\label{eq:manufactured_solution}
\vb*{u}_{\text{an}} = A x^2 (x - 1)^2 y^2 (y - 1)^2 \vb{e}_z.
\end{equation}
Using this manufactured solution, any goal functional $\mathbcal{L}_{\text{an}}$ can be evaluated. Solving the primal problem for this linear shell example gives $\vb*{u}_{\text{num}}$, which can be used to compute the DWR error estimate of $\Delta\mathbcal{L}_\text{num}$.\\

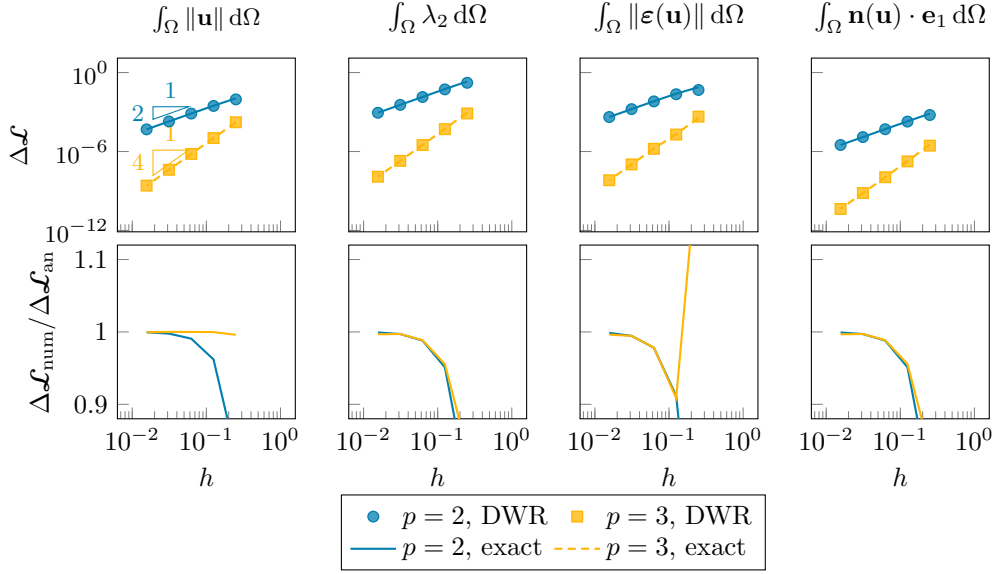
\begin{figure}
\centering
\begin{tikzpicture}
  \begin{groupplot}[
      group style={
          group name=my plots,
          group size=4 by 2,
          xlabels at=edge bottom,
          xticklabels at=edge bottom,
          ylabels at=edge left,
          yticklabels at=edge left,
          vertical sep=5pt,
          horizontal sep=20pt,
          group name=group
      },
      width=0.3\linewidth,
      height=0.2\textheight,
      xlabel={$h$},
      xmin=1e-2, xmax=1e0,
      ymin=1e-11, ymax=1e0,
      enlarge x limits = true,
      enlarge y limits = true,
      tickpos=left,
      xmode = log,
      ymode = log,
      enlargelimits,
      restrict expr to domain={x}{1e-5:2.5e-1},
  ]
\nextgroupplot[title = {$\int_\Omega\Vert\vb{u}\Vert\dd{\Omega}$},ylabel={$\Delta\mathbcal{L}$}]
\addplot+[style=p2*,only marks,fill opacity=1.0] table[header=true,x expr = 2^(-\thisrowno{0}),y expr = abs(\thisrowno{1}), col sep = comma, skip first n=3]{Data/Linear/example_shell3D_DWR_r6_e1_g1_c9.csv};
\addplot+[style=p2*,no markers,fill opacity=1.0] table[header=true,x expr = 2^(-\thisrowno{0}),y expr = abs(\thisrowno{2}), col sep = comma, skip first n=3]{Data/Linear/example_shell3D_DWR_r6_e1_g1_c9.csv};
\addplot+[style=p3*,only marks,fill opacity=1.0] table[header=true,x expr = 2^(-\thisrowno{0}),y expr = abs(\thisrowno{1}), col sep = comma, skip first n=2]{Data/Linear/example_shell3D_DWR_r6_e2_g1_c9.csv};
\addplot+[style=p3*,no markers,fill opacity=1.0] table[header=true,x expr = 2^(-\thisrowno{0}),y expr = abs(\thisrowno{2}), col sep = comma, skip first n=2]{Data/Linear/example_shell3D_DWR_r6_e2_g1_c9.csv};
\logLogSlopeTriangleRev{0.4}{0.2}{0.72}{2}{col1}
\logLogSlopeTriangleRev{0.4}{0.2}{0.47}{4}{col2}

\nextgroupplot[title = {$\int_\Omega\lambda_2\dd{\Omega}$}]
\addplot+[style=p2*,only marks,fill opacity=1.0] table[header=true,x expr = 2^(-\thisrowno{0}),y expr = abs(\thisrowno{1}), col sep = comma, skip first n=3]{Data/Linear/example_shell3D_DWR_r6_e1_g2_c1.csv};
\addplot+[style=p2*,no markers,fill opacity=1.0] table[header=true,x expr = 2^(-\thisrowno{0}),y expr = abs(\thisrowno{2}), col sep = comma, skip first n=3]{Data/Linear/example_shell3D_DWR_r6_e1_g2_c1.csv};
\addplot+[style=p3*,only marks,fill opacity=1.0] table[header=true,x expr = 2^(-\thisrowno{0}),y expr = abs(\thisrowno{1}), col sep = comma, skip first n=2]{Data/Linear/example_shell3D_DWR_r6_e2_g2_c1.csv};
\addplot+[style=p3*,no markers,fill opacity=1.0] table[header=true,x expr = 2^(-\thisrowno{0}),y expr = abs(\thisrowno{2}), col sep = comma, skip first n=2]{Data/Linear/example_shell3D_DWR_r6_e2_g2_c1.csv};

\nextgroupplot[title = {$\int_\Omega\Vert\bm{\varepsilon}(\vb{u})\Vert\dd{\Omega}$}]
\addplot+[style=p2*,only marks,fill opacity=1.0] table[header=true,x expr = 2^(-\thisrowno{0}),y expr = abs(\thisrowno{1}), col sep = comma, skip first n=3]{Data/Linear/example_shell3D_DWR_r6_e1_g3_c9.csv};
\addplot+[style=p2*,no markers,fill opacity=1.0] table[header=true,x expr = 2^(-\thisrowno{0}),y expr = abs(\thisrowno{2}), col sep = comma, skip first n=3]{Data/Linear/example_shell3D_DWR_r6_e1_g3_c9.csv};
\addplot+[style=p3*,only marks,fill opacity=1.0] table[header=true,x expr = 2^(-\thisrowno{0}),y expr = abs(\thisrowno{1}), col sep = comma, skip first n=2]{Data/Linear/example_shell3D_DWR_r6_e2_g3_c9.csv};
\addplot+[style=p3*,no markers,fill opacity=1.0] table[header=true,x expr = 2^(-\thisrowno{0}),y expr = abs(\thisrowno{2}), col sep = comma, skip first n=2]{Data/Linear/example_shell3D_DWR_r6_e2_g3_c9.csv};

\nextgroupplot[title = {$\int_\Omega\vb{n}(\vb{u})\cdot \vb{e}_1\dd{\Omega}$}]
\addplot+[style=p2*,only marks,fill opacity=1.0] table[header=true,x expr = 2^(-\thisrowno{0}),y expr = abs(\thisrowno{1}), col sep = comma, skip first n=3]{Data/Linear/example_shell3D_DWR_r6_e1_g7_c0.csv};
\addplot+[style=p2*,no markers,fill opacity=1.0] table[header=true,x expr = 2^(-\thisrowno{0}),y expr = abs(\thisrowno{2}), col sep = comma, skip first n=3]{Data/Linear/example_shell3D_DWR_r6_e1_g7_c0.csv};
\addplot+[style=p3*,only marks,fill opacity=1.0] table[header=true,x expr = 2^(-\thisrowno{0}),y expr = abs(\thisrowno{1}), col sep = comma, skip first n=2]{Data/Linear/example_shell3D_DWR_r6_e2_g7_c0.csv};
\addplot+[style=p3*,no markers,fill opacity=1.0] table[header=true,x expr = 2^(-\thisrowno{0}),y expr = abs(\thisrowno{2}), col sep = comma, skip first n=2]{Data/Linear/example_shell3D_DWR_r6_e2_g7_c0.csv};

\nextgroupplot[ymode=normal,ymin = 0.9, ymax = 1.1,ylabel={$\Delta\mathbcal{L}_{\text{num}}/\Delta\mathbcal{L}_{\text{an}}$},legend to name={CommonLegend},legend columns=2, transpose legend]
\addlegendimage{style=p2*,only marks}\addlegendentry{$p=2$, DWR}
\addlegendimage{style=p2*,mark=none}\addlegendentry{$p=2$, exact}
\addlegendimage{style=p3*,only marks}\addlegendentry{$p=3$, DWR}
\addlegendimage{style=p3*,mark=none}\addlegendentry{$p=3$, exact}
\addplot+[style=p2*,no markers,solid,fill opacity=1.0] table[header=true,x expr = 2^(-\thisrowno{0}),y index = {3}, col sep = comma, skip first n=3]{Data/Linear/example_shell3D_DWR_r6_e1_g1_c9.csv};
\addplot+[style=p3*,no markers,solid,fill opacity=1.0] table[header=true,x expr = 2^(-\thisrowno{0}),y index = {3}, col sep = comma, skip first n=2]{Data/Linear/example_shell3D_DWR_r6_e2_g1_c9.csv};

\nextgroupplot[ymode=normal,ymin = 0.9, ymax = 1.1]
\addplot+[style=p2*,no markers,solid,fill opacity=1.0] table[header=true,x expr = 2^(-\thisrowno{0}),y index = {3}, col sep = comma, skip first n=3]{Data/Linear/example_shell3D_DWR_r6_e1_g2_c1.csv};
\addplot+[style=p3*,no markers,solid,fill opacity=1.0] table[header=true,x expr = 2^(-\thisrowno{0}),y index = {3}, col sep = comma, skip first n=2]{Data/Linear/example_shell3D_DWR_r6_e2_g2_c1.csv};

\nextgroupplot[ymode=normal,ymin = 0.9, ymax = 1.1]
\addplot+[style=p2*,no markers,solid,fill opacity=1.0] table[header=true,x expr = 2^(-\thisrowno{0}),y index = {3}, col sep = comma, skip first n=3]{Data/Linear/example_shell3D_DWR_r6_e1_g3_c9.csv};
\addplot+[style=p3*,no markers,solid,fill opacity=1.0] table[header=true,x expr = 2^(-\thisrowno{0}),y index = {3}, col sep = comma, skip first n=2]{Data/Linear/example_shell3D_DWR_r6_e2_g3_c9.csv};

\nextgroupplot[ymode=normal,ymin = 0.9, ymax = 1.1]
\addplot+[style=p2*,no markers,solid,fill opacity=1.0] table[header=true,x expr = 2^(-\thisrowno{0}),y index = {3}, col sep = comma, skip first n=3]{Data/Linear/example_shell3D_DWR_r6_e1_g7_c0.csv};
\addplot+[style=p3*,no markers,solid,fill opacity=1.0] table[header=true,x expr = 2^(-\thisrowno{0}),y index = {3}, col sep = comma, skip first n=2]{Data/Linear/example_shell3D_DWR_r6_e2_g7_c0.csv};

    \end{groupplot}
\path (group c1r2.south west) -- node[below=25pt]{\ref*{CommonLegend}} (group c4r2.south east);

\end{tikzpicture}
\caption{Linear static analysis of a clamped plate with a uniformly distributed load according to \cref{eq:manufactured_load} according to a manufactured solution from \cref{eq:manufactured_solution}. The top row provides $\Delta\mathbcal{L}$ against the uniform mesh size $h$. The bottom row presents the efficiency of the error estimators against the mesh size $h$. The markers represent error estimates computed via the DWR and the lines represent the exact error, i.e. the error of the numerical solution with respect to the analytical solution. The results are given for goal functionals (from left to right) $\mathbcal{L}=\int_\Omega\Vert\vb{u}\Vert\dd{\Omega}$ (displacement norm), $\mathbcal{L}=\int_\Omega\lambda_2\dd{\Omega}$ (second principal stretch), $\mathbcal{L}=\int_\Omega\Vert\bm{\varepsilon}(\vb{u})\Vert\dd{\Omega}$ (membrane strain norm) and $\mathbcal{L}=\int_\Omega\vb{n}(\vb{u})\cdot \vb{e}_1\dd{\Omega}$ (first component of membrane force).}
\label{fig:linear_static}
\end{figure}

In \cref{fig:linear_static}, the results for the linear shell problem are given. The title of each column represents the goal-function that is used for error estimation in this column. The top row provides the errors $\Delta\mathbcal{L}_{\text{an}}$ and $\Delta\mathbcal{L}_{\text{num}}$ with respect to an uniformly refined mesh size.  As can be seen in this figure, $\Delta\mathbcal{L}_{\text{num}}$ quickly converges to $\Delta\mathbcal{L}_{\text{an}}$ for different spline orders $p$. In addition, the bottom row of \cref{fig:linear_static} provides the efficiency of the error estimator, given by $\Delta\mathbcal{L}_{\text{num}}/\Delta\mathbcal{L}_{\text{an}}$. These figures confirm convergence of the DWR estimates to the analytical goal functional errors for all considered goal functionals. Only for the membrane strain norm goal functional the error estimate for coarse meshes is inaccurate. This can possibly be explained by the in-plane shear strain that cancels out over the whole domain but which does contribute in the norm $\Vert\bm{\varepsilon}\Vert$\\

Concluding, the linear shell benchmarks shows that for different goal functionals the DWR method provides accurate estimation of the error $\Delta\mathbcal{L}$ of the goal functional $\mathbcal{L}$ starting at relatively small mesh sizes of $h<10^{-1}$.

\subsection{Modal analysis of a circular plate}\label{subsec:modal}
As a next example, the vibration modes of a circular plate with clamped boundary conditions are computed. The geometry with boundary conditions is illustrated in \cref{fig:setupModal}. The circular plate has unit diameter, a thickness $t=10^{-2}\:[\text{m}]$, Young's modulus $E=10^6\:[\text{Pa}]$, Poisson's ratio $\nu=0.3$ and density $\rho=1\:[\text{kg}/\text{m}^3]$. The analytical solutions for the eigenfrequencies of the circular plate are obtained by
\begin{equation}
    \omega_n = \gamma_n^2\sqrt{D/\rho t},
\end{equation}
where $\gamma_n$ is the $n^\text{th}$ root of the equation $(I_{m-1}(\gamma R)-m/RI_{m}(\gamma R))J_m(\gamma R) - I_m(\gamma R)(J_{m-1}(\gamma R)-m/RJ_m(\gamma R))=0$ following from a separation of variables solution \cite{Verhelst2019}, $R$ is the radius of the plate and $D=Et^3/(12(1-\nu^2))=9.16\cdot 10^{-8}$ is the flexural rigidity. As stated in \cref{subsec:eigenvalueproblem}, the goal functional eigenvalue problems is given in \cref{eq:eiggoal}, and requires the eigenvalue problem in \cref{eq:eigenvalueproblem} to be solved with linear operators $\mathbcal{A}(\vb{v},\bm{\phi})=\mathbcal{W}'(\vb{0},\vb{v},\bm{\phi})$ and $\mathbcal{B}(\vb{v},\bm{\phi})=\mathbcal{M}(\vb{v},\bm{\phi})$, see \cref{eq:mass,eq:linearisationN}.\\

\begin{figure}
 \centering
 \begin{minipage}{0.55\linewidth}
 \resizebox{\linewidth}{!}
 {
 \def\N{20}
\tdplotsetmaincoords{60}{-30}
\begin{tikzpicture}[tdplot_main_coords,scale=2]
\node (Z) at (0,0,0) {Z};
\filldraw[thick,black,top color=black!10,bottom color=black!50,shading angle=30] (Z) circle [radius=2];
\coordinate (O) at (Z);
\draw[-latex](O)--++(0.5,0,0) node[right]{$x$};
\draw[-latex](O)--++(0,0.5,0) node[left]{$y$};
\draw[-latex](O)--++(0,0,0.5) node[above]{$z$};

\draw[-latex] (Z.center) -- node[midway,right]{$R$} (1.41,1.41,0);
\node (label) at (-1.41,1.41,0) {};
\fill (label)  circle [radius=0.02];
\draw (label.center) to[in=150,out=20] ++(0.15,-1,0) node[right]{$\mathbf{u}=\nabla \mathbf{u}\cdot\mathbf{n}=\mathbf{0}$};
\end{tikzpicture}
 }
 \end{minipage}
 \hfill
 \begin{minipage}{0.35\linewidth}
 \begin{tabular}{llll}
 $E$    &$=$& $10^6$ & $[\text{Pa}]$ \\
 $\nu$ &$=$& $0.3$ & $[\text{-}]$ \\
 $\rho$ &$=$& $1$ & $[\text{kg}/\text{m}^3]$ \\
 $R$    &$=$& $0.5$ & $[\text{m}]$\\
 $t$     &$=$& $0.01$ & $[\text{m}]$\\
\end{tabular}
 \end{minipage}
 \caption{Geometry and parameters for a vibrating circular plate with a clamped boundary.}
 \label{fig:setupModal}
\end{figure}

\begin{figure}
 \centering
 \input{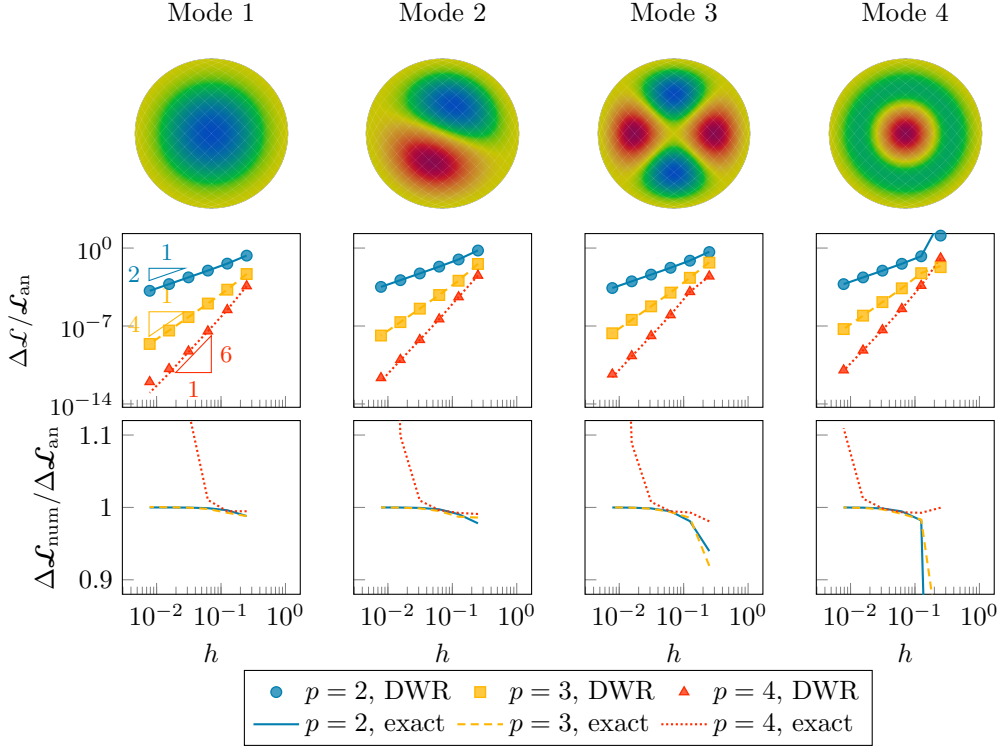}
 \caption{Modal analysis of a circular plate. The top row provides the mode shapes from mode 1 up to 4. The mid row provides $\Delta\mathbcal{L}$ against the uniform mesh size $h$ and the bottom row provides the efficiency of the error estimators against the mesh size $h$. The markers represent error estimates computed via the DWR and the lines represent the exact error, i.e. the error of the numerical solution with respect to the analytical solution. The eigenfrequencies for mode 1 up to 4 are, respectively, $\omega_1=9.56\cdot 10^{-4}\:[\text{Hz}]$, $\omega_2=4.14\cdot 10^{-3}\:[\text{Hz}]$, $\omega_3=1.11\cdot 10^{-2}\:[\text{Hz}]$ and $\omega_4=1.45\cdot 10^{-2}\:[\text{Hz}]$.}
 \label{fig:modal}
 \end{figure}

\Cref{fig:modal} presents the first four eigenmodes in the top row. Furthermore, $\Delta\mathbcal{L}_{\text{an}}$ and $\Delta\mathbcal{L}_{\text{num}}$ as a function of the element size for uniformly refined domain are given in the middle row for the first four eigenmodes. These plots show that the approximation for the error $\Delta\mathbcal{L}_{\text{num}}$ approximates $\Delta\mathbcal{L}_{\text{an}}$. In the bottom row of \cref{fig:modal} the efficiencies also show that the approximation converges to an efficiency equal to 1. However, for the $p=4$ line, the efficiency degrades when the `exact' error obtained by the analytical solution approaches values around $10^{-11}$, which is attributed to the approximation of the roots $\gamma_i$ and the precision of the eigenvalue solver.

Concluding, the modal analysis benchmark shows that the DWR method provides accurate estimation of the eigenfrequency error for different considered mode shapes.

\subsection{Linear buckling analysis of a square plate}\label{subsec:buckling}
\begin{figure}
 \centering
 \begin{minipage}{0.55\linewidth}
 \resizebox{\linewidth}{!}
 {
 \def\N{20}
\tdplotsetmaincoords{60}{-30}
\begin{tikzpicture}[tdplot_main_coords,scale=2]
\node (A) at (-1,-1,0){};
\node (B) at (1,-1,0){};
\node (C) at (1,1,0){};
\node (D) at (-1,1,0){};
\filldraw[black,top color=black!10,bottom color=black!50,shading angle=30](A.center)--node [midway](AB){} (B.center)-- node [midway](BC){} (C.center)-- node [midway](CD){}(D.center)--cycle node [midway](DA){};

\foreach \y in {0,...,\N}
{
	\draw[latex-] (B)++(0,2*\y/\N,0) --++(1/4,0,0) node (B\y){};
	\draw[latex-] (D)++(2*\y/\N,0,0) --++(0,1/4,0) node (D\y){};
}
\draw (B0.center) --(B\N.center);
\draw (D0.center)--(D\N.center) node[above right] {$\sigma t$};

\coordinate (O) at (-1.2,-1.2,0);
\draw[-latex](O)--++(0.5,0,0) node[right]{$x$};
\draw[-latex](O)--++(0,0.5,0) node[left]{$y$};
\draw[-latex](O)--++(0,0,0.5) node[above]{$z$};


\draw[latex-latex](D.west)--node[midway,below left]{$W$}(A.west);
\draw[latex-latex](C.south east)--node[midway,below right]{$L$}(D.south east);

\fill (AB.center) circle [radius=0.02];
\draw (AB.center) to[in=80,out=-30]++(0.2,-0.5,0) node[below] {$u_y=u_z=0$};

\fill (BC.center) circle [radius=0.02];
\draw (BC.center) to[in=180,out=0]++(0.5,0.5,0) node[right] {$u_z=0$};

\fill (CD.center) circle [radius=0.02];
\draw (CD.center) to[in=300,out=120]++(0.6,0.5,0) node[above] {$u_z=0$};

\fill (DA.center) circle [radius=0.02];
\draw (DA.center) to[in=0,out=-180]++(-0.5,-0.5,0) node[left] {$u_x=u_z=0$};
\end{tikzpicture}
 }
 \end{minipage}
 \hfill
 \begin{minipage}{0.35\linewidth}
 \begin{tabular}{llll}
 $E$    &$=$& $10^6$ & $[\text{Pa}]$ \\
 $\nu$ &$=$& $0.3$ & $[\text{-}]$ \\
 $L$    &$=$& $1$ & $[\text{m}]$\\
 $W$   &$=$& $1$ & $[\text{m}]$\\
 $t$     &$=$& $0.01$ & $[\text{m}]$\\
\end{tabular}
 \end{minipage}
 \caption{Geometry and parameters for the plate buckling example. A distributed load of $\sigma t $ is acting on two boundaries in two different directions, and the other boundaries are simply supported and fixed in out-of-plane direction.}
 \label{fig:setupBuckling}
\end{figure}
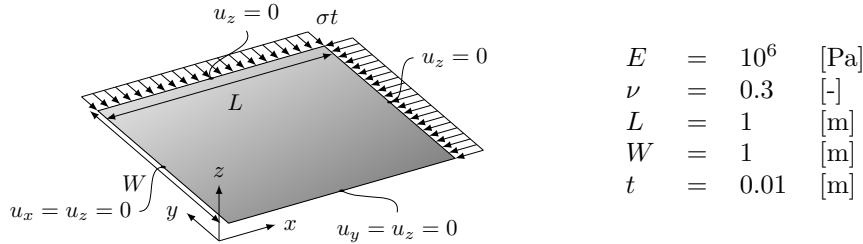

Similar to modal analysis, DWR error estimation for buckling analysis also relies on the formulations in \cref{subsec:eigenvalueproblem}. The difference with the modal analysis error estimation is that the buckling analysis error estimation involves the solution of a pre-buckling solution and no mass matrix. As an example for buckling analysis, a square simply supported plate is considered, see \cref{fig:setupBuckling}, with a Saint-Venant Kirchhoff constitutive law with Young's modulus $E=10^6\:[\text{Pa}]$ and Poisson's ratio $\nu=0.3\:[\text{-}]$. The dimensions of the plate are $L\times W \times t=1\times1\times 0.01\:[\text{m}^3]$. The plate is subject to a distributed line load of $\sigma t$ in both directions. The analytical solution for the buckling load with $m$ half waves in $x$-direction and $n$ half waves in $y$-direction for a square plate with sides $L$ and with equal loads is given in \cite{Jones2006} and reads:
\begin{equation}
 \sigma^{m,n}_{c}t = \frac{D \pi^2}{L^2}\qty(m^2 + n^2),
\end{equation}
with $D=Et^3/12(1-\nu^2)$ the flexural rigidity of the plate. Using this expression, the first four unique modes are, indexed in ascending order: $\sigma^{1,1}_{c}t = \sigma^1_{c}t = 1.808\:[\text{N}/\text{m}]$, $\sigma^{2,1}_{c}t = \sigma^{1,2}_{c}t = \sigma^2_{c}t = 4.519\:[\text{N}/\text{m}]$, $\sigma^{2,2}_{c}t = \sigma^3_{c}t = 7.230\:[\text{N}/\text{m}]$, $\sigma^{3,1}_{c}t = \sigma^{1,3}_{c}t = \sigma^4_{c}t = 9.038\:[\text{N}/\text{m}]$.

\begin{figure}
	\centering
	\input{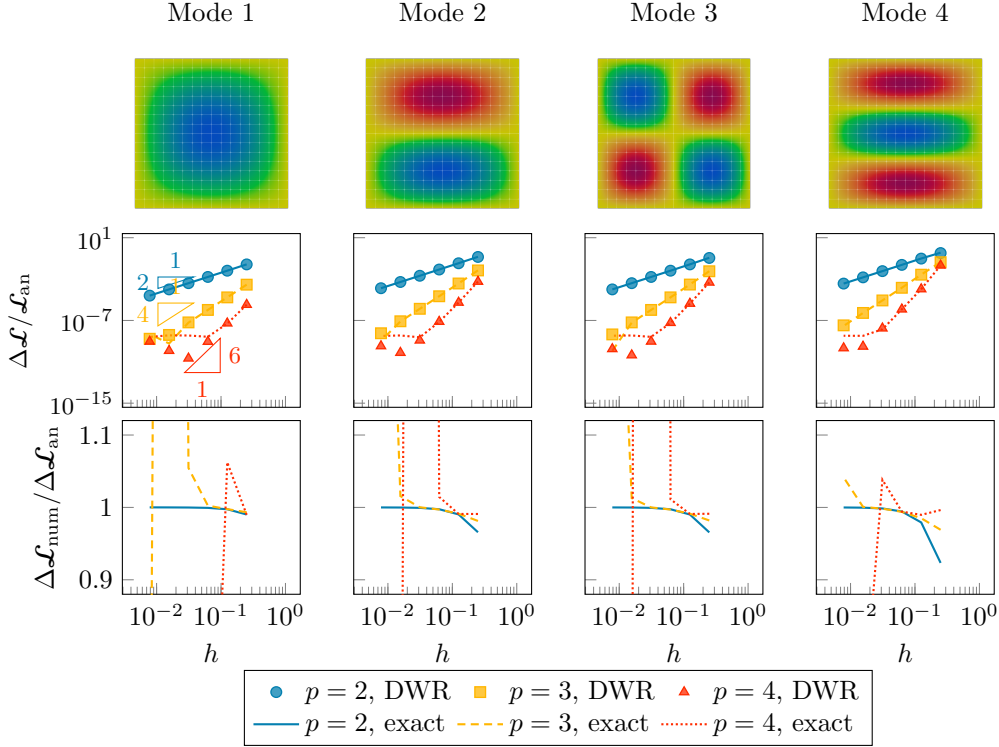}
	 \caption{Buckling analysis of a square plate with simply supported boundary conditions. The top row provides the mode shapes from mode 1 up to 4. The bottom row provides $\Delta\mathbcal{L}$ against the uniform mesh sise $h$. The markers represent error estimates computed via the DWR and the lines represent the exact error, i.e. the error of the numerical solution with respect to the analytical solution. The critical loads for unique modes 1 up to 4 are, respectively, $\sigma^{1,1}_{c}t = \sigma^1_{c}t = 1.808\:[\text{N}/\text{m}]$, $\sigma^{2,1}_{c}t = \sigma^{1,2}_{c}t = \sigma^2_{c}t = 4.519\:[\text{N}/\text{m}]$, $\sigma^{2,2}_{c}t = \sigma^3_{c}t = 7.230\:[\text{N}/\text{m}]$, $\sigma^{3,1}_{c}t = \sigma^{1,3}_{c}t = \sigma^4_{c}t = 9.038\:[\text{N}/\text{m}]$.}
	 \label{fig:buckling}
 \end{figure}

\Cref{fig:buckling} depicts the analytical error $\Delta\mathbcal{L}_{\text{an}}$ and the DWR error estimate $\Delta\mathbcal{L}_{\text{num}}$ as a function of the mesh size $h$ for uniform refinements. Both errors are normalised with the analytical value of the critical buckling load. As can be seen in this figure, the DWR prediction of the error converges with a rate of convergence of $2(p-1)$ for all degrees $p$ until it reaches values of around $10^{-10}$ after which the errors stagnate and increase again (in particular for $p=4$). This behaviour is similar to the behaviour observed in \cite{Liu2021c} and can be attributed to the round-off errors as discussed. These errors occur when the number of degrees of freedom is large enough and the machine precision is limited. In case of a buckling simulation, where the non-linear stiffness operator is constructed on an initial solution of a linear simulation, the influence of round-off errors is expected to occur sooner. In addition, it can also be seen that the error computed using the analytical solution stagnates. This is due to the fact that the numerical approximation of the critical buckling load shows small variations depending on the solution to the linear problem that is solved to obtain the tangential stiffness matrix to compute the generalised eigenvalue problem for buckling. For the results presented in \cref{fig:buckling}, the load $\sigma t=10^{-4}\:[\text{N}/\text{m}]$ was used to compute the buckling linearisation.\\

\subsection{Non-linear analysis of a pinched thin plate}\label{subsec:non-linearShell}
\begin{figure}
 \centering
 \begin{minipage}{0.55\linewidth}
 \resizebox{\linewidth}{!}
 {
 \def\N{20}
\tdplotsetmaincoords{60}{-30}
\begin{tikzpicture}[tdplot_main_coords,scale=2]
\node (A) at (-1,-1,0){};
\node (B) at (1,-1,0){};
\node (C) at (1,1,0){};
\node (D) at (-1,1,0){};
\node (Z) at (0,0,0){};
\filldraw[black,top color=black!05,bottom color=black!20,shading angle=30](A.center)--node [midway](AB){} (B.center)-- node [midway](BC){} (C.center)-- node [midway](CD){}(D.center)--cycle node [midway](DA){};

\coordinate (O) at (-1.2,-1.2,0);
\draw[-latex](O)--++(0.5,0,0) node[right]{$x$};
\draw[-latex](O)--++(0,0.5,0) node[left]{$y$};
\draw[-latex](O)--++(0,0,0.5) node[above]{$z$};

\filldraw[black,top color=black!10,bottom color=black!50,shading angle=30] (B.center) -- node[midway](CCD){} (AB.center) -- node[midway](ABZ){} (Z.center) -- node[midway](BCZ){} (BC.center);


\draw[latex-latex](D.west)--node[midway,below left]{$W$}(A.west);
\draw[latex-latex](C.north)--node[midway,above left]{$L$}(D.north);

\fill (A.center) circle [radius=0.02];
\draw (A.center) to[in=80,out=-30]++(0.2,-0.5,0) node[below] {$\mathbf{u}=\mathbf{0}$};

\fill (B.center) circle [radius=0.02];
\node[below right] at (B.center) {$\mathbf{u}=\mathbf{0}$};

\fill (C.center) circle [radius=0.02];
\node[above right] at (C.center) {$\mathbf{u}=\mathbf{0}$};

\fill (D.center) circle [radius=0.02];
\node[above left] at (D.center) {$\mathbf{u}=\mathbf{0}$};

\fill (Z.center) circle [radius=0.01];
\draw[-latex,thick] (Z.center)--++(0,0,0.5) node[right]{$P$};

\node[above,rotate=-40] at (ABZ.center) {symm};

\node[above,rotate=16] at (BCZ.center) {symm};
\end{tikzpicture}
 }
 \end{minipage}
 \hfill
 \begin{minipage}{0.35\linewidth}
 \begin{tabular}{llll}
 $E$    &$=$& $1.0$ & $[\text{MPa}]$ \\
 $\nu$ &$=$& $0.3$ & $[\text{-}]$ \\
 $P$    &$=$& $4\cdot10^{-7}$ & $[\text{N}]$ \\
 $L$    &$=$& $1$ & $[\text{mm}]$\\
 $W$   &$=$& $1$ & $[\text{mm}]$\\
 $t$     &$=$& $10^{-3}$ & $[\text{mm}]$\\
\end{tabular}
 \end{minipage}
 \caption{Geometry and parameters for a square thin plate subject to a point load $P$ in the middle. The plate is fully constrained in every corner. Because the problem is symmetric, only a quarter of the domain is modelled. Hence, symmetry conditions are applied. On the $x$-aligned symmetry axis, this implies that $u_y=\pdv{u_z}{y}=\pdv{u_x}{y}=0$ and on the $y$-aligned symmetry axis this implies that $u_x=\pdv{u_z}{x}=\pdv{u_y}{x}=0$.}
 \label{fig:setupNonlinear}
\end{figure}
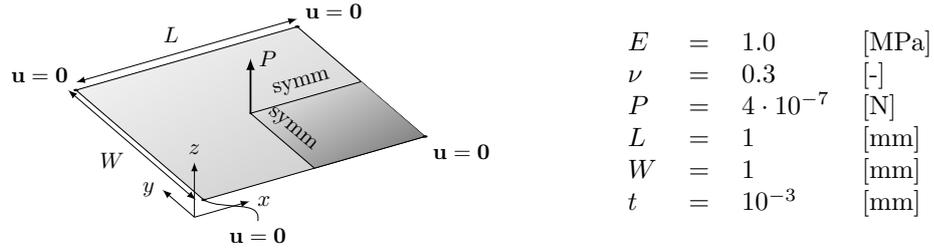

In a next example, we consider a square membrane subject to a point-load in the middle and with corners fixed in all directions, see \cref{fig:setupNonlinear}. The membrane is modelled with a Saint-Venant Kirchhoff constitutive law with Young's modulus $E=1.0\:[MPa]$ and a Poisson ratio $\nu=0.3$. The thickness of the membrane is $t=10^{-3}\:[\text{mm}]$ and the length and with are $L\times W = 1\times1\:[\text{mm}]$. The simulation is performed on a quarter of the domain, employing symmetry conditions as depicted in \cref{fig:setupNonlinear}. A load of $P=4\cdot10^{-7}\:[\text{N}]$ is applied in the center of the sheet. The static load case is solved using an arc-length method to ensure convergence of the solution. Furthermore, an adaptive refinement strategy is employed with admissible refinement. The jump parameter $m$ is set to 2 and the maximum number of refinement levels is $8$ or $11$, which equals a tensor basis level with $2^{8}\times2^{8}=256\times256$ or $2^{11}\times2^{11}=2048\times2048$ elements, respectively. The refinement parameter is set to $\rho_r=0.5$. The goal functionals considered in this case are based on displacements as well as on principal stresses:
\begin{equation}\label{eq:goalFunctionalsNonlinear}
\begin{aligned}
\mathbcal{L}(\vb{u}) &=\Vert\vb{u}(\vb{x}_P)\Vert,\\
\mathbcal{L}(\vb{u}) &=\int_\Omega \bm{\sigma}_p\cdot\vb{e}_y\dd{\Omega}.
\end{aligned}
\end{equation}
\RevTWO{Instead of being a domain-integrated goal functional, the first goal functional is evaluated on the point $\vb{x}_P$ where the force $P$ is applied, see \cref{fig:setupNonlinear}.}

\begin{figure}
\centering
\includegraphics[width=0.45\linewidth]{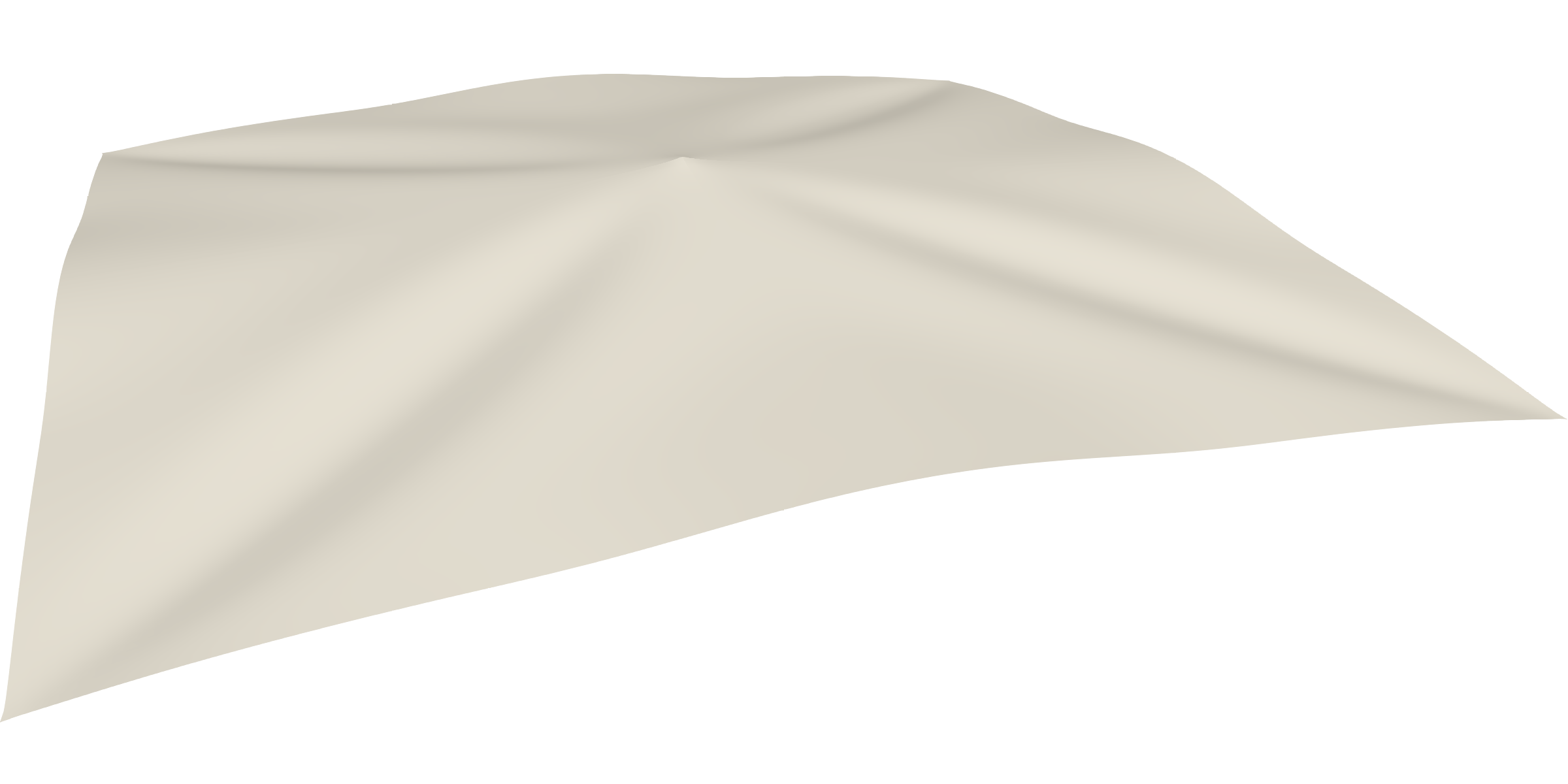}
\caption{Deformed surface from the benchmark presented in \cref{fig:setupNonlinear}. The result is the last solution from the adaptive meshing routine with deformation norm goal-functional of which the results are presented in \cref{fig:PointLoad_errors}.}
\label{fig:PointLoad_converged}
\end{figure}

\Cref{fig:PointLoad_converged} presents the deformed membrane for the last step of the adaptive simulation. Furthermore, \cref{fig:PointLoad_errors} presents the estimated error $\Delta\mathbcal{L}$ given the goal functionals in \cref{eq:goalFunctionalsNonlinear} for the uniform refinement as well as for the adaptive refinement simulation with maximum level $8$ or $11$. Moreover, \cref{fig:elementErrorsNonlinear_g1C9,fig:elementErrorsNonlinear_g6C1} provide the absolute element-wise errors for the adaptive refinement simulation and for the uniform refinement series for both considered goal functionals. The contour lines in these error fields represent the vertical deflection of the sheet.\\

From all provided results, it can be observed that the adaptive mesh provides for both goal functionals an efficient converging mesh, where the accuracy per degree of freedom is higher compared to uniformly refined meshes. However, it can also be seen that the total error $\Delta\mathbcal{L}$ is not strictly decreasing for both goal functionals. The bottom plots of \cref{fig:PointLoad_errors} indicate that this point is closely related to the maximum depth that is reached: the percentage of the total element error $e$ that can still be refined rapidly decreases, meaning that the only elements that are still available for refinement are the ones that have insignificant contribution to the total element error $e$, deeming refinement of these elements meaningless. It can be seen from \cref{fig:elementErrorsNonlinear_g1C9,fig:elementErrorsNonlinear_g6C1} that the maximum element depth is reached in the corner where the sheet is fixed \RevTWO{and in the corner where the load is applied}. After the maximum level is reached, the refined elements start distributing over the diagonal of the domain. For both goal functionals, \cref{fig:PointLoad_errors} shows that some further decrease in the total error $\Delta\mathbcal{L}$ can be gained after the maximum error is reached, but that it is most effective to increase the maximum refinement depth. Comparing the error fields and meshes for both goal functionals (see \cref{fig:elementErrorsNonlinear_g1C9,fig:elementErrorsNonlinear_g6C1}), it can be seen that the second-principal stress-based error field shows slightly wider error bands in the finest depicted uniform refinement error fields than the displacement-based error estimator. As a consequence, the corresponding adaptive meshes show that the elements indeed tend to be broader distributed along the diagonal of the domain in case of the displacement-driven refinement (\cref{fig:elementErrorsNonlinear_g1C9}).\\

\begin{figure}
\centering
        \begin{tikzpicture}
    \begin{groupplot}[
        group style={
            group name=group,
            group size=2 by 2,
            xlabels at=edge bottom,
            xticklabels at=edge bottom,
            ylabels at=edge left,
            vertical sep=5pt,
            horizontal sep=0.1\linewidth
        },
        width=0.45\linewidth,
        height=0.2\textheight,
        xlabel={\#DoFs},
        tickpos=left,
        xmode = log,
        grid = major,
        log origin=infty,
        xmin = 1e1,
        xmax = 1e5,
        enlarge x limits=true,
    ]


    \nextgroupplot[ylabel = {$\Delta\mathbcal{L}$},ymode = log, title = {$\mathbcal{L} = \Vert\vb{u}(\vb{x}_P)\Vert$},legend to name={CommonLegend}, legend columns=3,ymax=1e-3,ymin=1e-6]
    \addplot+[style=p2*,skip coords between index={6}{999}] table[header=true,x expr = \thisrowno{7},y index = {1}, col sep = comma]{Data/PointLoadNonlinear/NonlinearAdaptive_g1C9_point_max8.csv};\addlegendentry{Adaptive (8)}
    \addplot+[style=p3*,skip coords between index={9}{999}] table[header=true,x expr = \thisrowno{7},y index = {1}, col sep = comma]{Data/PointLoadNonlinear/NonlinearAdaptive_g1C9_point.csv};\addlegendentry{Adaptive (11)}
    \addplot+[style=p4*] table[header=true,x expr = \thisrowno{7},y index = {1}, col sep = comma]{Data/PointLoadNonlinear/NonlinearUniform_g1C9_point.csv};\addlegendentry{Uniform}

    \addplot+[style=p2,skip coords between index={0}{6}] table[header=true,x expr = \thisrowno{7},y index = {1}, col sep = comma]{Data/PointLoadNonlinear/NonlinearAdaptive_g1C9_point_max8.csv};\addlegendentry{..., blocked}
    \addplot+[style=p3,skip coords between index={0}{8}] table[header=true,x expr = \thisrowno{7},y index = {1}, col sep = comma]{Data/PointLoadNonlinear/NonlinearAdaptive_g1C9_point.csv};\addlegendentry{..., blocked}

    \addplot[forget plot,only marks,black,mark indices={8,10,18,21},mark=o] table[header=true,x expr = \thisrowno{7},y index = {1}, col sep = comma]{Data/PointLoadNonlinear/NonlinearAdaptive_g1C9_point_max8.csv};
    \addplot[forget plot,only marks,black,mark indices={8,10,18,21},mark=square] table[header=true,x expr = \thisrowno{7},y index = {1}, col sep = comma]{Data/PointLoadNonlinear/NonlinearAdaptive_g1C9_point.csv};
    \addplot[forget plot,only marks,black,mark indices={1,2,3,4},mark=triangle] table[header=true,x expr = \thisrowno{7},y index = {1}, col sep = comma]{Data/PointLoadNonlinear/NonlinearUniform_g1C9_point.csv};

    \nextgroupplot[ymode = log,title = {$\mathbcal{L} = \int_\Omega \bm{\sigma}_p\cdot\vb{e}_y\dd{\Omega}$}]
    \addplot+[style=p2*,skip coords between index={6}{999}] table[header=true,x expr = \thisrowno{7},y index = {1}, col sep = comma]{Data/PointLoadNonlinear/NonlinearAdaptive_g6C1_max8.csv};
    \addplot+[style=p3*,skip coords between index={9}{999}] table[header=true,x expr = \thisrowno{7},y index = {1}, col sep = comma]{Data/PointLoadNonlinear/NonlinearAdaptive_g6C1.csv};
    \addplot+[style=p4*] table[header=true,x expr = \thisrowno{7},y index = {1}, col sep = comma]{Data/PointLoadNonlinear/NonlinearUniform_g6C1.csv};

    \addplot+[style=p2,skip coords between index={0}{6}] table[header=true,x expr = \thisrowno{7},y index = {1}, col sep = comma]{Data/PointLoadNonlinear/NonlinearAdaptive_g6C1_max8.csv};
    \addplot+[style=p3,skip coords between index={0}{8}] table[header=true,x expr = \thisrowno{7},y index = {1}, col sep = comma]{Data/PointLoadNonlinear/NonlinearAdaptive_g6C1.csv};

    \addplot[forget plot,only marks,black,mark indices={8,10,18,21},mark=o] table[header=true,x expr = \thisrowno{7},y index = {1}, col sep = comma]{Data/PointLoadNonlinear/NonlinearAdaptive_g6C1_max8.csv};
    \addplot[forget plot,only marks,black,mark indices={8,10,18,21},mark=square] table[header=true,x expr = \thisrowno{7},y index = {1}, col sep = comma]{Data/PointLoadNonlinear/NonlinearAdaptive_g6C1.csv};
    \addplot[forget plot,only marks,black,mark indices={1,2,3,4},mark=triangle] table[header=true,x expr = \thisrowno{7},y index = {1}, col sep = comma]{Data/PointLoadNonlinear/NonlinearUniform_g6C1.csv};

    \nextgroupplot[ylabel = {Ref. \% of $e$},ymode = log]
    \addplot+[style=p2*,skip coords between index={6}{999}] table[header=true,x expr = \thisrowno{7},y expr = \thisrowno{11}/(\thisrowno{10}+\thisrowno{11}), col sep = comma]{Data/PointLoadNonlinear/NonlinearAdaptive_g1C9_point_max8.csv};
    \addplot+[style=p2,skip coords between index={0}{5}] table[header=true,x expr = \thisrowno{7},y expr = \thisrowno{11}/(\thisrowno{10}+\thisrowno{11}), col sep = comma]{Data/PointLoadNonlinear/NonlinearAdaptive_g1C9_point_max8.csv};

    \addplot+[style=p3*,skip coords between index={9}{999}] table[header=true,x expr = \thisrowno{7},y expr = \thisrowno{11}/(\thisrowno{10}+\thisrowno{11}), col sep = comma]{Data/PointLoadNonlinear/NonlinearAdaptive_g1C9_point.csv};
    \addplot+[style=p3,skip coords between index={0}{8}] table[header=true,x expr = \thisrowno{7},y expr = \thisrowno{11}/(\thisrowno{10}+\thisrowno{11}), col sep = comma]{Data/PointLoadNonlinear/NonlinearAdaptive_g1C9_point.csv};

    \nextgroupplot[ymode = log]
    \addplot+[style=p2*,skip coords between index={6}{999}] table[header=true,x expr = \thisrowno{7},y expr = \thisrowno{11}/(\thisrowno{10}+\thisrowno{11}), col sep = comma]{Data/PointLoadNonlinear/NonlinearAdaptive_g6C1_max8.csv};
    \addplot+[style=p2,skip coords between index={0}{5}] table[header=true,x expr = \thisrowno{7},y expr = \thisrowno{11}/(\thisrowno{10}+\thisrowno{11}), col sep = comma]{Data/PointLoadNonlinear/NonlinearAdaptive_g6C1_max8.csv};

    \addplot+[style=p3*,skip coords between index={9}{999}] table[header=true,x expr = \thisrowno{7},y expr = \thisrowno{11}/(\thisrowno{10}+\thisrowno{11}), col sep = comma]{Data/PointLoadNonlinear/NonlinearAdaptive_g6C1.csv};
    \addplot+[style=p3,skip coords between index={0}{8}] table[header=true,x expr = \thisrowno{7},y expr = \thisrowno{11}/(\thisrowno{10}+\thisrowno{11}), col sep = comma]{Data/PointLoadNonlinear/NonlinearAdaptive_g6C1.csv};
    \end{groupplot}
    \path (group c1r2.south east) -- node[midway,below=25pt]{\ref*{CommonLegend}} (group c2r2.south west);
    \end{tikzpicture}
    \caption{Estimated error convergence (top) and the percentage of the total element error $e$ that is available for refinement (bottom) against the number of degrees of freedom (DoFs) for adaptively and uniformly refined meshes with respect to the goal functionals from \cref{eq:goalFunctionalsNonlinear}. The markers labeled with a black border are the markers for which the mesh is plotted in \cref{fig:elementErrorsNonlinear_g1C9,fig:elementErrorsNonlinear_g6C1}. The filled markers represent points where refinement is not blocked, and empty markers points where refinement is blocked because the maximum refinement depth is reached. Note that the errors are computed before refinement, hence blocked elements in iteration $i$ have effect on the error computation in iteration $i+1$.}
\label{fig:PointLoad_errors}
\end{figure}
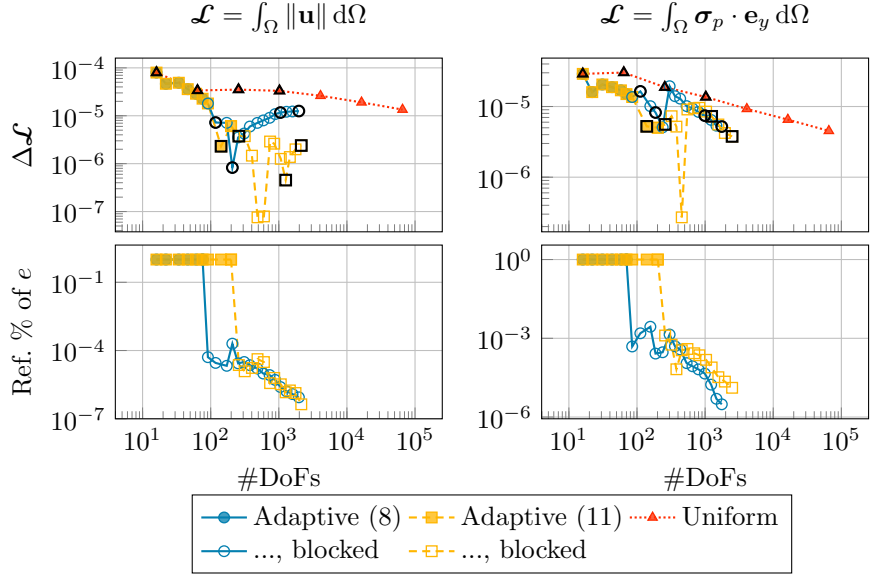

\begin{figure}
\centering
\begin{subfigure}{\linewidth}
\centering
\input{Figures/Results_NonlinearStaticPointLoad_MeshesUniform_g1C9}
\caption{Uniform refinement}
\end{subfigure}

\begin{subfigure}{\linewidth}
\centering
\input{Figures/Results_NonlinearStaticPointLoad_Meshes_g1C9_max8}
\caption{Adaptive refinement with maximum level 8.}
\end{subfigure}

\begin{subfigure}{\linewidth}
\centering
\input{Figures/Results_NonlinearStaticPointLoad_Meshes_g1C9}
\caption{Adaptive refinement with maximum level 11.}
\end{subfigure}
\caption{Normalised element error values $e_k/\Delta\mathcal{L}$ for uniformly (a) and adaptively (b, c) refined meshes using goal-function $\mathbcal{L}(\vb{u}) = \Vert\vb{u}(\vb{x}_P)\Vert$. The meshing steps increase from left to right. The contour lines represent the displacement of the membrane, with intervals of $0.1\:[\text{mm}]$. The bottom right corner of the pictures indicates the fixed corner and the top left corner is the corner where the load is applied.}
\label{fig:elementErrorsNonlinear_g1C9}
\end{figure}

\begin{figure}
\centering
\begin{subfigure}{\linewidth}
\centering
\input{Figures/Results_NonlinearStaticPointLoad_MeshesUniform_g6C1}
\caption{Uniform refinement.}
\end{subfigure}

\begin{subfigure}{\linewidth}
\centering
\input{Figures/Results_NonlinearStaticPointLoad_Meshes_g6C1_max8}
\caption{Adaptive refinement with maximum level 8.}
\end{subfigure}

\begin{subfigure}{\linewidth}
\centering
\input{Figures/Results_NonlinearStaticPointLoad_Meshes_g6C1}
\caption{Adaptive refinement with maximum level 11.}
\end{subfigure}
\caption{Normalised element error values $e_k/\Delta\mathcal{L}^2$ for uniformly (a) and adaptively (b, c) refined meshes using goal-function $\mathbcal{L}(\vb{u}) = \int_\Omega \bm{\sigma}_p\cdot\vb{e}_y\dd{\Omega}$. The meshing steps increase from left to right. The contour lines represent the displacement of the membrane, with intervals of $0.1\:[\text{mm}]$. The bottom right corner of the pictures indicates the fixed corner and the top left corner is the corner where the load is applied.}
\label{fig:elementErrorsNonlinear_g6C1}
\end{figure}

Concluding, the non-linear shell benchmark shows that an adaptive meshing strategy provides accurate solutions on meshes with a small number of degrees of freedom compared to uniform meshes. Furthermore, the benchmarks show the importance of refinement levels, meaning that convergence of the adaptive meshing strategy vanishes as soon as the elements contribution to a large extent to the total error are on the lowest allowed level. This stresses the relevance of spline constructions that allow for deep levels of refinement with moderate computational costs.

\subsection{Snap-through instability of a cylindrical roof}\label{subsec:roof}
In order to present results of the adaptive isogeometric method developed in this paper for quasi-static problems, hence completing full cycles in \cref{fig:flowchart}, the well-known benchmark of a collapsing cylindrical roof \cite{Sze2004} subject to a point-load is considered. The goal of this benchmark problem is to evaluate whether the presented adaptive isogeometric method provides DoF-wise efficient solutions compared to solutions with uniform refinements.\\

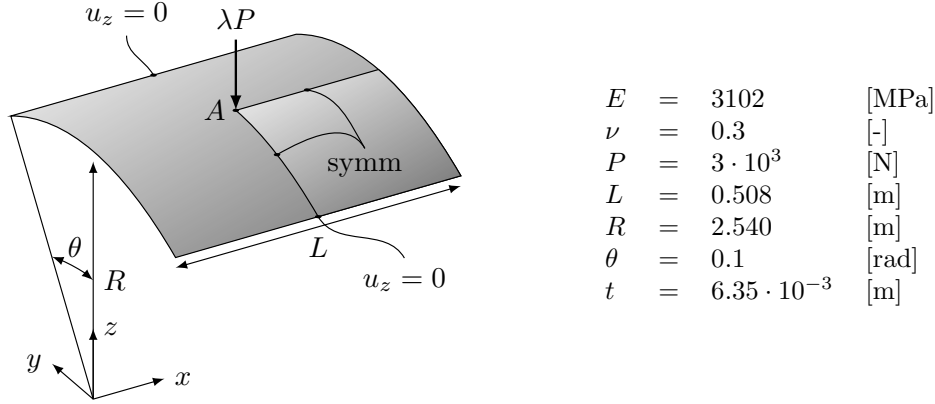
\begin{figure}
 \centering
 \begin{minipage}{0.55\linewidth}
 \resizebox{\linewidth}{!}
 {
 \def\N{20}
\def\mZ{5}
\tdplotsetmaincoords{60}{-30}
\begin{tikzpicture}[tdplot_main_coords,scale=2]
\filldraw[black,top color=black!10,bottom color=black!50,shading angle=30] (1,0,0) node (A){} .. controls (1,0.333333333333333,\mZ*0.033445462325053) and (1,0.666666666666667,\mZ*0.05016819348758) .. (1,1,\mZ*0.05016819348758) -- (1,1,\mZ*0.05016819348758) .. controls (1,1.333333333333333,\mZ*0.05016819348758) and (1,1.666666666666667,\mZ*0.033445462325053) .. (1,2,0) node (B){} -- node[midway](labeltop){} (-1,2,0) node (C){} .. controls (-1,1.666666666666667,\mZ*0.033445462325053) and (-1,1.333333333333333,\mZ*0.05016819348758) .. (-1,1,\mZ*0.05016819348758) -- node(mid){} (-1,1,\mZ*0.05016819348758) .. controls (-1,0.666666666666667,\mZ*0.05016819348758) and (-1,0.333333333333333,\mZ*0.033445462325053) .. (-1,0,0) node (D){}--  node[midway](labelbot){} cycle;


\draw[latex-latex] (A.south) -- (D.south) node[midway, below]{$L$};
\coordinate (O) at (-1,1,-1.5);
\draw[-latex](O)--++(0.5,0,0) node[right]{$x$};
\draw[-latex](O)--++(0,0.5,0) node[left]{$y$};
\draw[-latex](O)--++(0,0,0.5) node[right]{$z$};
\draw[latex-] (mid)-- node[midway](mR){} (O);
\draw(O)-- node[midway](mphi){} (-1,2,0);

\node[right] at (mR){$R$};

\draw[latex-latex] (mR.center) arc [start angle=0,end angle=22,radius=1] node[midway,above]{$\theta$};

\filldraw[black,top color=black!10,bottom color=black!50,shading angle=30] (1,0,0) .. controls (1,0.333333333333333,\mZ*0.033445462325053) and (1,0.666666666666667,\mZ*0.05016819348758) .. (1,1,\mZ*0.05016819348758) --  node[midway] (Aver){} (0,1,\mZ*0.05016819348758) .. controls (0,0.666666666666667,\mZ*0.05016819348758) and (0,0.333333333333333,\mZ*0.033445462325053) .. (0,0,0) node[midway] (Ahor){}-- cycle;

\fill (labeltop) circle [radius=0.02];
\draw (labeltop.center) to[in=-90,out=120]++(0.1,0.5,0) node[above]{$u_z=0$};
\fill (labelbot) circle [radius=0.02];
\draw (labelbot.center) to[in=120,out=-60]++(0.2,-0.7,0) node[below]{$u_z=0$};

\fill (Ahor.center) circle [radius=0.02];
\draw (Ahor.center) to[in=120,out=30]++(0.5,-0.2,0) node(label){};
\fill (Aver.center) circle [radius=0.02];
\draw (Aver.center) to[in=120,out=-10] (label.center){};

\node[below] at (label) {symm};

\fill (0,1,\mZ*0.05016819348758) circle [radius=0.02];
\node[left] at (0,1,\mZ*0.05016819348758) {$A$};
\draw[latex-,thick] (0,1,\mZ*0.05016819348758)--++(0,0,0.5) node[above]{$\lambda P$};
\end{tikzpicture}
 }
 \end{minipage}
 \hfill
 \begin{minipage}{0.35\linewidth}
 \begin{tabular}{llll}
 $E$    &$=$& $3102$ & $[\text{MPa}]$ \\
 $\nu$ &$=$& $0.3$ & $[\text{-}]$ \\
 $P$    &$=$& $3\cdot10^3$ & $[\text{N}]$ \\
 $L$    &$=$& $0.508$ & $[\text{m}]$\\
 $R$    &$=$& $2.540$ & $[\text{m}]$\\
 $\theta$&$=$& $0.1$ & $[\text{rad}]$\\
 $t$     &$=$& $6.35\cdot10^{-3}$ & $[\text{m}]$\\
\end{tabular}
 \end{minipage}
 \caption{Geometry and parameters for cylindrical roof with a point-load $P$ in the middle. The bottom-right corner of each domain corresponds to the point $A$: the reference point for which the $z$-displacements are plotted. The roof is free on the curved edges and simply supported ($u_z=0$) on the straight edges. As a consequence, the problem is symmetric and a quarter of the domain is modelled. On the $x$-aligned symmetry axis, this implies that $u_y=\pdv{u_z}{y}=\pdv{u_x}{y}=0$ and on the $y$-aligned symmetry axis this implies that $u_x=\pdv{u_z}{x}=\pdv{u_y}{x}=0$.}
 \label{fig:setupRoof}
\end{figure}

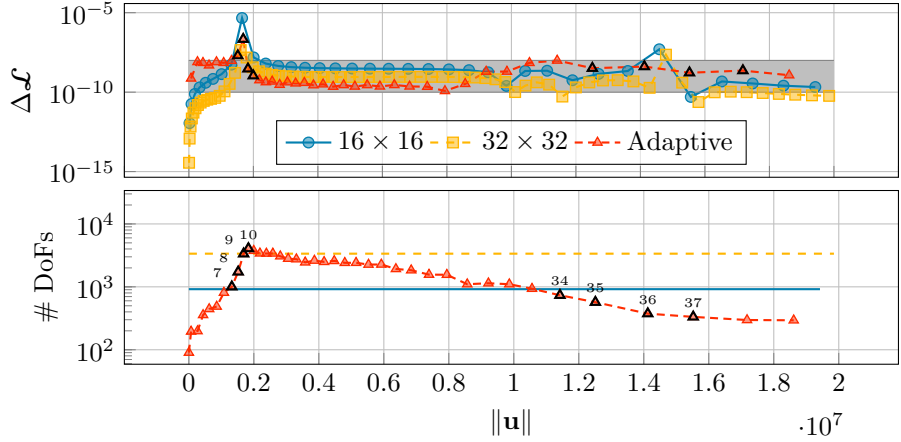
\begin{figure}
\centering
\begin{tikzpicture}
    \begin{groupplot}[
        group style={
            group name=my plots,
            group size=1 by 2,
            xlabels at=edge bottom,
            xticklabels at=edge bottom,
            ylabels at=edge left,
            yticklabels at=edge left,
            vertical sep=5pt,
            horizontal sep=5pt
        },
        width=0.9\linewidth,
        height=0.2\textheight,
        tickpos=left,
        grid = major,
        xlabel={$\Vert\vb{u}\Vert$},
        ymode=log,
        restrict x to domain=0:2e7,
        enlargelimits,
    ]
    \nextgroupplot[ylabel={$\Delta\mathbcal{L}$}, legend columns=3,legend style={at={(rel axis cs: 0.5,0.2)},anchor=center}]
    \addplot+[fill opacity=1.0] table[header=true,x expr = \thisrowno{2},y index = {4}, col sep = comma]{Data/ArcLengthRoof/RoofUniform_r4_g9C9.csv};\addlegendentry{$16\times16$}
    \addplot+[fill opacity=1.0] table[header=true,x expr = \thisrowno{2},y index = {4}, col sep = comma]{Data/ArcLengthRoof/RoofUniform_r5_g9C9.csv};\addlegendentry{$32\times32$}
    \addplot+[fill opacity=1.0] table[header=true,x expr = \thisrowno{2},y index = {4}, col sep = comma]{Data/ArcLengthRoof/RoofAdaptive_g9C9.csv};\addlegendentry{Adaptive}
    \addplot [forget plot,only marks,black,mark indices={8,9,10,11,35,36,37,38},mark=triangle] table[header=true,x expr = \thisrowno{2},y index = {4}, col sep = comma]{Data/ArcLengthRoof/RoofAdaptive_g9C9.csv};
    \addplot[ultra thin ,gray,no markers,domain=0:2e7,samples=2,name path=A] {1e-8};
    \addplot[ultra thin ,gray,solid,no markers,domain=0:2e7,samples=2,name path=B] {1e-10};
    \addplot[gray!50,opacity=20] fill between [of=A and B];

    \nextgroupplot[ylabel={\# DoFs},ymin=1e2, ymax = 2e4,xlabel={$\Vert\vb{u}\Vert$},]
    \addplot+[fill opacity=1.0,no markers] table[header=true,x expr = \thisrowno{2},y index = {3}, col sep = comma]{Data/ArcLengthRoof/RoofUniform_r4_g9C9.csv};
    \addplot+[fill opacity=1.0,no markers] table[header=true,x expr = \thisrowno{2},y index = {3}, col sep = comma]{Data/ArcLengthRoof/RoofUniform_r5_g9C9.csv};
    \addplot+[fill opacity=1.0] table[header=true,x expr = \thisrowno{2},y index = {3}, col sep = comma]{Data/ArcLengthRoof/RoofAdaptive_g9C9.csv};
    \addplot [forget plot,only marks,black,mark indices={8,9,10},mark=triangle,
                            nodes near coords={\coordindex},
                            every node near coord/.append style={anchor=south east,font=\tiny},
                        ] table[header=true,x expr = \thisrowno{2},y index = {3}, col sep = comma]{Data/ArcLengthRoof/RoofAdaptive_g9C9.csv};
    \addplot [forget plot,only marks,black,mark indices={11,35,36,37,38},mark=triangle,
                            nodes near coords={\coordindex},
                            every node near coord/.append style={anchor=south,font=\tiny},
                        ] table[header=true,x expr = \thisrowno{2},y index = {3}, col sep = comma]{Data/ArcLengthRoof/RoofAdaptive_g9C9.csv};
    \end{groupplot}
    \end{tikzpicture}
\label{fig:roof_results2}
\caption{Errors (top) and number of degrees of freedom (DoFs, bottom) for the goal functional $\mathbcal{L}=\int_\Omega \Vert \vb{m}(\vb{u})\Vert\dd{\Omega}$ with tolerances $\text{tol}_r=10^{-10}$ and $\text{tol}_c=10^{-8}$ for the collapsing roof subject to the displacement norm $\Vert\vb{u}\Vert$. The gray region represents the region $\Delta\mathbcal{L}\in[\text{tol}_r,\text{tol}_c]$ where refinement and coarsening are performed and where the adaptivity iterations are terminated. Above the gray region only refinement is performed and below the gray region only coarsening is performed.}
\label{fig:roof_results}
\end{figure}

The geometry for the benchmark problem is presented in \cref{fig:setupRoof}. Here, the radius of the roof $R=2.540[\text{m}]$, the angle is $\theta=0.1\:[\text{rad}]$ and the length and thickness are, respectively, $L=0.508\:[\text{m}]$ and $t=6.35\cdot 10^3\:[\text{mm}]$. Moreover, the material properties are $E=3102\:[\text{MPa}]$ and $\nu=0.3\:[\text{-}]$ for a Saint-Venant Kirchhoff material. Only a quarter of the roof is modelled, as depicted in \cref{fig:setupRoof}, because of symmetry. The simulation is performed using a Crisfield arc-length method \cite{Crisfield1981} with arc-length $\Delta L=25$ and a zero force-scaling. The goal functional that is used for error evaluation and adaptivity is based on the norm of the flexural strain tensor over the whole domain $\mathbcal{L}=\int_\Omega \Vert \vb{m}(\vb{u})\Vert\dd{\Omega}$. The jump parameter for admissible meshing is set to $m=2$. The mesh will be refined when $\Delta \mathbcal{L}>\text{tol}_r$ and coarsened when $\Delta \mathbcal{L}<\text{tol}_c$. As discussed in \cref{sec:algorithm}, $\text{tol}_r<\text{tol}_c$ such that the mesh is refined and coarsened simultaneously when $\Delta\mathbcal{L}\in[\text{tol}_r,\text{tol}_c]$, which is also the condition for termination. The maximum number of mesh adaptivity iterations is set to 5 in this case. The tolerances $\text{tol}_r$ and $\text{tol}_c$ are determined based on the results of uniformly refined simulations with $16\times16$ (918 DoFs) and $32\times32$ elements (3366 DoFs) by taking a wide band around the error envelopes in \cref{fig:roof_results2} excluding peaks. The tolerances are $(\text{tol}_r,\text{tol}_c)=(10^{-10},10^{-8})$. It should be noted that these tolerances can also be based on requirements in engineering, or they can be determined during the computations; both are beyond the scope of this paper. The refinement parameter is set to $\rho_r=0.5$,  the coarsening parameter to $\rho_c=0.05$ and the maximum refinement level is $11$. The adaptive and uniform meshes are modelled with bi-cubic B-spline basis functions (i.e. $p=3$).\\

\begin{figure}
\centering
\begin{tikzpicture}
    \begin{axis}[
        width=0.6\linewidth,
        height=0.3\textheight,
        tickpos=left,
        grid = major,
        xlabel = {$w_A$},
        ylabel = {$\lambda P\:(\times 10^3)$},
        zlabel = {$\Vert\vb{u}\Vert$},
        restrict x to domain=0:50,
        restrict y to domain=-1:1,
        restrict z to domain=0:2e7,
        xmin = 0, xmax = 50,
        ymin = -1, ymax = 1,
        zmin = 0, zmax = 2e7,
        view = {20}{40},
        legend pos=outer north east,
    ]
     \addplot3[black!50,mark=x,only marks,fill opacity=1.0,on layer=axis foreground] table[header=true,x expr = \thisrowno{0},y expr= \thisrowno{1}/1000,z expr = {0}, col sep = comma]{Data/ArcLengthRoof/RoofRef.csv};\addlegendentry{Reference \cite{Sze2004}}

    \pgfplotsset{cycle list shift=1}
    \addplot3+[solid,fill opacity=1.0,on layer=axis foreground,mark=none] table[header=true,x expr = -\thisrowno{0},y expr = {1},z index = {2}, col sep = comma]{Data/ArcLengthRoof/RoofAdaptive_g9C9.csv};\addlegendentry{Adaptive}
    \pgfplotsset{cycle list shift=0}
    \addplot3+[solid,fill opacity=1.0,on layer=axis foreground,mark=none,forget plot] table[header=true,x expr = {0},y expr= \thisrowno{1}/1000,z index = {2}, col sep = comma]{Data/ArcLengthRoof/RoofAdaptive_g9C9.csv};

    \addplot3+[solid,fill opacity=1.0,on layer=axis foreground,mark=none,forget plot] table[header=true,x expr = -\thisrowno{0},y expr= \thisrowno{1}/1000,z expr = {0}, col sep = comma]{Data/ArcLengthRoof/RoofAdaptive_g9C9.csv};

    \addlegendimage{black,only marks,mark=triangle}\addlegendentry{Mesh points}
    \addplot3 [solid,fill opacity=1.0,on layer=axis foreground,only marks,black,mark indices={8,9,10,11,35,36,37,38},mark=triangle] table[header=true,x expr = -\thisrowno{0},y expr = {1},z index = {2}, col sep = comma]{Data/ArcLengthRoof/RoofAdaptive_g9C9.csv};
    \addplot3 [solid,fill opacity=1.0,on layer=axis foreground,only marks,black,mark indices={8,9,10,11,35,36,37,38},mark=triangle] table[header=true,x expr = {0},y expr= \thisrowno{1}/1000,z index = {2}, col sep = comma]{Data/ArcLengthRoof/RoofAdaptive_g9C9.csv};
    \addplot3 [solid,fill opacity=1.0,on layer=axis foreground,only marks,black,mark indices={8,9,10,11,35,36,37,38},mark=triangle] table[header=true,x expr = -\thisrowno{0},y expr= \thisrowno{1}/1000,z expr = {0}, col sep = comma]{Data/ArcLengthRoof/RoofAdaptive_g9C9.csv};

    \addplot3+[black,fill opacity=0.5,solid,mark=none,forget plot] table[header=true,x expr = -\thisrowno{0},y expr= \thisrowno{1}/1000,z index = {2}, col sep = comma]{Data/ArcLengthRoof/RoofAdaptive_g9C9.csv};
    \end{axis}
    \end{tikzpicture}
\caption{Projection the adaptively refined result from the commonly used $\lambda P,w_A$-space \cite{Sze2004} onto the displacement norm $\Vert\vb{u}\Vert$. The solid lines correspond to the results obtained by the adaptively refined mesh, the triangular markers correspond to points of which the mesh is provided in \cref{fig:roof_meshes} and the cross-markers indicate the result from \cite{Sze2004}.}
\label{fig:roof_coords}
\end{figure}
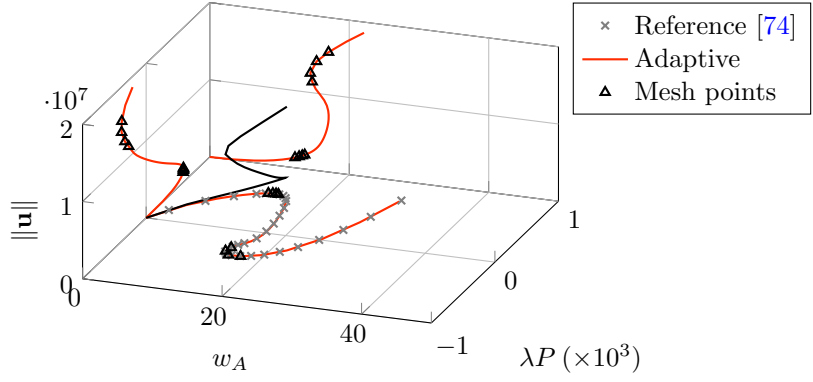

In \cref{fig:roof_coords}, the results of a simulation of the collapsing roof for the adaptive mesh are plotted in a $\Vert{\vb{u}}\Vert$, $\lambda P$, $w_A$-space. Reference solutions for the present benchmark problem are typically given in the $\lambda P$, $w_A$-space, but since the solution curve is not bijective an alternative coordinate $\Vert{\vb{u}}\Vert$ is used to represent solutions for this benchmark problem. This is motivated by the projection of the solutions $\lambda P $ and $w_A$ projected against $\Vert\vb{u}\Vert$ in \cref{fig:roof_coords}. The results of \cite{Sze2004} are provided as a reference.\\

\begin{figure}
 \centering
   \begin{subfigure}{\linewidth}
  \centering
  \input{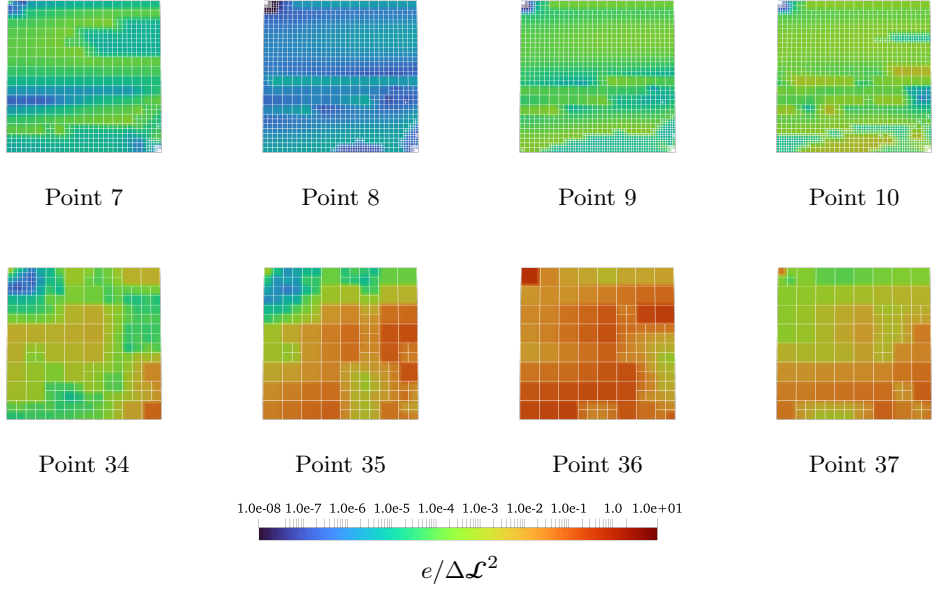}
  \label{}
   \end{subfigure}
\caption{Meshes corresponding to the points marked in \cref{fig:roof_coords,fig:roof_results} by the black-bordered marks. The point $A$ marks the point where the load is applied in \cref{fig:setupRoof}. The top row of meshes corresponds to the first limit point in \cref{fig:roof_coords} and the bottom row corresponds to the second limit points. The meshes are ordered from left to right for increasing solution steps. The elements are coloured according to the squared error $e_k$ normalised by the total error $\Delta\mathbcal{L}$, i.e. $e_k/\Delta\mathbcal{L}^2$.}
\label{fig:roof_meshes}
\end{figure}

In \cref{fig:roof_results}, the error and the number of degrees of freedom are plotted against $\Vert\vb{u}\Vert$. The error envelopes for the uniform meshes show that a large peak occurs around $\Vert\vb{u}\Vert=0.2$, relating to the first limit point of the collapse of the roof, as seen by the markers in \cref{fig:roof_coords}. Additionally, it can be seen that the present algorithm providing mesh adaptivity manages to keep the error within specified bounds (see \cref{fig:roof_results}, top), except on the peak just before $\Vert\vb{u}\Vert$ where the maximum number of adaptivity iterations is insufficient. Furthermore, it has consistently less degrees of freedom than the uniform mesh with $32\times32$ elements. In \cref{fig:roof_meshes} a selection of meshes is provided. The meshes are provided as series of 4 consecutive meshes around the limit points of the solution curve, as indicated in \cref{fig:roof_coords}, and are depicted in increasing order from left to right for the first (top) and second (bottom) limit points. The black-bordered markers in \cref{fig:roof_coords,fig:roof_results} indicate the points of which the meshes are shown. From the first row of meshes in \cref{fig:roof_meshes}, it can be seen that the first limit point requires relatively fine meshes and that the elements start concentrated around point $A$ and its diagonal opposite and spread out on the bottom symmetry boundary as the snapping takes place. Furthermore, in the bottom row of \cref{fig:roof_meshes}, it can be seen that the second limit point does not require many elements, hence the number of elements is slowly decreasing throughout this section of the load-displacement curve. For a complete overview of the mesh in each load-step, we refer to Video 1 in the supplementary material of this paper.\\

Concluding, this example shows that the goal-adaptive meshing procedure is capable of keeping the error in terms of a goal functions within pre-defined bounds for a solution stepping simulation with limit-point instabilities. Throughout the simulation, the procedure keeps a relatively high efficiency per degree of freedom compared to uniform meshes. It should be noted, however, that the adaptive refinement iterations require a higher computational demand. Therefore, the next and final example provides a procedure where no adaptivity iterations are performed.\\

\newpage

\subsection{Wrinkling analysis}\label{subsec:wrinkling}
As a last example, the wrinkling analysis of a thin membrane subject to a tensile load is considered. This problem is inspired by \cite{Cerda2002,Cerda2003,Panaitescu2019} and was previously modelled using isogeometric Kirchhoff-Love shells in \cite{Verhelst2021}, to which we refer for a detailed problem set-up. The goal of this benchmark in the present paper is to demonstrate the use of the goal-adaptive meshing procedure on a bifurcation problem with geometric and material non-linearities.\\

\begin{figure}
\centering
 \begin{minipage}{0.55\linewidth}
 \resizebox{\linewidth}{!}
 {
\begin{tikzpicture}[scale=0.7]
\footnotesize
\def\B{4}
\def\L{6}
\def\strain{2};
\def\overlap{0.02}

\fill[opacity=0.15, bottom color=black!10, top color=black!05, draw=black] (0,0) -- (\L,0) -- (\L,\B) -- (0,\B) --cycle;

\filldraw[draw=black,bottom color=black!20, top color=black!10] (0,0) to[out=15,in=180] (\L/2+\strain/2,0.5) to[out=0,in=165] (\L+\strain,0) -- (\L+\strain,\B) to[out=205,in=0] (\L/2+\strain/2,\B-0.5) to[out=180,in=-15] (0,\B) -- cycle;

\filldraw[dashed,draw=black,bottom color=black!30, top color=black!20,thick] (\L/2+\strain/2,\B/2)  -- (\L+\strain,\B/2) -- (\L+\strain,\B) to[out=205,in=0] (\L/2+\strain/2,\B-0.5) -- cycle;

\draw[latex-latex] (0,\B+0.1)--node [above]{$(1+\epsilon)L$}(\L+\strain,\B+0.1);
\draw[latex-latex] (0,-0.1)--node [below]{$L$}(\L,-0.1);
\draw[latex-latex] (0.1,0)--node [right]{$W$}(0.1,\B);

\node at (3*\L/4+3*\strain/4,\B)[rotate=0,below]{$\Gamma_1$};
\node at (\L+\strain,3*\B/4)[rotate=0,left]{$\Gamma_2$};
\node at (3*\L/4+3*\strain/4,\B/2)[rotate=0,above]{$\Gamma_3$};
\node at (\L/2+\strain/2,3*\B/4)[rotate=0,right]{$\Gamma_4$};

\draw[-latex,thick] (\L+\strain,\B/2) -- (\L+\strain+1,\B/2) node[right]{$P$};

\filldraw[bottom color=black!80,top color=black!70] (\overlap,-\overlap) -- (\overlap,\B+\overlap) -- (-0.2,\B+\overlap) -- (-0.2,-\overlap)--cycle;
\filldraw[bottom color=black!80,top color=black!70] (\L-\overlap+\strain,-\overlap) -- (\L-\overlap+\strain,\B+\overlap) -- (\L+0.2+\strain,\B+\overlap) -- (\L+0.2+\strain,-\overlap)--cycle;

\fill[color=black,opacity=0.1](0.5*\L+0.5*\strain,\B/2) ellipse [x radius=\L/2,y radius=0.1];
\fill[color=black,opacity=0.2](0.5*\L+0.5*\strain,\B/2+0.25) ellipse [x radius=\L/2-0.2,y radius=0.1];
\fill[color=black,opacity=0.2](0.5*\L+0.5*\strain,\B/2-0.25) ellipse [x radius=\L/2-0.2,y radius=0.1];
\fill[color=black,opacity=0.05](0.5*\L+0.5*\strain,\B/2+0.5) ellipse [x radius=\L/2-0.5,y radius=0.1];
\fill[color=black,opacity=0.05](0.5*\L+0.5*\strain,\B/2-0.5) ellipse [x radius=\L/2-0.5,y radius=0.1];
\fill[color=black,opacity=0.15](0.5*\L+0.5*\strain,\B/2+0.75) ellipse [x radius=\L/2-0.9,y radius=0.1];
\fill[color=black,opacity=0.15](0.5*\L+0.5*\strain,\B/2-0.75) ellipse [x radius=\L/2-0.9,y radius=0.1];

\draw[-latex] (0.5,0.5) -- node[right,inner sep=8pt] {$x$}(1.0,0.5);
\draw[-latex] (0.5,0.5) -- node[above,inner sep=8pt] {$y$}(0.5,1.0);

\end{tikzpicture}
 }
 \end{minipage}
 \hfill
 \begin{minipage}{0.35\linewidth}
 \begin{tabular}{llll}
 $c_1$    &$=$& $3.13\cdot10^5$ & $[\text{Pa}]$ \\
 $c_2$ &$=$& $1.24\cdot10^5$ & $[\text{Pa}]$ \\
 $L$    &$=$& $280$ & $[\text{mm}]$\\
 $W$    &$=$& $140$ & $[\text{mm}]$\\
 $t$     &$=$& $0.14$ & $[\text{mm}]$\\
 $P$    & $=$& \multicolumn{2}{l}{Variable}
\end{tabular}
 \end{minipage}
 \caption{Geometry and parameters for the wrinkling problem. The boundary $\Gamma_1$ is free, the boundary $\Gamma_2$ is fixed in $x$ and $y$ direction and rotations around the $y$-axis are fixed, the boundaries $\Gamma_3$ and $\Gamma_4$ have symmetry conditions applied.}
\label{fig:wrinkling_setup}
\end{figure}
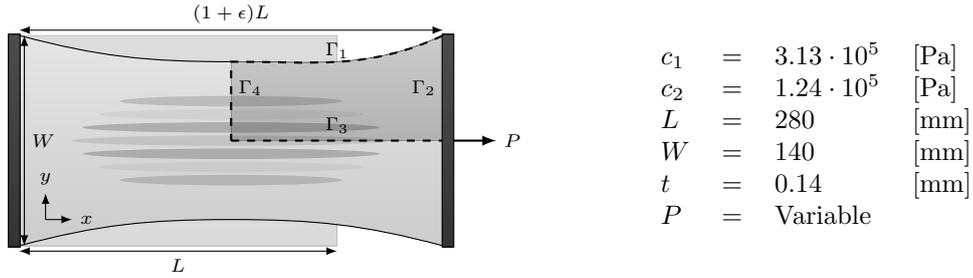

Given the geometry in \cref{fig:wrinkling_setup}, a quarter of the domain is considered with a symmetry boundary condition on $\Gamma_4$ and an antisymmetry condition on $\Gamma_3$. Furthermore, the boundary $\Gamma_1$ is free and on $\Gamma_2$ the $x$-displacement is constant and a horizontal load is applied on this side. The sheet has dimensions $L=280\:[\text{mm}]$, $W=140\:[\text{mm}]$ and $t=0.14\:[\text{mm}]$ such that $L/W=2$ and $t/W=10^3$. Furthermore, the material is modelled using a Mooney-Rivlin material model with strain energy density function
\begin{equation}
 \Psi(\vb{C}) = \frac{c_1}{2}\qty(I_1-3) + \frac{c_2}{2}\qty(I_2-3)
\end{equation}
and with parameters $c_1=3.16\cdot 10^5\:[\text{Pa}]$ and $c_1=1.24\cdot 10^5\:[\text{Pa}]$. Reference solutions are given by \cite{Verhelst2021} for isogeometric Kirchhoff-Love shell analysis and using ANSYS and LS-DYNA FEA models, respectively. Furthermore, \cite{Panaitescu2019} provides experimental data of the maximum amplitude with respect to the strain of the sheet. In the present paper, the reference simulations are performed on uniform cubic meshes with $32\times32$ and $64\times64$ elements, respectively.\\

For the adaptive simulation, a THB spline mesh with initially $32\times32$ elements is used and mesh adaptivity is activated after wrinkling initiation since the errors in the pre-wrinkling regime are small due to of the lack of out-of-plane deformations of the sheet. The goal-functional is a displacement-based functional on the $z$-component, i.e. $\mathbcal{L}(\vb{u})=\int_\Omega\vb{u}\cdot\vb{e}_z\dd{\Omega}$. The tolerances for refinement and coarsening are $\text{tol}_r=10^{-14}$ and $\text{tol}_c=10^{-10}$, respectively, and they are chosen based on the error envelope of the uniform refinement. The adaptive meshing parameters are chosen as $(\rho_r,\rho_c)=(0.5,0.005)$. These parameters are chosen based on the behaviour of the global error in the first load-steps afer bifurcation. Contrary to the previous example in \cref{subsec:roof}, there are no refinement iterations performed within the load step. This means that the refinement  and coarsening operations are performed after the load step based on the magnitude of the error compared to the tolerances. Furthermore, when the error is below a tolerance $\rho_{c,\text{min}}$, the initial mesh is used again. This is done to prevent further coarsening after re-stabilisation of the wrinkles, i.e. the moment when the wrinkles have dissapeared. For the refinement algorithm, a maximum depth of the THB grid is fixed to $11$ levels and the jump parameter is set to $m=2$. The wrinkling simulation is performed using a Crisfield arc-length method \cite{Crisfield1981} with a quadratic procedure to compute the mode shape at the bifurcation \cite{Wriggers1988}. This procedure further described in \cite{Verhelst2019}.\\

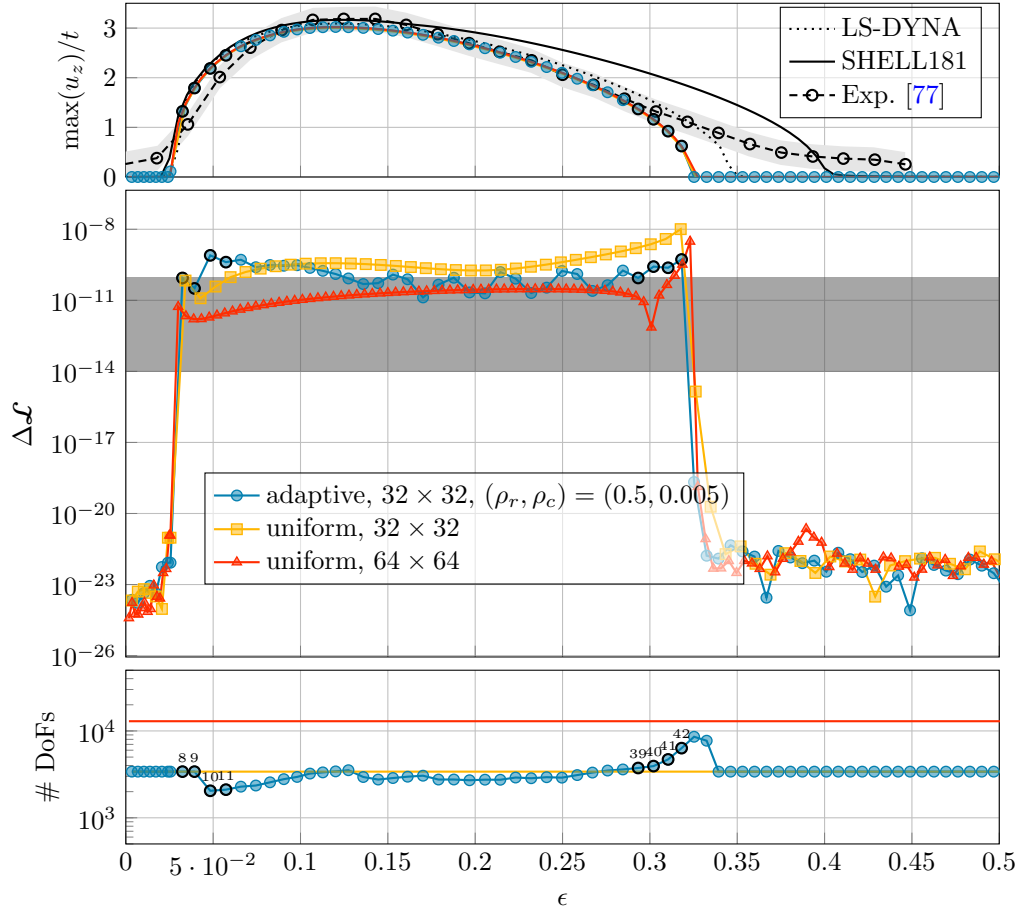
\begin{figure}
\centering
\begin{tikzpicture}
    \begin{groupplot}[
        group style={
            group name=my plots,
            group size=1 by 3,
            xlabels at=edge bottom,
            xticklabels at=edge bottom,
            ylabels at=edge left,
            yticklabels at=edge left,
            vertical sep=5pt,
            horizontal sep=5pt
        },
        width=\linewidth,
        height=0.2\textheight,
        tickpos=left,
        grid = major,
        xlabel=$\epsilon$,
        xmin = 0,
        xmax = 0.5,
    ]
    \nextgroupplot[
        ylabel=$\max(u_z)/t$,
        ymin = 0,
        ymax = 3.5,
    ]
    \addplot+[no marks,black,dotted,mark=diamond ,mark size=1] table[header=true,x index = {0},y index = {1}, col sep = comma]{Data/ArcLengthWrinkling/LSDynaFull.csv};\addlegendentry{LS-DYNA}
    \addplot+[no marks,black,solid,mark=+,mark size=1] table[header=true,x index = {0},y index = {1}, col sep = comma]{Data/ArcLengthWrinkling/ANSYS.csv};\addlegendentry{SHELL181}
    \addplot+[mark=o,mark size=2,black,mark options={black,solid}] table[header=true,x index = {0},y index = {1}, col sep = comma]{Data/ArcLengthWrinkling/PanaitescuExperiment.csv};\addlegendentry{Exp. \cite{Panaitescu2019}}
    \addplot [name path=upper,draw=none,mark=none,forget plot] table[header=true,x index = {2},y index = {3}, col sep = comma]{Data/ArcLengthWrinkling/PanaitescuExperiment.csv};
    \addplot [name path=lower,draw=none,mark=none,forget plot] table[header=true,x index = {4},y index = {5}, col sep = comma]{Data/ArcLengthWrinkling/PanaitescuExperiment.csv};
    \addplot [fill=black!10,forget plot] fill between[of=upper and lower];

    \pgfplotsset{cycle list shift=-3}
    \addplot+[] table[header=true,x index = {1},y index = {2}, col sep = comma]{Data/ArcLengthWrinkling/Wrinkling_r5_adaptive_R050_C0005.csv};
    \addplot [forget plot,only marks,black,mark indices={9,10,11,12,40,41,42,43},mark=o] table[header=true,x index = {1},y index = {2}, col sep = comma]{Data/ArcLengthWrinkling/Wrinkling_r5_adaptive_R050_C0005.csv};
    \addplot+[no marks,solid] table[header=true,x index = {1},y index = {2}, col sep = comma]{Data/ArcLengthWrinkling/Wrinkling_r5.csv};
    \addplot+[no marks,solid] table[header=true,x index = {1},y index = {2}, col sep = comma]{Data/ArcLengthWrinkling/Wrinkling_r6.csv};

    \nextgroupplot[ylabel={$\Delta\mathbcal{L}$},ymode=log,height=0.4\textheight,legend style={at={(axis cs: 0.045,1e-23)},anchor=south west}]
    \fill[black!70,opacity=0.5] (axis cs:0,1e-14)-- (axis cs:0.5,1e-14)--(axis cs:0.5,1e-10)--(axis cs:0,1e-10)--cycle;


    \addplot+[] table[header=true,x index = {1},y index = {4}, col sep = comma]{Data/ArcLengthWrinkling/Wrinkling_r5_adaptive_R050_C0005.csv};\addlegendentry{adaptive, $32\times32$, $(\rho_r,\rho_c)=(0.5,0.005)$}
    \addplot [forget plot,only marks,black,mark indices={9,10,11,12,40,41,42,43},mark=o] table[header=true,x index = {1},y index = {4}, col sep = comma]{Data/ArcLengthWrinkling/Wrinkling_r5_adaptive_R050_C0005.csv};
    \addplot+[,solid] table[header=true,x index = {1},y index = {4}, col sep = comma]{Data/ArcLengthWrinkling/Wrinkling_r5.csv};\addlegendentry{uniform, $32\times 32$}
    \addplot+[,solid] table[header=true,x index = {1},y index = {4}, col sep = comma]{Data/ArcLengthWrinkling/Wrinkling_r6.csv};\addlegendentry{uniform, $64\times 64$}

%
%
%
%
    \nextgroupplot[ylabel={\# DoFs},ymode=log,ymin=5e2, ymax=5e4]
    \addplot+[] table[header=true,x index = {1},y index = {3}, col sep = comma]{Data/ArcLengthWrinkling/Wrinkling_r5_adaptive_R050_C0005.csv};
    \addplot [forget plot,only marks,black,mark indices={9,10,11,12,40,41,42,43},mark=o,
                nodes near coords={\coordindex},
                every node near coord/.append style={anchor=south,font=\tiny},
                ] table[header=true,x index = {1},y index = {3}, col sep = comma]{Data/ArcLengthWrinkling/Wrinkling_r5_adaptive_R050_C0005.csv};
    \addplot+[no markers,solid] table[header=true,x index = {1},y index = {3}, col sep = comma]{Data/ArcLengthWrinkling/Wrinkling_r5.csv};
    \addplot+[no markers,solid] table[header=true,x index = {1},y index = {3}, col sep = comma]{Data/ArcLengthWrinkling/Wrinkling_r6.csv};
    \end{groupplot}
    \end{tikzpicture}
\caption{The non-dimensional maximal out-of-plane deformation $\max(u_z)/t$ (top), the goal functional error $\Delta\mathbcal{L}$ (mid) and the number of degrees of freedom of the computational mesh (bottom) with respect to the strain of the sheet $\epsilon$. The markers with a black border represent the points of which the meshes are provided in \cref{fig:wrinkling_results_meshes}. The colored lines are the solutions obtained by the present model with a Mooney-Rivlin material model with uniform or adaptive meshes. The solid line in the top figure is a SHELL181 result obtained using ANSYS, the dotted line in the top figure is a result obtained using the fully integrated shell in LS-DYNA and the dashed line with markers represents experimental data obtained by \cite{Panaitescu2019}.}
\label{fig:wrinkling_results}
\end{figure}

In \cref{fig:wrinkling_results}, the results for an adaptive isogeometric wrinkling simulation are provided. The top figure provides the normalised wrinkling amplitude with respect to the strain of the sheet ($\epsilon$), compared to the IGA and ANSYS SHELL181 element reference results from \cite{Verhelst2021} as well as the experimental results from \cite{Panaitescu2019}. The figure in the middle provides the error estimate in terms of the goal functional $\mathbcal{L}$ with respect to the strain of the sheet and the bottom figure provides the number of degrees of freedom with respect to the strain of the sheet. Firstly, the results from both uniform meshes show that the error estimate $\Delta\mathbcal{L}$ is close to zero when the wrinkles initiate, since the sheet is perfectly flat. As soon as the wrinkles form, the error estimate becomes non-zero and it peaks at the moment of re-stabilisation (i.e. the moment when the amplitude vanishes). After re-stabilisation, the error estimate is low, but slightly higher than before wrinkling, probably because the sheet is not numerically flat.\\

The adaptive meshing simulations show that even with zero inner-iterations for mesh adaptivity the adaptive mesh provides accurate results for a relatively small number of degrees of freedom. Just after wrinkling initiates, the adaptive meshing error peaks due to the coarsening of a large number of elements (as can be seen by the drop of degrees of freedom in the bottom figure). However, the mesh adaptively refines until in the gray region in the middle figure, where combined refinement and coarsening imply that the error balances around $\text{tol}_c$, i.e. the upper bound of the marked region. Towards the end of the wrinkling phase, the error increases for the uniformly refined mesh, explaining the increase in number of degrees of freedom for the adaptive mesh. Nevertheless, the adaptive mesh provides more accurate results compared to the $64\times64$ uniform mesh with less degrees of freedom.\\

\begin{figure}
\centering
\input{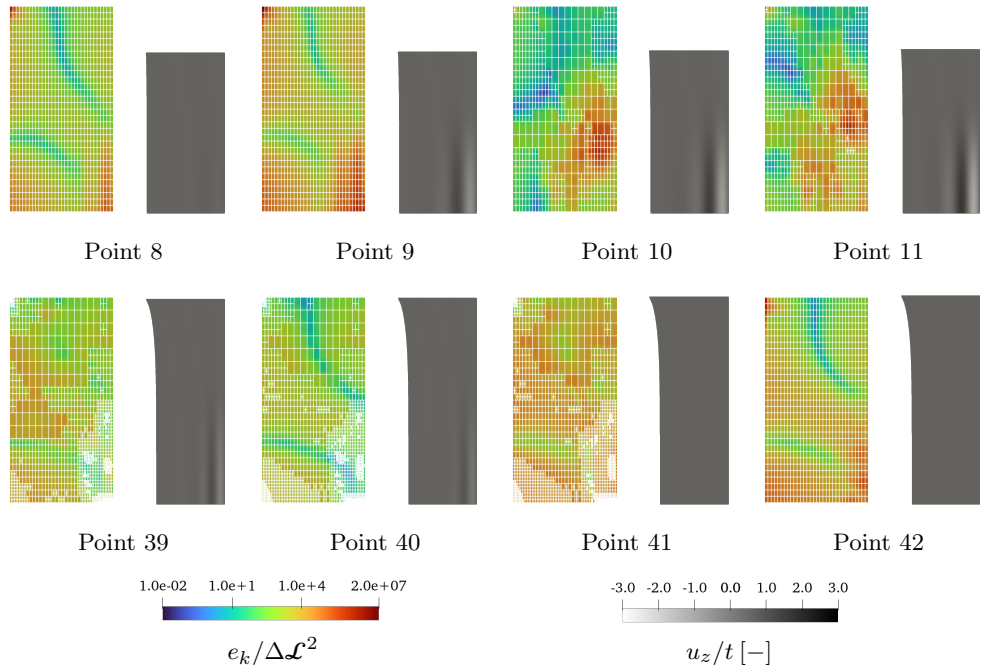}
\caption{Normalised element errors (left) and the normalised out-of-plane deformation (right) for the wrinkling benchmark plotted on the undeformed geometry with the corresponding elements.}
\label{fig:wrinkling_results_meshes}
\end{figure}

In \cref{fig:wrinkling_results_meshes} a selection of meshes from both adaptive simulations are provided, specifically for the wrinkling initiation and re-stabilisation points. For a complete overview of the mesh in each load-step, we refer to Video 2 in the supplementary material of this paper. The evolution of the meshes shows that the mesh elements concentrate around the wrinkles (bottom-right) and in the top-left corner, which represents the corner between the clamped edge (top) and the free edge. The former is expected since the employed goal functional is based on out-of-plane deformations. However, the fact that the mesh concentrates around the top-left corner is non-intuitive, but tells that this area is important to reduce the global error in terms of the out-of-plane deformations. Furthermore, it can be seen that around the re-stabilisation of the wrinkles the mesh concentrates around the bottom-left corner, indicating that this corner is of importance in accurately modelling the wrinkling amplitudes in the whole domain.\\

Concluding, the wrinkling benchmark shows the potential of mesh adaptivity for such applications. With fewer degrees of freedom, the THB-spline mesh is able to approximate the solution around a pre-defined error, eventhough the selection of the refinement and coarsening parameters and the tolerances has not been optimised in this study.

\clearpage

\section{Conclusions}\label{sec:conclusions}
This paper presents an adaptive method for isogeometric Kirchhoff-Love shells. The main contributions of this paper are a goal-adaptive error estimator for isogeometric Kirchhoff-Love shells using the Dual-Weighted Residual method and a slightly modified suitably graded refinement scheme taking into account refined elements in the definition of the coarsening neighborhood.\\

Using the Dual-Weighted Residual method and given a pre-defined goal functional (e.g the second-principal stress integrated over the domain), an estimator for the error in terms of this goal functional can be defined. The adjoint problem that needs to be solved on the original mesh and on a nested degree-elevated (`enriched') mesh has been defined for the isogeometric Kirchhoff-Love shell. In addition, the operators for modal and linear buckling analysis have been derived, implying an additional generalised eigenvalue problem to be solved on the enriched mesh. For suitable grading, the works of \cite{Buffa2016,Buffa2022,Bracco2018a,Carraturo2019} have been closely followed. In order to be able to refine and coarsen in the same iteration, refined elements have been added to the original definition of the coarsening neighborhood.\\

To assess the proposed adaptive isogeometric method for Kirchhoff-Love shells, few numerical benchmark problems have been evaluated. Linear static analysis provided an analytical solution has been used to evaluate the DWR error estimators. The eigenvalue problems for modal and buckling analyses have been evaluated on respectively the problem of circular plate vibration and square plate buckling. Based on the linear, modal and buckling analysis with analytical solutions, it can be concluded that the DWR estimator for Kirchhoff-Love shells can be used with several goal functionals and in several applications, as it provides high accuracy with respect to the exact errors.\\

Using the problem of a pinched membrane, the error estimator has been used to adaptively refine a mesh in a non-linear example. From this example, it can be concluded that eventhough the challenging non-linear problem, high accuracy per DoF can be obtained compared to uniformly refined meshes based on different goal functionals. Lastly, the adaptive isogeometric method from the present paper has been evaluated in solution-stepping problems for structural instabilities. Firstly, the method was applied to limit-point instability problem of a collapse of a cylindrical roof. Here, inner adaptivity iterations were performed for each load-step until the error was located in a desired interval. Again, the present method provides a high accuracy per degree of freedom. In addition, it was shown that the method is indeed able to provide adaptive meshes with respect to a pre-defined interval for a given goal functional. The last benchmark problem involves a tension-wrinkling bifurcation instability of a thin membrane. In this benchmark problem, adaptive meshing has been applied in the post-buckling regime based on wrinkling amplitudes. Here, no adaptivity iterations within the load steps were performed, showing that the method is still able to provide good adaptivity with respect to the pre-defined tolerances. Also, the results of the wrinkling error estimators show that the error peaks in the re-stabilisation phase, where the highest deviation with experimental results are observed.\\

Future developments for the present method include the application on multi-patch domains, both coupled with penalty methods as well as with globally continuous bases as presented in \cite{Farahat2023a} to handle more complex geometries. Furthermore, structural dynamics have been left out of the scope of this paper, since the DWR for dynamic problems requires backwards-in-time evaluation of the adjoint problem, which is ideally combined with parallel-in-time methods like ParaReal or MGRIT \cite{Falgout2014}. Lastly, future work can be done on the (adaptive) determination of the adaptive meshing parameters. On the one hand, one can apply the present method on real-world engineering applications, taking realistic goal functionals and margins into the adaptivity algorithm. On the other hand, advanced schemes for triggering pure or combined refinement or coarsening together with their parameters can be futher investigated.

\backmatter

\bmhead{Supplementary information}

\bmhead{Video 1}: Mesh evolution of for the collapsing roof. 
\bmhead{Video 2}: mesh evolution for the wrinkling example.

\bmhead{Declarations}

\bmhead{Funding}
H.M. Verhelst, M. Möller and J.H. Den Besten are greatful to the faculties of Electrical Engineering, Mathematics and Computer Science (EEMCS) and of Mechanical, Materials and Maritime Engineering (3mE) of Delft University of Technology for the financial support to conduct this research. A. Mantzaflaris is greatful for the financial of the EU ITN network GRAPES.

\bmhead{Conflict of interest} The authors declare that they have no known competing financial interests or personal relationships that could have appeared to influence the work reported in this paper.

\bmhead{Author contributions} The authors contributed to this paper in the following way. \textbf{H.M. Verhelst}: Conceptualization, Formal analysis, Investigation, Methodology, Software, Validation, Visualization, Writing -- original draft, Writing -- review \& editing. \textbf{A. Mantzaflaris}: Software, Writing -- review \& editing. \textbf{M. Möller}:  Conceptualization, Software, Writing -- review \& editing, Funding acquisition, Project administration, Supervision. \textbf{J.H. Den Besten}:  Conceptualization, Writing – review \& editing, Funding acquisition, Project administration, Supervision.


\bibliography{references}

\end{document}